\documentclass[amscd,amssymb,verbatim,12pt]{amsart}
\usepackage{epsfig}
\usepackage{graphicx}
\usepackage{pdfpages}
\usepackage{amssymb,amsmath}

\usepackage[all]{xy}

\theoremstyle{plain}

\newtheorem{thm}{Theorem}[section]
\newtheorem{lemma}[thm]{Lemma}

\newtheorem{prop}[thm]{Proposition}

\newtheorem*{theorem.nonumber}{Theorem}

\theoremstyle{remark}
\newtheorem{remark}[thm]{Remark}
\newtheorem{definition}[thm]{Definition}
\newtheorem{example}[thm]{Example}

\newtheorem{question}[thm]{Question}

\def \={\ = \ }
\def \+{\ +\ }
\def \-{\ - \ }

\def \b|{\big |}

\def \g1{\Gamma_1}
\def \<{\langle}
\def \>{\rangle}

\def \R{\Bbb R}

\def \Z{\Bbb Z}
\def \CC{\mathcal C}
\def \D{\mathcal D}

\def \Hom{\text{Hom}\,}

\def \Aut{\text{Aut\,}}

\def \germ{\text{germ\,}}
\def \Sph{\text{Sph\,}}

\begin{document}
\title
{The rates of growth in an acylindrically hyperbolic group}

\author
{Koji Fujiwara }
\email{kfujiwara@math.kyoto-u.ac.jp}
\address{Department of Mathematics, Kyoto University,
Kyoto, 606-8502, Japan}

\thanks{ The author is
    supported in part by Grant-in-Aid for Scientific Research
    (No. 15H05739, 20H00114).  }

%\subjclass{53C23}

\begin{abstract} 
Let $G$ be an acylindrically hyperbolic group
 on a $\delta$-hyperbolic
space $X$. Assume there exists $M$ such that for any finite generating set 
$S$ of $G$, the set $S^M$ contains a hyperbolic element on $X$. Suppose that 
$G$ is equationally Noetherian. Then we show 
the set of the growth rates of $G$ is well-ordered (Theorem \ref{main}).
The conclusion was known for hyperbolic groups, and this is a generalization. 

Our result applies to all lattices in simple Lie groups of rank-1 (Theorem \ref{main.ex}), and more generally,
some family of relatively hyperbolic groups (Theorem \ref{thm.rel.hyp}).
It also applies to the fundamental group, of exponential growth, 
of a closed orientable $3$-manifold except for the case that the manifold has Sol-geometry (Theorem \ref{3manifold}).
A potential application is a mapping class group, to
which the theorem applies if it is equationally Noetherian.

\end{abstract}
\maketitle

\section{Introduction}
\subsection{Definitions and results}

Let $G$ be a finitely generated group with a finite generating set $S$. 
We always assume that $S=S^{-1}$. 
Let $B_n(G,S)$ be the set of elements in $G$ whose word lengths are
at most $n$ with respect to the generating set $S$. 
We also denote $S^n$ instead of $B_n(G,S)$. Let $\beta_n(G,S)=|B_n(G,S)|$.
The {\it exponential growth rate} of $(G,S)$ is defined to be:
$$e(G,S)= \lim_{n \to \infty} \beta_n(G,S)^{\frac{1}{n}}.$$

A finitely generated group $G$ has {\it exponential growth} if there exists a finite generating set $S$ such that 
$e(G,S) >1$. The group $G$ has {\it uniform exponential growth} if there exists $c>1$, such that for every 
finite generating set $S$, $e(G,S)\ge c$. 

Given a finitely generated group $G$, we define:
$$e(G)= \inf_{|S| < \infty} e(G,S),$$
where the infimum is taken over all the finite generating sets $S$ of $G$. Since there are finitely generated groups that have
exponential growth, but do not have uniform exponential growth \cite{Wilson}, the infimum, $e(G)$, is not always obtained
by a finite generating set of a finitely generated group.

Given a finitely generated group $G$, we further define the following set in $\R$:
$$\xi(G)=\{e(G,S)| |S|<\infty\},$$
where $S$ runs over all the finite generating sets of $G$.
The set $\xi(G)$ is always countable. 

A non-elementary hyperbolic group contains a non-abelian free group, hence, it has exponential growth. In fact, a non-elementary
hyperbolic group has uniform exponential growth \cite{Koubi}. 
Recently it is proved that $\xi(G)$ of a non-elementary hyperbolic group $G$
is well-ordered (hence, in particular, has the minimum), \cite{FS}.
It was new even for free groups. 

In this paper, we deal with larger classes of groups. We state a main result.
See Definition \ref{def.acylindrical} for the definition of 
acylindricity. See Definition \ref{def.EN} for the definition 
of equational Noetherianity. 

\begin{thm}[Well-orderedness for acylindrical actions]\label{main}
Suppose $G$ acts on a $\delta$-hyperbolic space $X$
acylindrically, 
and the action is non-elementary.
Assume that there exists a constant $M$ such that for any 
finite generating set $S$ of $G$, the set $S^M$ contains a hyperbolic 
element on $X$.
Assume that $G$ is equationally 
Noetherian.
Then, 
 $\xi(G)$ is a well-ordered set.
\end{thm}
In particular, $\inf \xi(G)$ is realized by some $S$, ie,
$e(G)=e(G,S)$.

%In the theorem, we may consider a geodesic space for $X$.
%We do not lose generality if we assume $X$ to be a graph
%since we can consider the $1$-skeleton of a Rips complex of a geodesic space $X$, then the assumption still holds. 

The theorem holds under a weaker assumption, namely, 
we may replace the acylindricity of the action by that $S^M$
contains a hyperbolic and WPD element (Theorem \ref{main.proof}).
See Definition \ref{def.wpd}
for the definition of WPD.
Theorem \ref{main} is an immediate consequence of 
Theorem \ref{main.proof}
by Lemma \ref{wpd.acyl}. See the explanation 
at the beginning of the section \ref{section.3}.

We give some applications.

\begin{thm}[Theorem \ref{rel.hyp.growth}]\label{thm.rel.hyp}
Let $G$ be a group that is hyperbolic relative to a collection of 
subgroups $\{P_1, \cdots, 
P_n\}$.
Suppose $G$ is not virtually cyclic, and not equal to $P_i$
for any $i$.
Suppose each $P_i$ is finitely generated and equationally Noetherian.
Then $\xi(G)$ is well-ordered. 

\end{thm}
As examples of this theorem, we prove:

\begin{thm}[Rank-1 lattices, Theorem \ref{lattice}]\label{main.ex}
Let $G$ be one of the following groups:
\begin{enumerate}
\item
A lattice in a simple Lie group of rank-1. 
\item
The fundamental group of a complete Riemannian manifold $M$
of finite volume such that there exist $a,b>0$ with 
$-b^2 \le K \le -a^2 <0$, where $K$ denotes the
sectional curvature. 
\end{enumerate}

Then $\xi(G)$ is well-ordered.

\end{thm}

Another family of examples are $3$-manifold groups.
\begin{thm}[Theorem \ref{3manifold}]
Let $M$ be a closed orientable $3$-manifold, and $G=\pi_1(G)$.
If $M$ is one of the following, then $G$ has exponential growth and $\xi(G)$ 
is well-ordered. 
\begin{enumerate}
\item
$M$ is not irreducible and $G$ is not isomorphic to $\Z_2 * \Z_2$.
\item
$M$ is irreducible, not a torus bundle over a circle,   and $M$ has a non-trivial JSJ-decomposition.
\item
$M$ admits  hyperbolic geometry. 
\item
$M$ is  Seifert fibered such that the base orbifold is hyperbolic.
\end{enumerate}
\end{thm}

Potential examples of application 
of the main theorem are mapping class groups. We discuss this class in Section \ref{section.mcg}.
See Example \ref{ashot} for non-examples. 

We also show some finiteness result as follows. This was known for 
hyperbolic groups too, \cite{FS}.
\begin{thm}[Theorem \ref{7.1}]
Suppose the same assumption holds for $G$ as in Theorem \ref{main}.
Then for any $\rho \in \xi(G)$, up to the action of $Aut(G)$, there are at most finitely many finite generating sets $S$ such that $e(G,S)=\rho$.
\end{thm}

As a part of the proof of the main result, we show a basic result
on the growth of a group, 
generalizing a result known for hyperbolic group in \cite{AL}.
Given a group $G$ that satisfies the assumption in Theorem \ref{main}, where we do not need 
that $G$ is equationally Noetherian, 
there exists a constant $A>0$ such that 
for any finite generating set $S$ of $G$, we have
$$e(G,S) \ge A|S|^A.$$
The constant $A$ depends only on $\delta$
and the acylindricity constants. 
See Proposition \ref{bound.generators}
for the statement. Examples include mapping class groups and 
rank-1 lattices, see Example \ref{ex.bound}.

We also discuss the set of growth of subgroups in a finitely generated group $G$.
Define 
$$\Theta(G)=\{e(H,S)| S \subset G, |S| < \infty, H=\<S\>, e(H,S)>1\}.$$
The set $\Theta(G)$ is countable and contains $\xi(G)$.
If $G$ is a hyperbolic group, it is 
known by \cite[Section 5]{FS} that 
$\Theta(G)$ is well-ordered.

Similarly, we prove:

\begin{thm}[Theorem \ref{6.1}]

Suppose $G$ is one of the groups in Theorem \ref{main.ex}.
Then $\Theta(G)$ is a well-ordered set. 

\end{thm}

We also prove:
\begin{thm}[Finiteness, Theorem \ref{lattice.finite}]
%\kf{put thm number for subgroups}
Let $G$ be one of the groups in Theorem \ref{main.ex}.
%Then for each $\rho \in \xi(G)$ there
%are at most finitely many finite generating
%sets $S$, up to $\Aut(G)$, s.t. $e(G,S)=\rho$.
Then for each $\rho \in \Theta(G)$, there are 
at most finitely many $(H,S)$, up to isomorphism, 
such that $S$ is a finite generating set of $H<G$
with  $e(H,S)=\rho$.
\end{thm}

This kind of finiteness is known for hyperbolic groups \cite{FS},
and we generalize it (Theorem \ref{7.1}), which 
implies the above theorem as examples.

Some more definitions are in order in the following section.
\subsection{Acylindrical actions}
To generalize the properness of an action, 
Bowditch \cite{Bowditch} introduced the following definition. 

\begin{definition}[acylindrical action]\label{def.acylindrical}
An action of a group $G$ on a metric space $X$ is {\it acylindrical}
if for any $\epsilon >0$, there exist $R=R(\epsilon)>0$ and $N=N(\epsilon)>0$ such that 
for all $x,y \in X$ with $d(x,y) \ge R$, the set
$$\{g\in G| d(x,g(x)) \le \epsilon, d(y,g(y)) \le \epsilon\}$$
contains at most $N$ elements. 
\end{definition}

A group $G$ is called an {\it acylindrically hyperbolic} group, \cite{O},
if it acts on some $\delta$-hyperbolic space $X$ 
such that the action is acylindrical and non-elementary. Here, we say the action is {\it elementary} if the limit set of $G$ in the Gromov 
boundary $\partial X$
has at most two points. 
If the action is non-elementary, it is known that $G$ contains hyperbolic isometries.
Non-elementary hyperbolic groups and non-virtually-abelian 
mapping class groups are examples of acylindrically hyperbolic
groups, \cite{Bowditch}. There are many other examples.

\subsection{Limit groups and equational Noetherianity}
Let $G$ be a group and $\Gamma$ a finitely generated (or countable) group.
Let $\Hom(\Gamma,G)$ be the set of all homomorphisms from $\Gamma$ to $G$.

A sequence of homomorphisms $\{f_n\}$ from $\Gamma$ to $G$ is {\it stable}
if for each $g\in \Gamma$, either $f_n(g)=1$ for all sufficiently large $n$;
or $f_n(g)\not=1$ for all sufficiently large $n$.
If the sequence is stable, then the {\it stable kernel}
of the sequence, $\underrightarrow\ker(f_n)$, is defined by
$$\underrightarrow\ker(f_n)=\{g \in \Gamma| f_n(g)=1  \text{ for all sufficiently large }n\}.$$

We call the quotient $\Gamma/\underrightarrow\ker(f_n)$ a $G$-{\it limit
group}, or the limit group over $G$, associated to $\{f_n\}$, and 
the homomorphism $f:\Gamma \to \Gamma/\underrightarrow\ker(f_n)$ the {\it limit homomorphism}.
We say the sequence $\{f_n\}$ converges to $f$.

Let $G$ be a group and $F(x_1,\cdots, x_\ell)$ 
the free group on $X=\{x_1, \cdots, x_\ell \}$.
For an element $s \in F(x_1,\cdots, x_\ell )$
and $(g_1, \cdots, g_\ell) \in G^\ell$, 
let $s(g_1, \cdots, g_\ell) \in G$ denote the element
after we substitute every $x_i$ with $g_i$ 
and $x_i^{-1}$ by $g_i^{-1}$ in $s$.
Given a subset $S \subset F(x_1,\cdots, x_\ell)$,
define
$$V_G(S)=\{(g_1, \cdots, g_\ell) \in G^\ell |s(g_1, \cdots, g_\ell)=1
\text{ for all } s \in S\}.$$
$S$ is called a system of equations (with $X$
the set of variables), and 
$V_G(S)$ is called the algebraic set over $G$ defined by $S$.
We sometimes suppress $G$ from $V_G(S)$.
\begin{definition}[Equationally Noetherian]\label{def.EN}
A group $G$ is {\it equationally Noetherian}
if for every $\ell \ge 1$ and every subset $S$ in $F(x_1,\cdots, x_\ell)$,
there exists a finite subset $S_0 \subset S$ such that 
$V_G(S_0)=V_G(S)$.
\end{definition}

\begin{remark}
This definition appears in for example \cite{GrH}.
There is another version of the definition that considers
$S \subset G*F(x_1,\cdots, x_\ell)$, which is originally in \cite{BMR}
and also in \cite{RW}.
They are equivalent, see \cite[Lemma5.1]{RW}.

\end{remark}

Examples of equationally Noetherian groups
include free groups, \cite{Guba}; linear groups, \cite{BMR}; hyperbolic groups without torsion, \cite{Sela1}, then possibly with torsion, \cite{RW};
and 
hyperbolic groups relative to equationally Noetherian subgroups, \cite{GrH}.

What is important for us is the following 
general principle.

\begin{lemma}[Basic principle]\label{basic}
Let $\eta: F \to L$ be the limit map of a sequence 
of  homomorphisms, $f_n:F \to G$.
Suppose $G$ is equationally Noetherian.
Then for sufficiently large $n$, $f_n$ factors 
through $\eta$, namely, there exists
a homomorphism $h_n:L \to G$ such that 
$h_n \circ \eta =f_n$.
\end{lemma}
\proof
Let $X=\{x_1, \cdots, x_\ell\}$ and 
suppose $F=F(X)$.
Let $R=\{r_i\} \subset F(X)$ be a set of defining relations for $L$.
In general this is an infinite set.
Each $r_i$ is a word on $X$, so that we can 
see $R$ as a system of equations with $X$
the variable set. 
%(This is special in the sense that no coefficients from $G$.)
Since $G$ is equationally Noetherian,
there is a finite subset $R_0 \subset R$
such that $V(R)=V(R_0)$, namely,
every solution (an element in $G^\ell$) for $R_0$ is
a solution for $R$.

Now, since $\eta(r_1)=1$ in $L$ for a large enough $n$, we have 
$f_n(r_1)=1$ in $G$ since $\eta$ is the limit of $\{f_n\}$.
By the same reason, since $R_0$ is a finite set, there exists $N$ such that 
for every $n \ge N$, we have 
$f_n(r_i)=1$ in $G$ for all $r_i \in R_0$.
In other words, if $n \ge N$, then $(f_n(x_1), \cdots, f_n(x_\ell))\in V(R_0)$.
But since $V(R_0)=V(R)$, 
this implies that if $n \ge N$, then 
$(f_n(x_1), \cdots, f_n(x_\ell))\in V(R)$, namely, 
$f_n(r_i)=1$ in $G$ for all $r_i \in R$.
Since $R$ is a system of defining relations
for $L$, it implies that each $f_n$ with $n \ge N$ factors through
$\eta:F \to L$.
\qed

Regarding Theorem \ref{main}, there exists a finitely generated group $G$ that acylindrically acts on a $\delta$-hyperbolic
space, in fact a simplicial tree $T$,
 in the non-elementary way such that for any finite generating set $S$, the set
 $S^2$ contains a hyperbolic isometry on $T$; and that  there is no 
finite generating set $S$ with $e(G)=e(G,S)>1$.
In particular $\xi(G)$ is not well-ordered.

The following example is pointed out by Ashot Minasyan.
A group $G$ is called {\it Hopfian} if every  surjective
homomorphism $f:G \to G$ is an isomorphism. 
\begin{example}\label{ashot}
Take a finitely generated non-Hopfian group, $G$, for example $BS(2,3)$, (a Baumslag-Solitar group). 
Put $H=G* \Bbb Z$. Then $H$ is non-Hopfian, ie,
there exists a surjection $f:H \to H$ that is not an isomorphism. It is a standard fact that 
a finitely generated equationally Noetherian group
is Hopfian, therefore $H$ is not equationally Noetherian,
so that Theorem \ref{main} does not apply to $H$. 

But all other assumptions in the theorem are satisfied. 
Let $T$ be the Bass-Serre tree for $G*\Bbb Z$. The tree
$T$ is 0-hyperbolic. 
The action of $H$ on $T$ is acylindrical and non-elementary.
Also, for any finite generating set $S$ of $H$, it is a well-known lemma (due to Serre, cf. \cite{BrF}) that $S^2$ contains a hyperbolic 
isometry on $T$. 

On the other hand, it is known, \cite{Sa}, that a finitely generated group $K$ that is a free product is {\it growth-tight}, namely,
for any surjective homomorphism $h:K \to K$ that is not 
an isomorphism, and for any finite generating set $S$ of $K$, 
$e(K,S) > e(K,h(S))$.
Now, it follows that $e(H)$ is not achieved by any $S$, since
if it did, then take such $S$. But then, the non-isomorphic,
surjective homomorphism $f$ in the above would imply $e(H,S) > e(H,f(S))$, 
a contradiction.
\end{example}
\vskip .05in
\noindent
{\bf Acknowledgement}.
This work is a continuation of \cite{FS}.
The author would like to thank Zlil Sela for generously 
offering 
numerous comments and suggestions throughout
the work.
He is grateful to Mladen Bestvina, Emmanuel Breuillard, Thomas Delzant, Daniel Groves and Wenyuan Yang for useful comments. 

He would like to thank the referee for reading the paper 
thoroughly and making many suggestions, by which 
the paper significantly improved not only in terms of the presentation,
but also  mathematically, for example on Proposition \ref{2.2}.

\section{Lower bound of a growth rate}

Although most statements in the paper are for geodesic spaces $X$, 
{\em we consider a graph for $X$ instead of a geodesic space
in the arguments throughout the paper unless we indicate otherwise.}
The advantage is that an infimum is achieved for
various notions, for example, $L(S)$,  in the  arguments. 
But, we do not lose generality if we assume $X$ to be a graph
since we can consider the $1$-skeleton of a Rips complex of a geodesic space with a group action. Also, 
the various assumptions we consider
 (such as the hyperbolicity of the space, 
acylindricity of the action)
remain valid, maybe with slightly difference constants.

\subsection{Hyperbolic isometries and axes}

%We recall some definitions and results from \cite{BrF}. 
Suppose a group $G$ acts on a geodesic space $X$ by isometries. 
Choose a base point $x\in X$ and for $g\in G$, put
$$L(g)= \inf_{x \in X} |x-g(x)|, \lambda(g) = \lim_{n\to \infty} \frac{|x-g^n(x)|}{n}.$$
$L(g)$ is called the {\it minimal displacement}.
$\lambda(g)$ is called the {\it translation length}, and does not depend
on the choice of $x$, and for any $n>0$ we have
$\lambda(g^n)=n\lambda(g)$.
We also have
$$ \lambda(g) \le L(g) \le \lambda(g) + 7\delta.$$
The first inequality is trivial and  we leave the second as an exercise
(for example, use \cite[Corollary 1]{AL}).
%\kf{check ref}
For a finite set $S \subset G$ and $x\in X$, define 
$$L(S,x)=\max_{s \in S} |x-s(x)|,$$
then define 
$$L(S)=\inf_{x \in X} L(S,x).$$

We recall a few definitions and facts from $\delta$-hyperbolic spaces.
%An isometry $g$ of a $\delta$-hyperbolic space $X$ is called {\it hyperbolic} if 
%$\lambda (g) >0$.
%If $g$ is hyperbolic, then there exists a bi-infinite path, called an {\it axis},
%$\gamma$ such that the Hausdorff distance between $\gamma$ and
%$g^n(\gamma)$ is at most $10 \delta$ for any $n \in \Bbb Z$, 
An isometry $g$ of a hyperbolic space $X$ is called {\it elliptic} if
the orbit of a point by $g$ is bounded, and 
{\it hyperbolic} if there are $x\in X$ and $C>0$ such that 
for any $n>0$, we have $|x-g^n(x)| > Cn$.
The element $g$ is hyperbolic iff $\lambda(g)>0$.

A hyperbolic isometry $g$ is associated with a bi-infinite quasi-geodesic, $\gamma$, called an  {\it axis} in $X$.
If there exists a bi-infinite geodesic $\gamma$
that is invariant by $g$, that would be an ideal choice
for an axis, but that is not always the case.

As a remedy,  if $L(g)\ge 10\delta$, 
take a point $x\in X$ where $L(g)$ is achieved.
Then take a geodesic $[x,g(x)]$ between $x$ and $g(x)$ and consider
the union of its $g$-orbit, which defines a $g$-invariant path,
(see for example \cite{De}).
If $L(g) < 10\delta$, then take $n>0$ such that 
$L(g^n) \ge 10\delta$ and apply the construction 
to $g^n$, and use this path for $g$, which  is not $g$-invariant. 
We denote this axis as $A(g)$ in this paper. 
Also, for $g^n$ with $g\not=0$, we  may also take $A(g)$
as an axis for $g^n$. 

For $g$, an axis $A(g)$ is not  unique, but
uniformly (over all hyperbolic $g$) quasi-geodesic, 
such that for  any two points $x,y \in A(g)$, 
the Hausdorff distance between the segment
between $x,y$ on $A(g)$ and 
a geodesic between $x,y$, $[x,y]$, is 
at most $10\delta$.
Also, if $Hd(A(g), h(A(g))$ for $h \in G$ is finite, 
then it is bounded by $10\delta$, where $Hd$ is the 
Hausdorff distance. 
We sometimes call $A(g)$ a
 {\it $10\delta$-axes}.
 We consider a direction on the $10\delta$-axis using 
 the action of $g$.
 %\kf{check}

A hyperbolic isometry 
$g$ defines two limit points in $\partial X$, the visual(Gromov) boundary of $X$,
by $g^{\infty}=\lim_{n \to \infty} g^n(x), g^{-\infty}=\lim_{n \to -\infty} g^n(x)$, where $x$ is a base point. 
We say two hyperbolic isometries $g,h$ are {\it independent}
if $\{g^{\pm \infty}\}$ and $\{h^{\pm \infty}\}$ are disjoint.
If the Hausdorff distance between two axes is finite, then we say
they are {\it parallel}.

\subsection{WPD elements}
We consider another version of properness of 
a group action that is weaker than acylindricity.

%\kf{changed}
\begin{definition}[WPD, uniformly WPD, $D$-WPD]\label{def.wpd}
Let $G$ act on a $\delta$-hyperbolic space $X$.
Suppose that  $g \in G$ is hyperbolic on $X$.
We say $g$ is {\it WPD}  if
there is a $10\delta$-axis, $\gamma$, of $g$ such that 
for any $\epsilon>0$, 
there exists $D=D(\epsilon)>0$ such that for any $x, y \in \gamma$ 
with $|x-y| \ge D \lambda(g)$, the number
of the elements in the following set is at most $D$:
\begin{equation}\label{set.wpd}
\{h \in G| |h(x)-x| \le \epsilon, |h(y)-y| \le \epsilon \}.
\end{equation}
In application we often take $y=g^D(x)$. 
If we want to make the function $D$ explicit, we say that 
$g$ is $D$-WPD, or WPD w.r.t. $D$. We say that $D$ is 
a function for WPD. 

If there is a function $D$ such that if a set of hyperbolic elements
in $G$ are $D$-WPD, then we say they are
{\it uniformly WPD}, or {\it uniformly $D$-WPD}. 
If all hyperbolic elements in $G$ are uniformly $D$-WPD, then 
we say the action is uniformly ($D$-)WPD. 
\end{definition}

Some remarks are in order.  
The notion of WPD (weak proper discontinuity) was introduced in \cite{BeF}, where the function $D$ is not used,
but the definitions are equivalent. 

If $g$ is $D$-WPD, then it is $D'$-WPD
for any function $D'$ such that $D'(\epsilon)
\ge D(\epsilon)$ for all $\epsilon$.
So, without loss of generality, we assume that $D(\epsilon)$ does not decrease when we increase $\epsilon$.
We often use the value $D(100\delta)$ in this paper, for example, 
see Lemma \ref{lower.bound}.
For convenience we also assume that $D(100\delta) \ge 50$, which we use in the 
proof of Lemma \ref{lemma7.4}.

The choice of a $10\delta$-axis is not important in the definition.
Also, one can use $C$-axis for any $C>0$. It only changes
the function $D(\epsilon)$. 
Uniformly WPD is  related to but weaker than
 the notion of weak acylindricity in \cite{De}.
 See also \cite[Example 1]{De} for the difference between
 acylindricity and variations of WPD.

For an acylindrical action, it is known (\cite[Lemma 2.1]{F}) that 
there exists $T>0$ such that for any hyperbolic element $g$, we have
$\lambda(g) \ge T$.
This holds for uniformly WPD actions too, and
the argument is same, but for the readers' convenience,
we prove it.

\begin{lemma}[Lower bound on $\lambda(g)$]\label{lower.bound}
Suppose $G$ acts on a $\delta$-hyperbolic space $X$.
Let $g \in G$ be hyperbolic with a $10\delta$-axis $\gamma$.
\begin{enumerate}
\item
If $g$ is $D$-WPD, then 
$\lambda(g) \ge \frac{50\delta}{D(100\delta)}$.
\item
If the action is acylindrical, then 
$\lambda(g) \ge \frac{50\delta}{N(150\delta)}$.
\end{enumerate}

\end{lemma}

\proof
(1) 
Set $\D=D(100\delta)$. 
Let $x, y \in \gamma$ with 
$|x-y| \ge \D \lambda(g)$. Then there must be some $n$ with 
$0\le n \le \D$ such that 
$|x-g^n(x) | > 100\delta$ or 
$|y-g^n(y) | > 100\delta$.
This is because otherwise, 
all the elements $1,g,\cdots,g^{\D}$, which are $\D+1$
distinct elements,  are contained in 
the set 
$$
\{h\in G||h(x)-x|\le 100\delta, \text{ and } |h(y)-y| \le 100\delta\}.
$$
This is impossible since the contains at most $\D$ elements. 
Now suppose, say, $x$ satisfies  $|x-g^n(x) | > 100\delta$. Then it  implies that $\lambda(g^n) \ge 50\delta$.
Since $n\le \D$, we have $\lambda(g) \ge \frac{50\delta}{\D}$.

(2)
The argument is similar to (1).
Take $x,y \in \gamma$ with $|x-y|= R(150\delta)$.
If $|x-g^n(x)| \le 100\delta$ for some $n$, then 
$|y-g^n(y)| \le 150\delta$.
This implies, by acylindricity, there must be $n$ with $1\le n \le N(150\delta)$
s.t. $|x-g^n(x)| > 100\delta$. It follows 
$\lambda(g^n) \ge 50\delta$, so that 
$\lambda(g) \ge \frac{50\delta}{N(150\delta)}$.
\qed

For a function $D(\epsilon)$, put
$$T= \frac{50\delta}{D(100\delta)},$$
then by Lemma \ref{lower.bound} (1) we have 
$\lambda(g) \ge T$ for a $D$-WPD element $g$. 

%We denote the lower bounds of $\lambda(g)$ in Lemma \ref{lower.bound} by $T$. 
%For convenience, we may retake $T$ to be in the form of
%$$T=\frac{1}{p}$$
%for some $p \in \Bbb N$, so that 
%$\frac{\delta}{T}=p\delta$.

 The following lemma is straightforward. 
 
 \begin{lemma}[Acylindricity implies uniform WPD]\label{wpd.acyl}
 If an action of $G$ on a $\delta$-hyperbolic space $X$
 is acylindrical, then it is uniformly WPD.
 \end{lemma}
 
 \proof
 %\kf{rewrote it}
 Suppose $g\in G$ is hyperbolic on $X$,
 and let $\gamma$ be a $10\delta$-axis.
 Let $R(\epsilon), N(\epsilon)$ be the 
 acylindricity constants.
 Also, let $T>0$ be a uniform bound for $\lambda(g)\ge T$ by Lemma \ref{lower.bound} (2). $T$ does not depend on $g$ nor $\epsilon$. 
 Suppose $\epsilon$ is given. 
 Let $K=K(\epsilon)$ be a smallest integer with $K\ge \frac{R(\epsilon)}{T}$. Then $\lambda(g^K) \ge R(\epsilon)$.
 The constant $K$ does not depend on $g$.
 Then for any $x\in \gamma$ and $n \ge K$, we have $|x-g^n(x)| \ge R(\epsilon)$.
 By acylindricity, there are at most $N(\epsilon)$ elements
 which simultaneously move each of $x, g^n(x)$ by at most $\epsilon$.
 Put $D(\epsilon)=\max\{K(\epsilon),N(\epsilon)\}$,
 then the action is uniformly $D$-WPD.
 \qed

We state a lemma which is useful for us. 
\begin{lemma}\label{lemma.wpd}
 Suppose there are at most $D$ elements 
 that satisfies the condition  (\ref{set.wpd}) in the definition \ref{def.wpd} for 
 $\epsilon=100\delta$ if $|x-y| \ge D \lambda (g)$. Then $g$ is WPD, and moreover, there is a function $D'$ with $D'(100\delta)=D$
 that depends only on $D, \delta$ such that 
 $g$ is $D'$-WPD.
 \end{lemma}
 
 \proof
 Let $\gamma$ be a $10\delta$-axis of $g$.
 Suppose $\epsilon >0$ is given. We may assume $\epsilon > 100 \delta$. Take $x,y \in \gamma$ such that 
 $|x-y| \ge D \lambda(g)+ 2 \epsilon + 1000 \delta$.
 Take $p,q \in \gamma$ between $x$ and $y$ with 
 $|x-p|=|y-q|=\epsilon +50\delta$.
 Then $|p-q| \ge D \lambda(g) +800\delta$.

 Let $J$ be the collection of elements in $G$ such that 
 $|x-j(x)|\le \epsilon$ and $|y-j(y)| \le \epsilon$. If $j \in J$, then 
 $|p-j(p)| \le \epsilon +30\delta$ and $|q-j(q)| \le \epsilon + 30 \delta$; and 
 $j(p), j(q)$ are in the $15\delta$-neighborhood of $\gamma$.
 It implies that $J$ contains a subset $J_0$ that 
 contains at most $(2\epsilon+200\delta)/(10\delta)$ elements such that 
 for any $j \in J$, one can find $j_0 \in J_0$ with 
 $|p-j^{-1}j_0(p)| \le 70 \delta$.
 To see it, consider points $p_1, p_2 \in \gamma$ with $|p_1-p_2|=2\epsilon +200\delta$ such that $p$ is the mid point of the segment between
 $p_1,p_2$ on $\gamma$, which we denote by $[p_1,p_2]$. Then the points $j(p)$ with $j \in J$
 is contained in the $30\delta$-neighborhood of $[p_1,p_2]$.
 One should imagine that this neighborhood is a narrow tube
 around $[p_1,p_2]$. 
 Now by a pigeon hole argument, one can find a desirable subset $J_0$.
 (Notice that $|p-j^{-1}j_0(p)|=|j(p)-j_0(p)|$, so that 
 one needs to find a point $j_0(p)$ near (i.e. at most $70\delta$) a given point $j(p)$, which
 is possible.)

 But $|p-j^{-1}j_0(p)| \le 70 \delta$ implies $|q-j^{-1}j_0(q)| \le 100 \delta$.
 This is because $|p-q| \ge 800 \delta$ and both $[j(p),j(q)]$
 and $[j_0(p),j_0(q)]$ are contained in the $20\delta$-neighborhood of $\gamma$. 
 We have shown that the element $j^{-1}j_0$
 moves both $p,q$ by at most $100\delta$. By our assumption,  there are at most $D$ possibilities for such element. 
 In conclusion,  $J$ contains at most $D\times (2 \epsilon +200\delta)/(10\delta)$ elements. We proved that $g$ is WPD.
 
 We compute a WPD-function for $g$, which 
 we denote by $D'$. By assumption, we may set
 $D'(100\delta)=D$. First, 
  $$D \times (2\epsilon +200\delta)/(10\delta)=D\times(\epsilon/(5\delta) +20).$$
  Next, 
  $$D+(2\epsilon +1000\delta)/\lambda(g)
   \ge D+D(2\epsilon + 1000 \delta)/(50\delta)
   =D(21+\epsilon/(25\delta))
   $$ by Lemma \ref{lower.bound} (1).
  So, if $\epsilon > 100\delta$, we set $$D'(\epsilon)= D \max\{
  20+\epsilon/(5\delta),   
  21+ \epsilon/(25\delta) \}
  =D \times (21+\epsilon/(5\delta) ).$$  
  \qed

\subsection{Elementary closure}

Suppose $G$ acts on a hyperbolic space  $X$ and let $g\in G$ be a 
hyperbolic isometry with an axis $\gamma$. The {\it elementary closure}
of $g$ is defined by
$$E(g)=\{h\in G| Hd(\gamma,h(\gamma)) < \infty \}.$$
It turns out that $E(g)$ is a subgroup of $G$.
Clearly, $\<g\> < E(g)$.

We denote the $a$-neighborhood of a subset $Y\subset X$ by
$N_a(Y)$.

\begin{lemma}[parallel axes]\label{2.1}
Suppose $G$ acts on a $\delta$-hyperbolic space $X$.
Let $g\in G$ be hyperbolic with a $10\delta$-axis $\gamma$.
Let $h\in G$, then 
\begin{enumerate}
\item
If $h \in E(g)$, then $Hd(\gamma, h(\gamma)) \le 50\delta$.
\item
Assume that $g$ is $D$-WPD.
If $h \not\in E(g)$, then the diameter of 
$h(\gamma) \cap N_{50\delta}(\gamma)$ is at most  
%$$ R(100\delta) + L(g)N(100\delta) + 100\delta.$$
$$2D(100\delta) L(g) + 100\delta.$$
\end{enumerate}
%\kf{rewriting}

\end{lemma}
This lemma is well-known in slightly different versions,
for example \cite[Lemma 2.2]{F} for acylindrical actions,
so the proof will be brief. 

In general if two axes have finite Hausdorff distance, 
then we say they are {\em parallel}.
%The property (2)  says that if two 
%axes $h(\gamma), k(\gamma)$ are ``parallel'' (in the sense 
%that one is contained in the $50\delta$-neighborhood of the other)
%along a segment longer than $2D L(g) + 100\delta$, then they are 
%parallel. 

\proof
(1) By definition, $Hd(\gamma, h(\gamma)) < \infty$. Since
both $\gamma, h(\gamma)$ are $10\delta$-axes, we have a desired bound. 

(2)
Suppose not.  Suppose that the direction of
$\gamma, h\gamma$ coincide along the parallel part. 
Set $\D=D(100\delta)$. 
 Take $x \in h(\gamma)$ 
 near one end of the intersection such that 
 $g^{n}(x)$ for $0\le n \le 2\D$ are in the $50\delta$-neighborhood of the intersection. This is possible since the intersection is long 
 enough. 
 Consider the points $x, g^\D(x)$. Then $|x-g^\D(x)| \ge \D \lambda(g)$. 
 Letting $(hgh^{-1})^{-n} g^{n}$ with $0\le n \le \D$
 act on $x, g^\D(x)$, we have 
 $$|x-(hgh^{-1})^{-n} g^{n}(x)| \le 100 \delta, 
 \, \, 
 |g^\D(x)-(hgh^{-1})^{-n} g^{n}(g^\D(x))| \le 100 \delta$$
 for all $0 \le n \le \D$. 
 
But since $g$ is $D$-WPD,  there are at most $\D$ such elements, so that 
it must be that 
 for some $n\not=m$, we have
 $(hgh^{-1})^{-n} g^{n}=(hgh^{-1})^{-m} g^{m}$, 
 so that $g^{n-m}$ and $h$ commute.
 It implies that $\gamma$ and $h(\gamma)$ are parallel, 
 a contradiction. 
 
 If the direction for $\gamma$ and $h\gamma$ are opposite,
 consider $(hgh^{-1})^n g^n$ instead of $(hgh^{-1})^{-n} g^{n}$ , and the rest is same.
 \qed

We quote a fact.

\begin{prop}[elementary closure]\label{prop.elementary}
Suppose $G$ acts on a $\delta$-hyperbolic space $X$.
Let $g\in G$ be a hyperbolic element on $X$. Assume $g$ is $D$-WPD.
Then,
\begin{enumerate}
\item
$E(g)$ is virtually $\Bbb Z$, and contains $\<g\>$ as a finite 
index subgroup. 
\item
If $h\in G$ is a hyperbolic element such that 
$g$ and $h$ are independent,
then $E(g) \cap E(h)$ is finite. 
\end{enumerate}
\end{prop}
\proof
(1) This is \cite[Lemma 6.5]{DGO}.

(2) If $E(g) \cap E(h)$ is infinite, then it contains $\< g^N \>$ for some 
$N>0$. But since $g$ and $h$ are independent, 
for a sufficiently large $m$,  we have $g^m \not\in E(h)$, impossible. 
\qed

Note that under the assumption of the proposition, 
the action of $G$ is elementary if and only if $G$ is virtually $\Z$, which is equivalent to  that $G=E(g)$ in this case. 

In general, if a group $H$ is virtually $\Bbb Z$, then 
there exists an exact sequence
$$1 \to F \to H \to C \to 1,$$
where $C$ is either $\Bbb Z$ or $\Bbb Z_2 * \Bbb Z_2$, 
and $F$ is finite. 
In the case that $H$ is $E(g)$ in the above, we denote the finite group $F$
by $F(g)$.
If $C=\Bbb Z$,
$F(g)$ is the set of elements
of finite order in $E(g)$.

The axis $\gamma$ defines two points in the ideal boundary of $X$,
which we denote $\{\gamma(\infty), \gamma(-\infty)\}$.
$E(g)$ is exactly the set of elements that leaves this set 
invariant. If $h \in E(g)$ swaps those two points,
we say it {\em flips} the axis since it flips the direction 
of the axis $\gamma$.

%For a pseudo-Anosov element $g \in MCG$, let $E(g)< MCG$
%be the elementary closure of $g$. $E(g)$ is virtually $\Z$.

If the action of $G$ on $X$ is $D$-uniform WPD, then 
%By $D$-uniform WPD 
%there exists a constant $c$ that depends on the action 
 for any hyperbolic $g$, we have $|F(g)| \le 2D(100\delta)$. 
 %\kf{isn's this also D? Maybe, but I keep it as is.}
 Moreover, if $C=\Bbb Z$, then $|F(g)| \le D(100\delta)$.
%Indeed, one can take $c=2D$, and $c=D$
%if $C=\Bbb Z$.

Note that $\<g\> < E(g)$. 
If the action is non-elementary, then $E(g)\not=G$, which
is equivalent to that $G$ is not virtually $\Bbb Z$,
so that $G$ contains two independent hyperbolic isometries. 

 We say a hyperbolic element $g$ is {\it primitive} if 
$C=\Bbb Z$  and 
each element $h\in E(g)$ is written as
$h=fg^n$ for some $f \in F(g)$ and $n \in \Bbb N$.

\subsection{About the constants}
From now on, there will be 
many constants to make
the argument concrete and precise.
In the argument we consider a sequence of generators
$S_n$.
It would be a good idea to keep in mind that 
there are two kinds of constants.

The first kind are those that are fixed once
the constants $\delta, D(100\delta), M$ are given by the action:
$$\delta, D(100\delta), M; T, k, m, b.$$
The constants $k,m$ will appear in Lemma \ref{2.4} and 
$b$ in the proof of Lemma \ref{4.2}.

The second kind are those that depend
on a generating set $S_n$ of $G$:
$$L(S_n), \lambda(g), L(g), \lambda(u), L(S_n^{2MD}), \Delta_n.$$
 If $L(S_n)$ is bounded, then all of 
the constants in the second will be bounded, but 
if not, then, roughly they all diverge in the same order
as $L(S_n)$.

We remark that if $L(S_n)$ diverges, then one way
to argue is to rescale $X$ by $\frac{1}{L(S_n)}$, 
then go to the limit, which is a tree, and use the geometry of the tree. 
This approach is the one taken in  \cite{FS}, where
only this case happens. But in our setting, the new feature is that 
possibly, $L(S_n)$ is bounded. In this paper, we use a unified approach. 

For convenience, we assume $\delta>0$ from now on. 
\subsection{Lower bound of a growth rate}
\begin{lemma}[hyperbolic element of large displacement]\label{2.3}
Let $X$ be a $\delta$-hyperbolic space and $S$ 
a finite set of isometries of $X$. Suppose that $L(S) \ge 30 \delta$.
Let $x\in X$ be such that $L(S)=L(S,x)$.
Then there is a hyperbolic element $g \in S^2$ such that 
$$L(S)-8\delta \le |x-g(x)|$$ and 
$$    |x-g(x)| - 16\delta \le L(g).$$

\end{lemma}

This is exactly same as a part of \cite[Lemma 7]{AL}
(in that paper, our element $g$ is denoted as $b$ in the proof), 
although their setting is that $X$ is a Cayley graph of 
a hyperbolic group $G$ and $S \subset G$.

The proof is same verbatim after a suitable translation of 
notions, so we omit it
(cf. \cite[Theorem 1.4]{BrF} for somewhat similar result).

Note that if $s \in S$ then $|x-s(x)| \le L(S)$;
and if $s\in S^2$ then $L(s) \le |x-s(x)| \le 2L(S)$ 
by the definition of $L(S)$ and the choice of $x$.

We summarize the properties of $g$ we use later:
\begin{itemize}
\item
$L(S) -24 \delta \le L(g) \le 2L(S).$
\item
The distance between $x$ and the $10\delta$-axis of $g$ is 
at most $ 20 \delta$.
\end{itemize}
The second one follows from 
$L(g)=\min_{y\in X} |y-g(y)| \ge |x-g(x)| - 16\delta$.

\begin{lemma}[free subgroup with primitive hyperbolic elements]\label{2.4}
Let $X$ be a $\delta$-hyperbolic geodesic space and $G$ a group acting 
on $X$. Assume $G$ is not virtually cyclic.
Let $D(\epsilon)$ be a function for WPD. 
%Then there exist constants $n,m,k$ with the following property:
Set
$$k=60D(200\delta), m=66k+4=3960D(200\delta)+4.$$

Let $S$ be a finite set that generates $G$ with $L(S) \ge 50 \delta$,
and  $x\in X$ with 
%\kf{changed 30d to 50d}
$L(S)=L(S,x)$. Suppose $g\in S^2$ is a hyperbolic element 
with a $10\delta$-axis $\gamma$ with $d(x,\gamma) \le 20 \delta$
such that $$L(g) \ge |x-g(x)| -16\delta \ge L(S)-24\delta.$$
Assume that $g$ is $D$-WPD. 

Then there exists $s\in S$ such that 
$g^k, sg^k s^{-1}$ are independent and freely generate a rank-$2$ free group $F<G$ with the following property.
There is a WPD-function $D'(\epsilon)$ with
$$D'(100\delta)=D(200\delta)$$
such that every
non-trivial element $h\in F$ is hyperbolic on $X$
and $D'$-WPD. Also, 
$h$ satisfies: 
$$\lambda(h) \ge 10(2D(200\delta) L(g) +100\delta) \ge 10L(S).$$

Moreover, there is an element 
$u \in F$ that satisfies:
\begin{enumerate}
\item $u$ has a $10\delta$-axis $\alpha$ with $d(x,\alpha) \le 50 \delta$.

\item
$u$ is primitive, namely, there exists
$$1 \to F(u) \to E(u) \to \Bbb Z \to 1$$
such that any element $h \in E(u)$ is written as
$h=f u^p$ with $f \in F(u)$ and $p \in \Bbb Z$.

\item
%\kf{maybe k=D}
$|F(u)| \le D'(100\delta)=D(200\delta)$.
%\kf{maybe $D(100\delta)$ is enough}
\item
$u \in S^m$.

\end{enumerate}

\end{lemma}

%We note that $n,m,k$ depends only on $\delta$ and $D$.
By $L(S) \ge 50\delta$, we have $L(g) \ge 26\delta$, $2L(g) \ge L(S)$, and $2\lambda(g) \ge L(g)$.
%Since $g^n, sg^ns^{-1}$ are independent, 
%the intersection of the $50\delta$-neighborhood of the $20\delta$-axis of $g$
%and the axis of $sgs^{-1}$ is at most 
%$R(100\delta)+\lambda(g)N(100\delta)+ 100\delta$.

\begin{remark}\label{2.5}
In the proof of Proposition \ref{bound.generators} we will apply Lemma \ref{2.4} to $S^{MD(200\delta)}$ instead of a generating set $S$ itself. In that case, the element 
$s\in S^{MD(200\delta)}$ that appears in the above lemma can be chosen from $S$
itself. We choose such $s$ in the beginning of the proof of the lemma and the rest of 
the argument is exactly same. 

\end{remark}

\proof
Set $\D=D(200\delta)$ and $T=\frac{50\delta}{\D}$.
Since $g$ is $D$-WPD, by
Lemma \ref{lower.bound} (1), we have 
$\lambda(g) \ge T$.
Remember that by our assumption, we always have $D(100\delta) \le D(200\delta)$, 
so that the estimate in the lemma holds for $D(200\delta)$ as well. 

Since $G$ is not virtually cyclic, there is an element $s\in S$ with 
$s\not\in E(g)$.
By Lemma \ref{2.1}(2), the diameter of the intersection $s(\gamma) \cap N_{50\delta}(\gamma)$
is at most 
$$ 2\D L(g) +100\delta.$$

In view of this, set $$k=10\left(4\D +\frac{100\delta}{T}\right)
=60\D  \ge 10.$$
%where $T$ is a lower bound of $\lambda(g)$ from Lemma \ref{lower.bound}.
Then
\begin{align*}
\lambda(g^k) \ge k\lambda(g) &\ge 10(4\D\lambda(g) + 100\delta)
\\
&\ge 10(2\D L(g) +100\delta).
\end{align*}
Here, we used $2\lambda(g) \ge L(g), \lambda(g) \ge T$.

Note that $sgs^{-1}$ is hyperbolic and $D$-WPD with a $10\delta$-axis $s(\gamma)$. 
Since $\lambda(g^k)$ is at least 10 times longer than the
above intersection,  $g^k$ and $sg^k s^{-1}$ 
freely generate a free group, $F$. Also, 
its non-trivial elements $h$ are hyperbolic and $$\lambda(h)  
\ge  \lambda(g^k) \ge 10(2\D L(g) +100\delta) \ge 10L(S).$$
For the last inequality, we used $2L(g)  \ge L(S)$.

We argue that there is a function $D'(\epsilon)$
with $D'(100\delta)=D(200\delta)$
such that every non-trivial $h\in F$ is uniformly $D'$-WPD.
We only need to argue for $\epsilon=100\delta$ to check WPD
by Lemma \ref{lemma.wpd}.
Let $\gamma$ be a $10\delta$-axis of $h$.
We will show that 
for any $ y, z \in \gamma$ with $|y-z| \ge 3 \lambda(h)$, 
there are at most $D(200\delta)$ elements $j \in G$ satisfying 
$$|j(y)-y| \le 100\delta \text{ and }  |j (z)-z|
\le 100\delta.$$
We denote the collection of those elements $j$ by $J$. 
We will use the fact that the axis of  $h$ is, roughly speaking, a concatenation 
of some translates of the segment $[x,g^k(x)]$
by elements in $F$ (cf. a more precise description of the axis of the 
element $u$ below.)
Also, we point out that without loss of generality, one may 
take a conjugate of $h$ in $G$ in the argument. 
Since $h$ is a word on $g^k$ and $sg^ks^{-1}$, by taking a conjugate of $h$ in $G$, one may assume that 
the segment $[x,g^k(x)]$ is contained in the $30\delta$-neighborhood
of the segment $[y,z]$ except for some small neighborhood of $x$
and $g^k(x)$. Taking a further conjugate of $h$ if necessary, 
one may assume that the segment $[x,g^{k/2}(x)]$
is contained in the $30\delta$-neighborhood of $[y,z]$. 
Here, we used $|y-z| \ge 3 \lambda(h)$.

That implies that one can find points $p,q \in [x, g^{k/2}(x)]$ 
with $|p-q| \ge \frac{k}{4}\lambda(g) =15\D \lambda(g)$
such that $|p-j(p)| \le 150 \delta$ and $|q-j(q)| \le 150 \delta$
if $j \in J$. 

But since $g$ is $D$-WPD where $\D=D(200\delta)$, 
the set $J$ contains at most $D(200\delta)$ elements. 
(Strictly speaking the points $p$ and $q$ are maybe not 
exactly on a $10\delta$-axis of $g$, but one can choose
nearby points of $p$ and $q$ on the axis and argue.)
We showed that $h$ is uniformly WPD.

%\kf{$\alpha$ is not axis in our sense}

\begin{figure}[h]
%\hspace*{-3.3cm}     
\begin{center}                                              
   \includegraphics[scale=0.4]{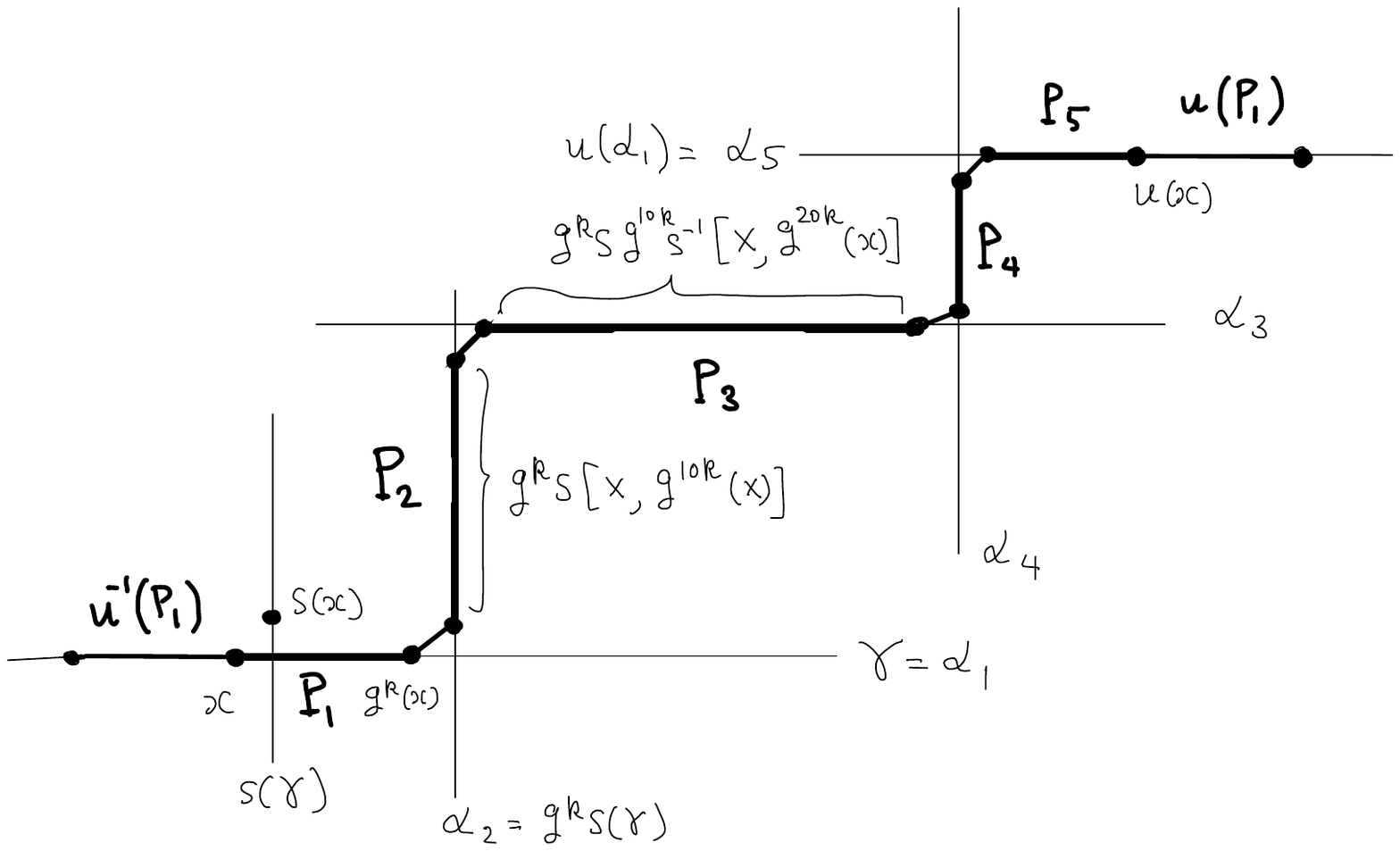}
   \end{center}
\caption{The thick line is $\bar \alpha$. It is a broken geodesic from $x$ to $u(x)$ with four short ``gaps''. The first gap is $[g^k(x), g^ks(x)]$.
}
\label{fig.separator}
\end{figure}

Now we argue for the moreover part. 
We define an element $u$ by 
$$u = g^k(sg^{10k}s^{-1}) g^{20k}(sg^ks^{-1})g^k.$$
%We will argue that  $u$ is hyperbolic
%by exhibiting  an axis, $\alpha$.

(1) To construct  an axis $\alpha$ of $u$, 
 consider the following (see Figure \ref{fig.separator}):
\begin{align*}
[x,g^k(x)] 
\cup& g^ks [x,g^{10k}(x)] \cup g^k(sg^{10k}s^{-1}) [x,g^{20k}(x)] 
\\
\cup& g^k(sg^{10k}s^{-1}) g^{20k}s [x,g^k(x)] 
\cup g^k(sg^{10k}s^{-1}) g^{20k}(sg^ks^{-1})[x,g^k(x)], 
\end{align*}
which consists of five geodesic segments, which we call {\it pieces}.
We name them as $P_1, \cdots, P_5$.
Their 
length are $U, 10U,20U,U,U$ for some constant $U>0$ maybe with some error
up to $100\delta$. Remember that 
$$U\ge \lambda(g^k) \ge 10L(S).$$

Note that $P_1$ 
is contained in $\gamma$, which is the axis of $g$.
Put $\alpha_1=\gamma$. 
Similarly, $P_2, \cdots, P_5$  are contained
in the axes of the conjugates of $g$ by $g^ks, g^k(sg^{10k}s^{-1}),
g^k(sg^{10k}s^{-1}) g^{20k}s, g^k(sg^{10k}s^{-1}) g^{20k}(sg^ks^{-1})$,
respectively. 
We will call those axes 
 $\alpha_2,\alpha_3,\alpha_4,\alpha_5$. 
 % Note that $P_i$ is (uniquely) contained in $\alpha_i$. 
 %Note that $\alpha_5=u(\alpha_1)$.
 
There are four short ($\le L(S)$) gaps between 
$P_i$ and $P_{i+1}$. 
We put (short) geodesics between the gaps and 
obtain a path, $\bar \alpha$.
Then take the union of the
$u$-orbit of $\bar \alpha$, which we call $\hat \alpha$.
This is an axis for $u$ since $g$ and $sgs^{-1}$ are independent
and their $10\delta$-axes stay close ($\le 50 \delta$) to 
each other along a short (compared to $U$) segment. 
%It follows that $u$  is hyperbolic.  

We call the image of a piece $P_i$ by $u^p (p\in \Bbb Z)$  also 
a piece. 
%The piece $u^p(P_i)$ is {\em contained} in the axes $u^p(\alpha_i)$. 
The axis $\hat \alpha$ is a sequence of 
pieces (with short gaps in between).
%Likewise, we denote the axis that contains $P_i$
%by $\alpha_i$ for $i \in \Bbb Z$.

By construction, $x \in \hat \alpha$, but maybe 
$\hat \alpha$ is not exactly a $10\delta$-axis for $u$, so 
we take a $10\delta$-axis, $\alpha$.
One can check that $d(x,\alpha) \le 20\delta$.
(This is the reason we put $g^k$ at the end of $u$.
Without $g^k$ at the end, $u$ is still hyperbolic.)
Also, $\hat \alpha$ and $\alpha$ stay close to 
each other in the sense that most part of each piece in $\hat \alpha$, except
for some short parts near the two ends, is in the $20\delta$-neighborhood of $\alpha$.

(2)
If two axes $\beta,\beta'$ are parallel, we write
$\beta \sim \beta'$.
Since $s \not\in E(g)$, $\alpha_1
\not\sim \alpha_2$. Also, $\alpha_2 \not\sim \alpha_3$, $\alpha_3 \not\sim \alpha_4$, and $\alpha_4 \not\sim \alpha_5$.

We first argue that $E(u)$ maps to $\Bbb Z=C$ in the exact
sequence after Proposition \ref{prop.elementary}.
Suppose $h \in E(u)$. By definition, $\alpha \sim h(\alpha)$. 
It suffices to show that $h$ does not flip
the direction of $\alpha$. 
For that, we examine the sequence of the lengths of the pieces
on $\alpha$ and on $h(\alpha)$. 
Here, we say that $h(P)$ is a piece if $P$ is a piece in $\alpha$. 

\begin{figure}[htbp]
\hspace*{-3.3cm}     
\begin{center}
                                                      
   \includegraphics[scale=0.4]{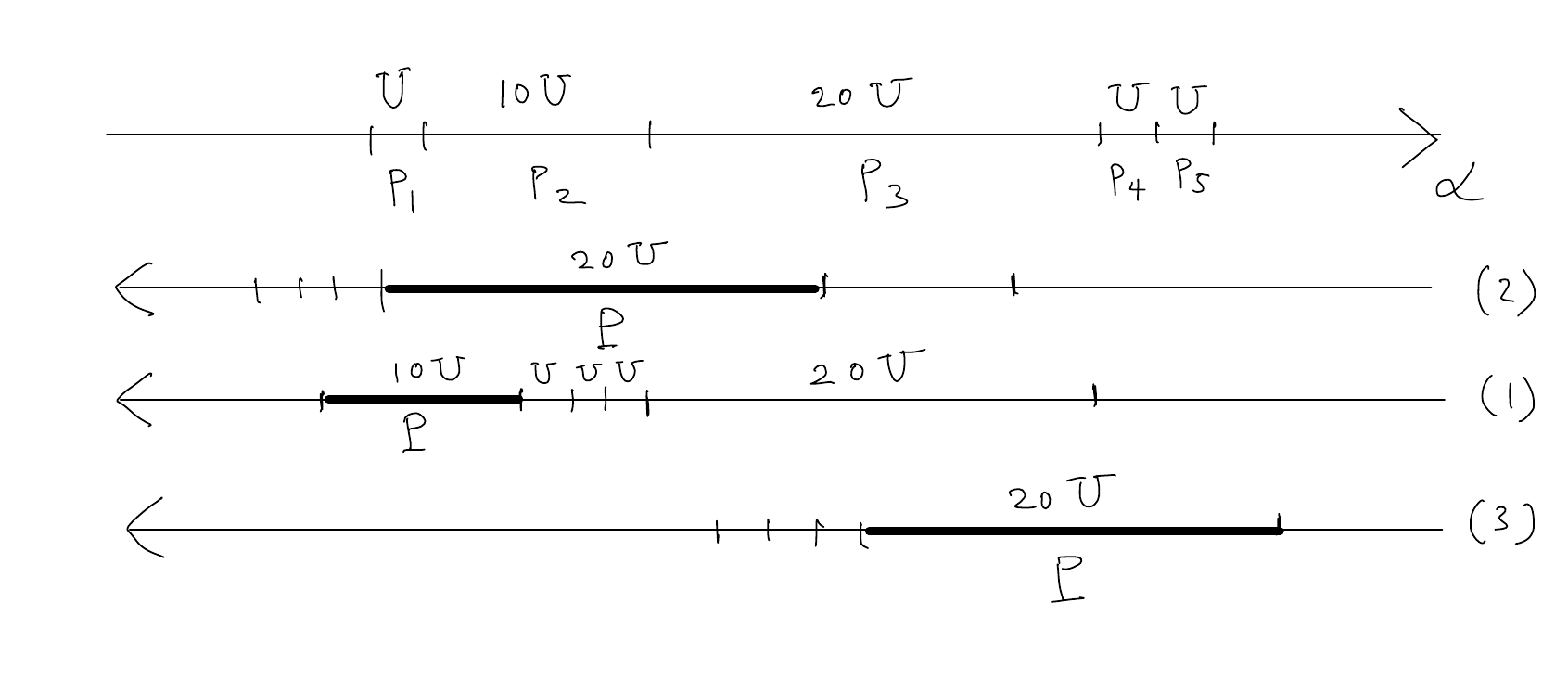}%
   \end{center}
\caption{The direction of $\alpha$ is to the right, and the direction of $h(\alpha)$ is to the left. The figure indicates the three positions of $P$.}
 \label{fig.noflip}
\end{figure}

The sequence for the part $P_1, \cdots, P_5$
is $U,10U,20U,U,U$, so that on $\alpha$, the sequence is (from the left to the right in Figure \ref{fig.noflip}):
$$\cdots; U,10U,20U,U,U; \,\, U,10U,20U,U,U; \cdots$$

Now, if the direction of $h(\alpha)$ was opposite to $\alpha$, then the sequence 
on $h(\alpha)$ would be (from the left to the right in the figure):
 $$\cdots; U,U, 20U, 10U,U;\, \,  U,U,20U,10U,U; \cdots$$

From those two sequences we observe that 
one of the followings must hold:
\begin{enumerate}
\item
there is a piece $P$ on $h(\alpha)$ that 
 intersects both $N_{50\delta}(P_1)$ and $
N_{50\delta}(P_2)$ at least $\frac{U}{2}$ in 
diameter, or
\item
there is a piece $P$ on $h(\alpha)$ that 
 intersects both $N_{50\delta}(P_2)$ and $
N_{50\delta}(P_3)$ at least $\frac{U}{2}$ in 
diameter, or
\item
there is a piece $P$ on $h(\alpha)$ that 
 intersects both $N_{50\delta}(P_3)$ and $
N_{50\delta}(P_4)$ at least $\frac{U}{2}$ in 
diameter.
\end{enumerate}
In the above, $P$ has length (approximately) $5U$ or $10U$.

Let $\beta$ be the axis that contains the piece $P$.
Then, since $\frac{U}{2}$ is at least
$5(2\D L(g)+ 100\delta)$, Lemma \ref{2.1}(2)
 would imply that either 
 (1) $\alpha_1 \sim \beta \sim \alpha_2$, 
 (2) $\alpha_2 \sim \beta \sim \alpha_3$, 
 or (3) $\alpha_3 \sim \beta \sim \alpha_4$, respectively. 
 In either case, we obtain a contradiction 
since $\alpha_1 \not\sim \alpha_2$, $\alpha_2 \not\sim \alpha_3$ and $\alpha_3 \not\sim \alpha_4$.
We showed that $h$ does not flip $\alpha$.
%It implies that $C=\Bbb Z$.

By the same reason, ie, the constrain from 
the combinatorics, $h(x)$ must
be close to $u^p(x)$ for some $p\in \Bbb Z$, namely, the distance is 
at most $2\D L(g)+ 100\delta$.
This implies that $\alpha_1 \sim (hu^{-p})(\alpha_1)$
and also, $\alpha_2 \sim (hu^{-p})(\alpha_2)$.
It then follows that 
$$hu^{-p} \in E(g), hu^{-p} \in E(g^nsgs^{-1}g^{-n}),$$
which implies $g^k(hu^{-p})g^{-k} \in (E(g) \cap E(sgs^{-1}))$. But the right hand side is a finite group, so that 
$hu^{-p}$ has finite order, therefore 
$hu^{-p} \in F(g)$.
And, of course, $hu^{-p} \in F(u)$
%Also, the conlusion  holds for all $h^n(x), n \in \Bbb Z$,
%and moreover, 
%with same $u^p(x)$, since the above bound
%is at least 10 times smaller than $\lambda(g^n)=U$.
%It implies that $hu^{-p}  \in F(u)$,
namely, there is $f \in F(u)$ with $h=fu^p$.

(3). 
This is a consequence of the $D'$-uniform WPD, and is mentioned in  the paragraphs
following Proposition \ref{prop.elementary}.

(4). Since $g\in S^2$ and $s\in S$, we have $u \in S^{14k+4}$.
We are done since $m = 66k+4$.
%
%(5) This holds for all non-trivial elements in the free group $F$.
\qed

\begin{prop}[Lower bound of growth]\label{bound.generators}

Suppose  $G$ acts on a $\delta$-hyperbolic space $X$.
Let $D(\epsilon)$ be a function for WPD.
Assume that there exists a constant $M$ such that for any 
finite generating set $S$ of $G$, the set $S^M$ contains a hyperbolic 
element $h$ that is $D$-WPD. 

Set
$$A=\frac{1}{79220M (D(200\delta))^3}>0.$$
 Then
%  there exists a constant $A>0$
%such that 
for any finite generating set $S$ of $G$, we have
$e(G,S) \ge A|S|^A$.

\end{prop}
This result is a generalization of the result on hyperbolic groups by \cite{AL}.
We adapt their argument to our setting, which is straightforward.

\proof
Set $\D=D(200\delta)$. 
Since the element $h \in S^M$ is $D$-WPD, we have
$\lambda(h) \ge \frac{50\delta}{\D}$
by Lemma \ref{lower.bound} (1).
It implies that $\lambda(h^\D) \ge 50 \delta$.
Since $h^\D \in S^{M\D}$, we have $L(S^{M\D}) \ge 50 \delta$.

%Let $T=\frac{50\delta}{D}>0$ be a constant from that gives a lower bound on
%$\lambda(h)$  of a hyperbolic element $h$ in $G$.
%Let $h \in S^M$ be a hyperbolic element. Then $\lambda(h) \ge T$. Since $D=\frac{50\delta}{T}$, we have
%$\lambda(h^D) \ge 50 \delta$.
%But $h \in S^{MD}$, we have $L(S^{MD}) \ge 50 \delta$.

Lemma \ref{2.3} applies to $S^{M\D}$ since $L(S^{M\D}) \ge 50 \delta$.
By Lemma \ref{2.3} applied to $S^{M\D}$ with $x$ such that 
$L(S^{M\D})=L(S^{M\D},x)$, there is $g \in S^{2M\D}$ such that 
$$L(S^{M\D}) - 8 \delta \le |x-g(x)|$$
and the $10\delta$-axis $\gamma$ of $g$ satisfies $d(x,\gamma) \le 20 \delta$.

Then Lemma \ref{2.4} applies to $S^{M\D}$ and $g \in S^{2M\D}$.
By Lemma \ref{2.4} applied to $S^{M\D}$ and $g$, there exists
$s \in S^{M\D}$ such that $\<g^k, sg^ks^{-1}\>=F$ and $F$ contains $u$ that 
is primitive such that 
$$|F(u)| \le \D; \, u \in S^{M\D m}; \,  L(u) \ge 10L(S^{M\D}),$$
where $k=60\D, m=3960\D +4$, and that 
$u$ is $D'$-WPD, where $D'(100\delta)=D(200\delta)$.

We note that the element $s$ can be chosen from $S$ since $G$
is not virtually cyclic. This can be easily seen in the proof of Lemma \ref{2.4}
since we only need to choose $s$ such that $s \not\in E(g)$. 
See Remark \ref{2.5}.

Now, take a maximal subset $W \subset S$ such that  
any two distinct elements $w,v \in W$ are in different 
$F(u)$-(right) cosets.
Then, $|W| \ge \frac{|S|}{\D} $. 

Let $\alpha$ be a $10\delta$-axis of $u$
with $d(x,\alpha) \le 50\delta$.
\\
{\it Claim 1}.
For distinct $v,w \in W$, $v\alpha, w\alpha$ 
are not parallel. Indeed, if they were parallel, then $w^{-1}v \in E(u)$.
Moreover, we will argue $w^{-1}v \in F(u)$, which is a contradiction since $w,v$ are in distinct $F(u)$-cosets. 
The reason for $w^{-1}v \in F(u)$ is that 
 since $w^{-1}v \in S^2 \subset S^{M\D}$, 
we have 
$$L(w^{-1}v) \le |w^{-1}v(x)-x| \le L(S^{M\D}).$$
But on the other hand, since $u$ is primitive,
if  $w^{-1}v \not\in F(u)$, then 
$$L(w^{-1}v) \ge L(u) -100 \delta \ge 6 L(S^{M\D}).$$
The last inequality is from 
 $L(u) \ge 10L(S^{M\D})$.
Those two estimates contradict. We showed the claim.

It implies that for distinct $v,w \in W$, the intersection of 
$v\alpha$ and the $50\delta$-neighborhood of $w\alpha$ is 
bounded by 
$$2 D'(100\delta)  L(u) +100\delta$$
by Lemma\ref{2.1} (2) since $u$ is $D'$-WPD.
Remember that $D'(100\delta)=D(200\delta)=\D$.
So, this bound is $2\D L(u) + 100\delta$.

Set $U=u^{20\D}$. 
Then $U \in S^{20M\D^2 m}$
and $L(U) \ge 19\D L(u)$ (maybe not quite $20\D L(u)$).
 $\alpha$ is an axis for $U$ as well. 
Set 
$$B=\{wUw^{-1}|w\in W\}.$$
\noindent
{\it Claim 2}. We have $|B|=|W|$ and $B$  freely generates a free group of rank $|W|$.

This is because for any $w \in W$, we have $|x-w(x)| \le L(S^{M\D})$
since $W \subset S^{M\D}$, 
which means the axis of $wUw^{-1}$, $w\alpha$, is close to $x$.
To be precise, close means that the distance is much smaller than 
$L(u)$ since $L(u) \ge 10 L(S^{M\D})$.
Also, for the  axes $w\alpha$ and $v\alpha$ of any distinct $v,w \in W$,  
the intersection of one with the $50\delta$-neighborhood of 
the other is $9$ times shorter than 
$L(U)$ since $L(U) \ge 19\D L(u)$.
In this setting, the usual ping-pong argument shows the claim.

Since $w \in S$, we have $B=\{wUw^{-1}|w\in W\} \subset S^{20M\D^2 m+2}$.
It follows that for any $n \in \Bbb N$, 
$$|S^{(20M\D^2 m+2)n}| \ge |B^n| \ge |B|^n = |W|^n \ge \frac{ |S|^n}{\D^n}.$$
It implies 
$$e(G,S) \ge \D^{-\frac{1}{20M\D^2 m+2}} |S|^{\frac{1}{20M\D^2 m+2}}.$$

%\frac{\log|S|-\log D}{MDm20D+2}$$
%Set $A=\frac{\log D}{20MD^2m+2}$, which is 
%$=\frac{\log D}{16880MD^3+2}.$
Since 
$$\min\left\{\D^{-\frac{1}{20M\D^2 m+2}}, \frac{1}{20M\D^2 m+2}\right\}=\frac{1}{20M\D^2 m+2},$$
which is at least $\frac{1}{79220M\D^3}$,
since $m=3960\D+4$.
Setting $$A=\frac{1}{79220M\D^3},$$
we have $e(G,S) \ge A|S|^A$.
This is a desired conclusion since $\D=D(200\delta)$.
\qed

\begin{example}\label{ex.bound}
Proposition \ref{bound.generators} applies to the following examples:
\begin{enumerate}
\item
Non-elementary hyperbolic groups (the original case in \cite{AL}).

\item
The mapping class groups of a compact orientable surface
$\Sigma_{g,p}$ with $3g+p \ge  4$.
See Section \ref{section.mcg}, where the assumptions are checked. 
\item
A lattice in a simple Lie group of rank-1.
See Theorem \ref{lattice}.

\item
The fundamental group of a complete Riemannian 
manifold of finite volume whose sectional 
curvature is pinched by two negative constants. 
See Theorem \ref{lattice}.
\end{enumerate}

\end{example}

\section{Well-orderedness}\label{section.3}
\subsection{Main theorem}
We prove the following theorem.
Note that Theorem \ref{main} immediately follows
from this theorem 
combined with Lemma \ref{wpd.acyl}, since the lemma
says that an acylindrical action is uniformly WPD, so that every
hyperbolic element is WPD. 
\begin{thm}[Well-orderedness for uniform WPD actions]\label{main.proof}

Suppose $G$ acts on a $\delta$-hyperbolic space $X$, and $G$ is not virtually cyclic. 
Let $D(\epsilon)$ be a function for WPD.
Assume that there exists a constant $M$ such that for any 
finite generating set $S$ of $G$, the set $S^M$ contains a hyperbolic 
element on $X$ that is $D$-WPD.
Assume that $G$ is equationally 
Noetherian.
Then, 
 $\xi(G)$ is a well-ordered set.
\end{thm}

\proof 
We will prove that $\xi(G)$ does not contain a strictly decreasing convergent sequence.
To argue by contradiction, 
suppose that there exists a sequence of finite generating sets $\{S_n\}$,
such that $\{e(G,S_n)\}$ is a strictly decreasing sequence and $\lim_{n \to \infty} e(G,S_n)=d$, for some $d>1$.

By Proposition \ref{bound.generators},  we may assume
that the cardinality of the generating sets $|S_n|$ from the decreasing sequence is bounded, and by possibly passing to a subsequence we may assume that the cardinality
of the generating sets is fixed, $|S_n|=\ell$.
 
Let $S_n=\{x_1^{(n)}, \cdots, x_\ell^{(n)}\}$.
Let $F$ be the free group of rank $\ell$ with a free generating set:
$S=\{s_1, \ldots, s_\ell\}$.
For each index $n$, we define a map: $f_n:F  \to G$, by setting: $f_n(s_i)=x_i^{(n)}$.
Since $S_n$ are generating sets, the map $f_n$ is an epimorphism for every $n$.
Note that $e(G,S_n)=e(G,f_n(S))$.
%We denote $F_\ell$ as $F$.

Since $F$ is countable, the sequence $\{f_n:F\to G\}$ subconverges to 
a surjective homomorphism $\eta:F\to L$. $L$ is called a {\it limit group} over the group $G$.

By assumption, $G$ is equationally Noetherian.
By the general principle (Lemma \ref{basic}),
there exists an epimorphism $h_n: L \to G$ such that 
by passing to a subsequence we may assume that all the homomorphisms 
$\{f_n\}$ factor through the limit epimorphism: $\eta: F \to L$, ie, 
 $f_n=h_n \circ \eta$. 
  
    $$
  \xymatrix{
  (F, S) \ar[d]_\eta \ar[dr]^{f_n}& \\
  (L, \eta(S)) \ar[r]_{h_n} & (G, f_n(S))
  }
  $$
 
  Notice that since $f_n=h_n \circ \eta$ for every  index $n$,
we have $e(G, f_n(S)) \leq e(L, \eta(S))$.
  
We will show the following key result:

\begin{prop}\label{1.2}
Suppose $G$ satisfies the assumption in Theorem \ref{main.proof}.
Let $(L,\eta(S))$ be the limit group over $G$ of a sequence $f_n:(F,S) \to (G, f_n(S))$,
where $F$ is a free group with a free generating set $S$ and 
$f_n(S)$ are generating sets of $G$. Then 
$\lim_{n \to \infty} e(G,f_n(S))=e(L,\eta(S))$.

\end{prop}

We postpone proving this proposition until the next section
and finish the proof of the theorem.

 We assumed that the sequence $\{e(G, f_n(S))\}$ is strictly decreasing, hence, it can not converge to an
upper bound of the sequence, $e(L,\eta(S))$.
But the proposition says that it must converge to $e(L,\eta(S))$, 
a contradiction. 
Theorem \ref{main.proof} is proved. 
\qed

\section{Continuity of the growth rate}\label{section.4}

We prove Proposition \ref{1.2}.
We already know that  $e(L,\eta(S))  \ge e(G, f_n(S))$ from the existence of the surjections $h_n$, which follows from that $G$ is equationally Noetherian. 

It suffices to show that given $\epsilon>0$, for a large enough $n$,
$$\log e(L,\eta(S)) -\epsilon\le \log e(G,f_n(S)).$$

The strategy of the proof of this is same as \cite{FS}.
We note that from now on, we do not use that $G$ is equationally Noetherian
in the proof. 

Since the proof is long and complicated, we first informally describe
the idea, which already appeared in \cite{FS}.
We want to show $e(G,f_n(S))$ is almost equal to $e(L, \eta(S))$ for a
large enough $n$. First of all, if we take a large enough $r$, 
then $B_r(L,\eta(S))$ contains elements roughly as many as
$e(L,\eta(S))^r$ by the definition. Fix such $r$. Then if we take 
$n$ large enough, $B_r(L,\eta(S))$ and $B_r(G,f_n(S))$
are idential via the map $h_n$ since $L$ is a limit group.
But it does not mean that $e(G,f_n(S))$ is almost equal to 
$e(L,\eta(S))$ since the growth of the balls in $(G,f_n(S))$
may decay if we take the radius larger. 
But it turns out that if we take $r$ large enough, then 
roughly speaking, the growth of the ball  of radius $r$ in $(G,f_n(S))$
is almost equal to $e(G,f_n(S))$. This is due to the well-known
``local-to-global'' principle in hyperbolic groups, and 
it is implemented by inserting ``separators'' in our argument. 
The threshold for the radius is 
given by $m$ in the proof (see Section \ref{section.4.5}).

We  explain the idea more in detail. 
By the definition of the growth rate, we have for all $r$,
$$\log e(L,\eta(S)) \le \frac{1}{r}\log   |B_r(L,\eta(S))|.$$
This is because the sequence $\{\log |B_r(L,\eta(S))|\}$
is sub-additive. 

Fix $r$ (we will choose $r$
sufficiently large in the argument we will give later).
Then choose $n$ large enough such that 
$h_n:B_r(L,\eta(S)) \to B_r(G,f_n(S))$ is injective. 
The following map is naturally induced from $h_n$  for each $q\in \Bbb N$:
$$B_r(L,\eta(S))^q \to B_{qr}(G,f_n(S)) \subset G$$
by mapping $(w_1, \cdots, w_q)$ to $h_n(w_1 \cdots w_q)$.
If there is an $r$ such that 
this map is injective for all $q$, then an easy computation would show the desired inequality for the $n$ by letting $q\to \infty$. 

But of course this map is not injective in general. For example,
the concatenation of $h_n(w_1), h_n(w_2), \cdots, h_n(w_q)$
may have lots of backtracks at the concatenation points. 
As a remedy, we insert elements $u_i$ of bounded length, called {\it separators}, 
and define a new map
sending $(w_1, \cdots, w_q)$ to $h_n(w_1u_1 w_2 u_2  \cdots w_qu_q)$.
This map is denoted by $\Phi_n$.
The separators are constructed in Lemma \ref{4.2}. 
We arrange that  the concatenation of elements, after we insert $h_n(u_i)$'s, is a uniform quasi-geodesic in $G$. This follows from
 that $G$ is a hyperbolic group. 
 
But it is still not the case that the map $\Phi_n:B_r(L,\eta(S))^q \to B_{q(r+b)}(G,f_n(S))$ is injective,
where $b$ is the bound of the length of the separators.
What we actually show is
that $\Phi_n$  is injective if we restrict it
to the $q$-tuples in some fixed portion of $B_r(L,\eta(S))$, 
(see Lemma \ref{feasible}), which is enough for our purpose.
This part is very technical. 
To argue that $\Phi_n$ is injective on the certain fixed portion, 
we use the action of $G$ on $X$, ie, we map a tuple by $\Phi_n$ to $G$, 
then let it act on $X$. Then we analyze the orbit of a base point in $X$.

%Notice that on $B_m(L,\eta(S))$, this action is faithful for large
%enough $n$ depending on $m$.

\subsection{Separators}\label{section.separator}
We review the setting. 
$X$ is a $\delta$-hyperbolic space, and $G$ acts on it.
%The action is $D$-uniformly WPD. 
We assume $\delta\ge1$.
%and $D\ge 50 \delta$ (Proposition \ref{bound.generators}).
For each $n$, $S_n$ is a finite generating set of $G$
such that $S_n^M$ contains a hyperbolic element that is 
$D$-WPD. 
Using the homomorphism $h_n:L \to G$, we let $L$ act on $X$.
We first construct separators as elements in $G$ then 
pull them back to $L$ by $h_n$.

Let $g \in S_n^M$ be a hyperbolic element
that is $D$-WPD. 
Set $$\D=D(200\delta).$$
Then  $\lambda(g)\ge \frac{50\delta}{\D}$
by Lemma \ref{lower.bound}(1).
It implies
$$100\delta \le  \lambda(g^{2\D}) \le L(S_n^{2\D M}).$$ 
%The constants $\delta, D, C, N$ only depends on $\Sigma$.

Fix $n$. 
Let $y_n \in X$ be a point where 
$L(S_n^{2\D M})$ is achieved. 
Put $$\Delta_n=100\delta +4\D L(S_n^{2\D M}).$$
%We sometimes suppress $n$ and write $\Delta$.

We define a germ w.r.t. the constant $\Delta_n$.
Recall that given  three points $x,y,z \in X$, the {\em Gromov product},
$(y,z)_x$, is defined as follows:
$$(y,z)_x = \frac{|x-y|+|x-z|-|y-z|}{2}.$$

\begin{definition}[germs, equivalent and opposite germs]
Let $[x,y]$ be a (directed) geodesic segment in $X$. Suppose that 
$|x-y|\ge 10\Delta_n$. Then, the initial segment of $[x,y]$
of length $10\Delta_n$ is called the {\it germ}
of $[x,y]$ at $x$, denoted by $germ([x,y])$.
If $|x-y| < 10\Delta_n$, then we define the germ to be empty. 

%The initial part of the (directed) segment $[y_n,g(y_n)]$ of length $100\delta$
%is called the {\it germ} at $y_n$, and the 
%initial part of the segment $[g(y_n), y_n]$ of length $100\delta$ is called the germ at $g(y_n)$.
%If $\lambda(g) < 100\delta$ then the germ is empty.
We say two non-empty germs, $[x,y], [x,z]$,  at a common point $x$, 
are {\it equivalent} if 
$(y,z)_x \ge 4 \Delta_n$, 
%$|y-z| \le |x-y| + |x-z| - 8 \Delta_n$, 
and {\it opposite}
if 
$(y,z)_x \le 2\Delta_n$.
%$|y-z| \ge |x-y| + |x-z| - 4 \Delta_n$.
\end{definition}

%Roughly speaking, two germs are equivalent
%if they are parallel (ie, stay close) along the initial parts
%of length at least $4\Delta_n$, and opposite if 
%the initial parallel parts are at most
%$2\Delta_n$ in length. 

We sometimes call the germ of $[y,x]$ at $y$ 
as the germ of $[x,y]$ at $y$.
If $\gamma$ is the germ  of $[x,y]$ at $x$, then for $g\in G$, the segment 
$g(\gamma)$ is the germ of $[g(x),g(y)]$ at $g(x)$. 
For $g \in G$, we consider the germ of $[y_n,g(y_n)]$ and
call it {\em the germ of $g$ at $y_n$}, and write $\germ(g)$.

%We say two segments $\alpha,\beta$ {\it virtually intersect}
%if their $10\delta$-neighborhoods intersect, and the intersection
%is called the  {\it virtual 
%intersection}.

Recall  from Lemma \ref{2.4} that $k=60\D, m=3964\D +4$.
%To avoid confusion, we temporarily  use the letter $t$ as the subscription of  $S_n$ instead of $S_n$.

We consider germs w.r.t. $\Delta_n$.

\begin{lemma}[The constant $b$ and separators. c.f. Lemma 2.4 \cite{FS}]\label{4.2}
There exists a constant $b$ with the following property,
where $b$ depends only on $\delta, M$ and $D(\epsilon)$. 
For every $n$, 
there exist primitive, hyperbolic  elements $u_1, u_2, u_3, u_4\in S_n^{b}$; 
and mutually opposite germs  $\gamma_1, \gamma_2,\gamma_3,\gamma_4$ of some elements  in $S_n^{2\D M}$
at $y_n$ that satisfy:
\begin{enumerate}
\item[(i)]  
{\em The germs are all at $y_n$ in the following. }
\begin{itemize}
\item
The germ of $[y_n,u_1(y_n)]$ is equivalent to $\gamma_1$ .
The germ of $[y_n,u_1^{-1}(y_n)]$   is equivalent to $\gamma_3$.

\item
The germ of $[y_n,u_2(y_n)]$ is equivalent to $\gamma_1$.
The germ of $[y_n,u_2^{-1}(y_n)]$ is equivalent  to $\gamma_4$.
\item
The germ of 
$[y_n,u_3(y_n)]$ is equivalent to $\gamma_2$.
The germ of  $[y_n,u_3^{-1}(y_n)]$ is equivalent to $\gamma_3$.
\item
 The germ of 
$[y_n,u_4(y_n)]$ is equivalent to $\gamma_2$.
The germ of  $[y_n,u_4^{-1} (y_n)]$ is equivalent to $\gamma_4$. 
\end{itemize}

%The germ of 
%$[y_n,u_{i,j}(y_n)]$ at $y_n$ is equivalent to the germ $\gamma_i$ and its germ at $u_{i,j}(y_n)$ is 
%equivalent to $u_{i,j}(\gamma_j)$.

\item[(ii)] 
For all $i$, 
 $$\lambda(u_{i}) \ge 100 \Delta_n, $$ 
and 
the distance from $y_n$ to the $10\delta$-axis 
of $u_{i}$ is at most $\Delta_n$.

\item[(iii)] For every $w \in G$, and all $i,j$ (possibly $i=j$), 
 if the
$20\delta$-neighborhood of  $[y_n,u_{i}(y_n)]$ intersects the segment
$[w(y_n),wu_{j}(y_n)]$, then the diameter of the intersection is bounded by: $\frac {1} {10} d(y_n,u_{i}(y_n))$ and $\frac {1} {10} d(y_n,u_{j}(y_n))$.
If $i=j$, we assume in addition  that $w \not\in F(u_i)$,
where $F(u_i)$ is the finite normal subgroup in $E(u_i)$.

\end{enumerate}
\end{lemma}

\begin{remark}\label{remark.cancellation}
(1) 
Regarding Lemma \ref{4.2} (iii), if $w\in F(u_i)$ then 
the Hausdorff-distance between $[y_n,u_i(y_n)]$
and $[w(y_n), wu_i(y_n)]$ is 
at most $\Delta_n + 100 \delta$
since $w$ moves any point on a $10\delta$-axis
of $u_i$ by at most $50\delta$ since otherwise
$w$ would be hyperbolic. 
\\
(2) As we will see in the proof it suffices to take  $b$ to be
\begin{equation}\label{length.separator}
b=343440 M (D(200\delta))^2  + 22.
\end{equation}

\end{remark}

The elements $u_i$ are called {\it separators} and 
the property (iii) is called the {\it small cancellation property} 
of separators. 

Note that we have $L(S_n^{2M\D})
\ge L(S_n)$.

\proof
Since $L(S_n^{2M\D}) \ge 100 \delta$
and $y_n$ is a point where $L(S_n^{2M\D})$ is achieved,
by Lemma \ref{2.3} (also see the properties at the bullets after the lemma), there is a hyperbolic element $g \in S_n^{4M\D}$
such that 
$$2L(S_n^{2M\D}) \ge L(g) \ge L(S_n^{2M\D})-24\delta \ge \frac{1}{2}L(S_n^{2M\D})
$$
and the $10\delta$-axis of $g$, $\gamma$, is at distance at most $10\delta$
from $y_n$.

By Lemma \ref{2.4} applied to $S_n^{2M\D}$ and $g$ and the $10\delta$-axis $\gamma$,  there exists $s\in S_n$
(see the remark \ref{2.5}) such that 
$g^k$ and $sg^ks^{-1}$ are independent hyperbolic elements,
which freely generate a free group $F$ whose non-trivial element,
$h$, satisfies
$$\lambda(h) \ge 10(2\D L(g)+100\delta) \ge 5 \Delta_n.$$

Recall that $k=60 D(200\delta)=60\D$. 
Note that $g^k, sg^ks^{-1} \in S_n^{4M \D k +2}$. 
The distance from $y_n$ to $s\gamma$
is at most $40\delta+L(S_n^{2M\D})$.
Also, the intersection of $\gamma$ and 
the $50\delta$-neighborhood of $s\gamma$ is at most 
$2\D L(g) +100\delta$ in length,
which is $\le \Delta_n$, by Lemma \ref{2.1} since $\gamma$ and $s\gamma$ are not parallel.

Consider the following four germs at $y_n$ w.r.t. the constant $\Delta_n$:
$$\gamma_1= \germ(g^k), \gamma_2=\germ(sg^ks^{-1}), 
\gamma_3= \germ(g^{-k}), \gamma_4= \germ(sg^{-k}s^{-1}).$$
Note that any two of them are mutually opposite.

%Therefore,  $S_n^{200N+2}$ contains
%two independent pseudo-Anosov elements, $g,h$, which generate a free group $F(g,h)$.
%Moreover, by the acylindricity, all non-trivial elements, $k$,  in $F(g,h)$ satisfies
% (\ref{large.trans}): $50 L(S_n^N) \le \lambda(k)$.
% Note that $E(g) \cap E(h)$ is a finite group. 

Then, there exist separators $u_i \in F$
such that $u_i \in S_n^b$, 
where
$$b=1431\cdot 4\D Mk  + 22=343440\D^2 M+22.$$

For example, set $w=g^k, z=sg^ks^{-1}$ and:
\begin{align*}
u_1&= wzw^2zw^3z \cdots w^{19}zw^{20}, \\
u_2&= w^{21} z w^{22} z \cdots w^{39}zw^{40}z, \\
u_3&= zw^{41}zw^{42}zw^{43}z  \cdots w^{59}zw^{60}, \\
u_4&= zw^{61} z w^{62} zw^{63}z  \cdots w^{79}zw^{80}z.
\end{align*}
We compute that they are in $S_n^b$ since $w \in S_n^{4\D Mk}, z \in S_n^{4\D Mk+2}$. It will be important that this number $b$ does
not depend on $S_n$.  (See the proof of Proposition \ref{1.2}.)

Because of the combinatorial reason, they are primitive hyperbolic elements. 
The argument is similar to the one we used to showed that $u$ is primitive
hyperbolic, using $\gamma$ and $s\gamma$ are
not parallel, in the proof of Lemma \ref{2.4}(2). We omit it. 

(i) By definition of $u_i$, 
$$\germ[y_n,u_1(y_n)]=\germ[y_n,u_2(y_n)]=\germ(w)=\gamma_1;$$
$$\germ[y_n,u_3(y_n)]=\germ[y_n,u_4(y_n)]=\germ(z)=\gamma_2;$$
$$u_1^{-1}\germ[u_1(y_n),y_n]=u_3^{-1} \germ[u_3(y_n),y_n]=\gamma_3;$$
$$u_2^{-1}\germ[u_2(y_n),y_n]=u_4^{-1}\germ[u_4(y_n),y_n]=\gamma_4.$$

Here, we mean that the four germs in each 
line are same or equivalent to each other.

(ii). The estimate for $\lambda(u_i)$ is straightforward  from the definitions of $u_i$.
The claim on the axes are also shown similarly to the case of $u$
in Lemma \ref{2.4}(1) and we omit it. 
We remark that the distance estimate on the axes differ since it  comes from the 
fact that some of the separators start or end in $w$, so that 
distance becomes larger. 

(iii) follows from the definition of $u_i$, namely, the combinatorial 
structure of the words. The argument is similar to show
that $u$ in Lemma \ref{2.4} is primitive (Lemma \ref{2.4}(2)), so we will be brief.
Suppose $i\not=j$, and the intersection was longer.
Then, because of  the combinatorial structure of the words $u_i$ and $u_j$,
it follows the axes $\gamma$ and $s\gamma$ would be parallel, by Lemma \ref{2.1},
which is impossible. We are done. 
Suppose $i=j$, and the intersection was longer. Then since $w \not\in F(u_i)$, by the same reason,
$\gamma$ and $s\gamma$ would be parallel, a contradiction. 
So we are done in this case too. 
\qed

%\item
Remember $u_i \in G$ depend on the index $n$ of $S_n$, so 
let's write them as $u_i(n)$. 
Now, since $h_n$ is surjective, let $\hat u_{i}(n) \in L$ be an element with 
$h_n(\hat u_{i}(n))=u_{i}(n)$. Note that the word length of $\hat u_{i}(n)$
in terms of $\eta(S)$ is also bounded by $b$.
The elements $\hat u_i(n) \in L$ are also called the {\it separators}
for $h_n$. 

In the following we may just write $u$ (instead of $\hat u$) to denote a separator for $h_n$
to simplify the notation. 
We note that $|F(h_n(u))| \le \D$ for any separator $u$ for any $h_n$
since $h_n(u)$ is primitive (see the comment after Proposition \ref{prop.elementary}).

\subsection{Forbidden elements}\label{section.forbidden}
Given $m$,  choose and fix $n$ large enough
such that
the map $h_n$ is injective on $B_{2m}(L,\eta(S))$.
This is possible since this set is finite in the limit group $L$.
Remember
$$\Delta_n=100\delta +4\D L(S_n^{2\D M}).$$

Given $w \in B_m(L,\eta(S))$ we choose one of the separators, $u\in L$,
for $h_n$ s.t.
the germ of $[y_n,h_n(w)(y_n)]$ at $h_n(w)(y_n)$
and the germ of 
$[h_n(w)(y_n), h_n(wu)(y_n)]$ at $h_n(w)(y_n)$
are opposite. (Geometrically, it implies that the 
concatenation $[y_n,h_n(w)(y_n)] \cup [h_n(w)(y_n), h_n(wu)(y_n)]$ 
is almost a geodesic.)
Such a separator $u$ exists by the 
property (i) in Lemma \ref{4.2}.
Here, if $|y_n-h_n(w)(y_n)| < 10 \Delta_n$, then 
we choose any separator $u$.
We say $u$ is {\it admissible} for $w$.
Note that 
$$|y_n-h_n(wu)(y_n)| \ge |y_n-h_n(w)(y_n)| + |h_n(w)(y_n) -h_n(wu)(y_n)| - 4\Delta_n.$$
Also, the Hausdorff distance between
$[y_n,h_n(wu)(y_n)]$ and $[y_n,h_n(w)(y_n)]\cup [h_n(w)(y_n),h_n(wu)(y_n)]$
is at most $2\Delta_n + 10\delta$.

Moreover, given $w,w' \in B_m$ then we can choose a separator 
$u$ such that 
$u$ is admissible for $w$ and $u^{-1}$ is admissible 
for $w'^{-1}$. We say $u$ is {\it admissible} for $w,w'$.
Note that 
\begin{align*}
|y_n-h_n(wuw')(y_n)| \ge & |y_n-h_n(w)(y_n)| + |h_n(w)(y_n) -h_n(wu)(y_n)|
\\
&+|h_n(wu)(y_n) -h_n(wuw')(y_n)| - 8 \Delta_n.
\end{align*}

As before, 
the Hausdorff distance between
$[y_n,h_n(wuw')(y_n)]$ and $[y_n,h_n(w)(y_n)]\cup [h_n(w)(y_n),h_n(wu)(y_n)] \cup [h_n(wu)(y_n),h_n(wuw')(y_n)]$
is at most $2\Delta_n +10\delta$.

For $q>0$, we define a map
$$\Phi_n:B_m(L,\eta(S))^q \to B_{q(m+b)}(G,f_n(S)) \subset G$$
by sending $(w_1, \cdots, w_q)$
to $h_n(w_1 u_1 \cdots w_q u_q)$,
where $u_i \in L$ are separators we choose that 
are admissible for $w_i, w_{i+1}$
for $h_n$.
Remember that $u_i \in B_b(L,\eta(S))$.

It is easy to show that $\Phi_n$  maps an element 
$\not=(1, \cdots, 1)$
to a non-trivial element in $G$ using the property
of the separators, but what we need
is more. 
We will argue that on a large portion, called
the set of  ``feasible elements'', 
 $\Phi_n$ is injective 
by showing the image of the base point $y_n \in X$
by those elements are all distinct.

%Let $y_n \in \CC$ be a point where $L(S_n)$ is achieved.

For the given $m$, 
we will define {\it forbidden} elements in $B_m(L,\eta(S))$,
which depend on $n$. We are assuming that  $n$ is large enough so that 
$h_n$ is injective on $B_{2m}(L,\eta(S))$.
\begin{definition}[forbidden elements and tails]\label{def.forbidden}
Given $w \in B_m(L,\eta(S))$, if there exist $w' \in B_m(L,\eta(S))$ and 
a separator $u(n)$ which is admissible for $w$ such that 
\begin{equation}\label{forbidden}
|h_n(w')(y_n)-h_n(wu(n)) (y_n)| \le \frac{1}{5}|y_n-h_n(u(n))(y_n)|
\end{equation}
then we say $w$ 
 is {\it forbidden} w.r.t. $n$ (or in terms
of $h_n)$.
We call the segment $[h_n(w)(y_n), h_n(wu(n))(y_n)]$
a {\it tail} of $w$. 
The tail depends on the choice of $u$, so if we want to 
specify it, we say the tail of the pair $(w,u)$. 
\end{definition}

%\begin{definition}[parallel tails]
%Fix $n$.
Let  $w_1,w_2$ be two (forbidden) elements and $u_1,u_2$ admissible  separators , respectively.
If  $u_1=u_2$
and $h_n(w_1 ^{-1} w_2) \in F(h_n(u_1))$, then  the Hausdorff-distance between the two tails $[h_n(w_1)(y_n), h_n(w_1u)(y_n)]$ and $[h_n(w_2)(y_n), h_n(w_2u)(y_n)]$ is at most $\Delta_n + 100 \delta$.
This is an immediate consequence of Remark \ref{remark.cancellation}.
In this case, we say that the tails of $(w_1,u_1)$ and $(w_2,u_2)$ are {\it parallel}.

%\end{definition}

%We record some facts we use later. 
%\kf{better to define $F(u)$ using $h_n$}
%\begin{remark}\label{4.5}
%\kf{check this remark}
%(1) If two tails $(w_1,u), (w_2,u)$ are parallel, then
% by Lemma \ref{4.2} (ii)
%since the element $h_n(v_1^{-1}v_2) \in F(h_n(u_1))$ moves a point on 
%a $10\delta$-axis of $u_1$ at most $\le 50\delta$, since 
%otherwise, $h_n(u_1)$ would be hyperbolic. 
%\kf{check constants} 

On the other hand, if two tails $(w_1,u_1)$ and $(w_2,u_2)$ are not parallel, then 
the intersection of one of the tails
with the $20\delta$-neighborhood of the other tail is bounded
by the $\frac{1}{10}$ of the length of each tail. This is 
by the small cancellation property of the separators (Lemma \ref{4.2} (iii)),

We record one immediate consequence we use later. 

\begin{lemma}[Parallel tails]\label{lemma.tail}
Assume that $h_n$ is injective on 
$B_{2m}(L,\eta(S))$. 
Suppose $w\in B_m(L,\eta(S))$ is forbidden w.r.t. $n$, 
and $u$ is a separator admissible for $w$. 
Then there are at most $\D$ possibilities for $w_1 \in B_m$,
including $w=w_1$, such that $w_1$ is forbidden, $u$ is admissible for $w_1$,
 and 
the tails for $(w,u)$ and $(w_1,u)$ are parallel.
\end{lemma}

\proof

%$A$ does not depend
%on $n$ since $h_n(u)$ is a common element for each 
%$u$.
%\kf{we may move this later where we actually use it. Keep it as is}
%Also, given $(v_1,u)$ with $v_1 \in B_m(L,\eta(S))$ and $u$ admissible with $v_1$, there are at most $D$ possibilities
%for $v_2 \in B_m(L,\eta(S))$ such that 
%$(v_1,u), (v_2,u)$ are parallel. 
Since the two tails are parallel, we have $h_n(w_1^{-1} w) \in F(h_n(u))$
by the definition that two tails are parallel. 
Recall that (see the paragraphs after Proposition \ref{prop.elementary}.)
$$|F(h_n(u))| \le \D$$ for  all separators $u$ since $u$ is primitive.
Therefore, we have at most $\D$ possibilities for $h_n(w_1^{-1}w)$.
But since $h_n$ is injective on $B_{2m}(L,\eta(S))$,
 we have at most $\D$ possibilities  for $w_1^{-1}w
\in B_{2m}(L,\eta(S))$, and we are done. 
\qed

\subsection{Ratio of forbidden elements}

The proof of the following lemma occupies this subsection. 
%\kf{here, we use hn is inj on B2m}

\begin{lemma}\label{4.6}
%Given $m$, take $n$ large enough such that $h_n$ is injective on $B_{2m}(L,\eta(S))$. 
Assume that $h_n$ is injective on 
$B_{2m}(L,\eta(S))$. 
Consider the forbidden/non-forbidden elements in $B_m(L,\eta(S))$ w.r.t. $n$. Then
$$|\{\text{the forbidden elements}\} | \le \D |\{\text{the non-forbidden elements}\}|.$$
\end{lemma}

We denote $B_m(L,\eta(S))$ as $B_m$. 
The strategy of the proof is to show that if $w \in B_m$ is forbidden,
it will force some other elements in $B_m$ to be non-forbidden. 

\proof
The proof has three parts. 
\\
{\em Part I}.  We first construct  a subset $C(w)$ in $B_m$ and a tree like graph $T(w)$ in $X$.
The construction is inductive. In each step, a subset and a tree like graph grows, which will end in finite steps.

Suppose $w \in B_m$ is forbidden. We explain the inductive steps
to construct $C(w)$ and $T(w)$. 
\\
{\it Step 0}. Set $C_0(w)=w$ and $T_0(w)=[y_n, h_n(w)(y_n)]$.
\\
{\it Step 1}. 
Since $w$ is forbidden, there exist $w' \in B_m$
and a separator $u$ that is admissible for $w$ that satisfy (\ref{forbidden}).
%In the $L$-neighborhood of the tail, $[h_n(w)(y_n), h_n(wu)(y_n)]$, there are 
%10 points $h_n(wu^{(k)})(y_n), 1 \le k \le 10$, which 
%we call the neighboring points for the tail.
%We call the elements $wu^{(k)} \in B_{m+100N}$
Let $w'=s_1 \cdots s_r, s_i \in \eta(S), r\le m$ be a shortest representative
w.r.t. the word metric by $\eta(S)$. 
Then this defines a sequence of points in $X$, which we call
a path, $\gamma$, from $y_n$ to $h_n(w')(y_n)$
as follows:
$$y_n, h_n(s_1)(y_n), h_n(s_1s_2)(y_n), \cdots, h_n(s_1\cdots s_r)(y_n).$$
The distance between any two adjacent points on $\gamma$ is at most $L(S_n^{2MD})$ 
since it is achieved at $y_n$ and $h_n(s_i) \in S_n
\subset S_n^{2MD}$.

Consider the nearest point projection, denoted by $\pi$, from
a point $x\in \gamma$ to 
the tail $[h_n(w)(y_n), h_n(wu)(y_n)]$. The nearest
points may not be unique, but we choose one for $\pi(x)$.
Then the distance between the projection
of any two adjacent points on $\gamma$
is $\le L(S_n^{2MD}) + 100\delta \le \frac{\Delta_n}{4}$.
Note that $|\pi(y_n)-h_n(w)(y_n)| \le 2\Delta_n + 150 \delta$
since $u$ is admissible for $w$.
Also, 
$|\pi(h_n(w')(y_n)) - h_n(wu)(y_n)| \le \frac{1}{5}|h_n(w)(y_n)-h_n(wu)(y_n)|+10\delta$ by (\ref{forbidden}).

\begin{figure}[htbp]
\hspace*{-3.3cm}     
\begin{center}
                                                      
   \includegraphics[scale=0.4]{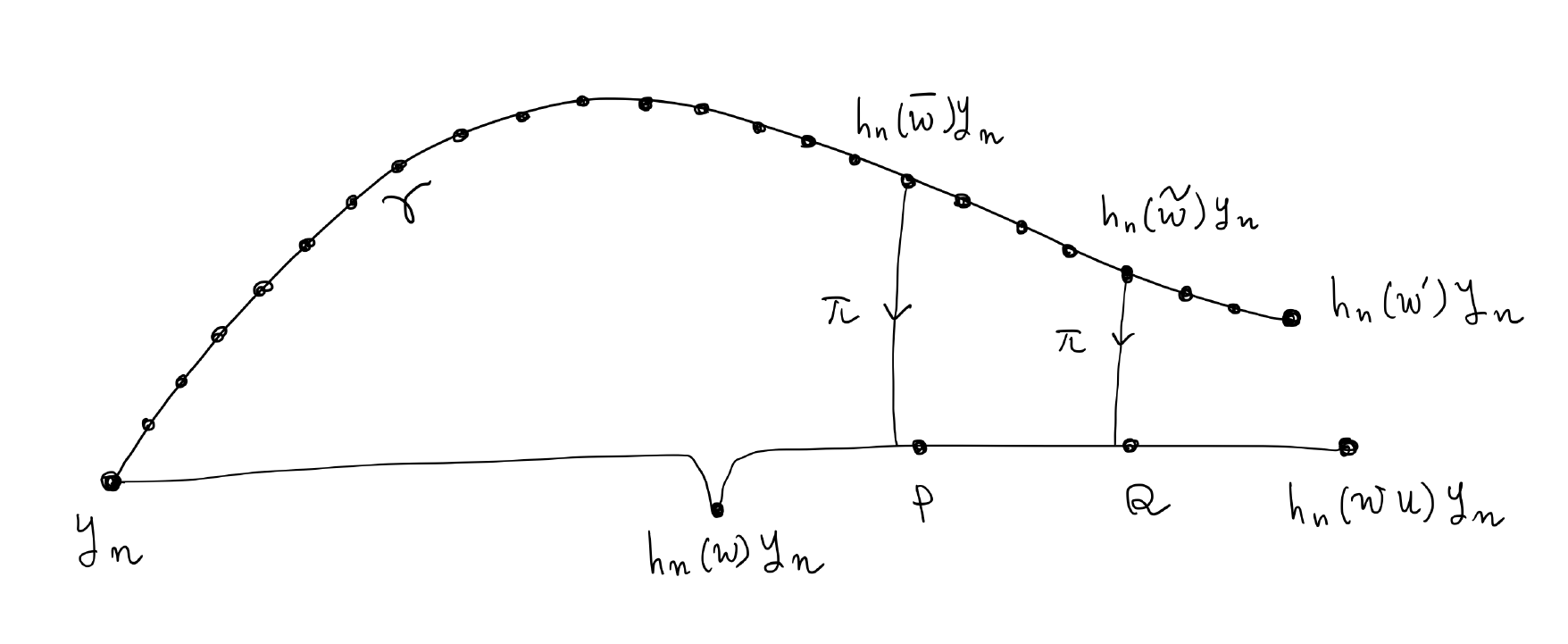}%
   \end{center}
\caption{If $w$ is forbidden, then we have two new candidates $\bar w, \tilde w$ for non-forbidden elements.}
 \label{figure.forbidden}
\end{figure}

By Lemma 4.2 (ii), 
the length of a tail is at least $100\Delta_n$.
Let $P, Q$ be the two points on the tail that trisect the tail into
three pieces of equal length, where $P$ is closer
to $h_n(w)(y_n)$ than $Q$ is. Each of the three
pieces has lengh at least $33\Delta_n$.
Let $h_n(s_1 \cdots s_p)(y_n)$ be a point on $\gamma$ whose projection is 
closest to $P$, and 
$h_n(s_1 \cdots s_q)(y_n)$  a point whose projection is closest
to $Q$. Then 
$$|\pi(h_n(s_1 \cdots s_p)(y_n))-P| \le L(S_n^{2MD})\le \frac{\Delta_n}{4} , 
$$
$$
|\pi(h_n(s_1 \cdots s_q)(y_n))-Q| \le L(S_n^{2MD}) \le \frac{\Delta_n}{4}.$$
 We denote $s_1 \cdots s_p,  s_1 \cdots s_q\in B_m$
 as $w_0, w_1$ and call them
  the {\it candidates} for 
non-forbidden elements obtained from $w$. They depend on $u, w'$ too. See Figure \ref{figure.forbidden}.

Consider the union of the three segements as follows:
\begin{align*}
\Pi(w) = [h_n(w)(y_n), h_n(wu)(y_n)]
&\cup  [h_n(w_0)(y_n), \pi(h_n(w_0)(y_n))]
\\
&\cup [h_n(w_1)(y_n), \pi(h_n(w_1)(y_n))].
\end{align*}

This depends not only $w$, but also $u, w'$.
$\Pi(w)$ is a tree
embedded in $X$. 
We call $\pi(h_n(w_0)(y_n))$ and $\pi(h_n(w_1)(y_n))$ the {\it branch points} of $\Pi(w)$, 
and $P,Q$ the {\it trisecting points} of (the tail of) $\Pi(w)$.
This union  is a disjoint union except for the two 
branch points. The distance betwee the branch points is 
at least $32 \Delta_n$. This in particular implies that 
$w_0 \not = w_1$. 

%Also
%$[y_n,h_n(w)(y_n)] \cap  N_{10\delta}(\Pi(w))$ is contained in the $(2\Delta_n+100\delta)$-neighborhood 
%of $h_n(w)(y_n)$, 
%since $u$ is admissible for $w$ and  the points $P,Q$
%are far from $h_n(w)(y_n)$. In this sense, the union $[y_n,h_n(w)(y_n)] \cup  \Pi(w)$
%looks like a tree, and we call it a tree like graph. 
%In particular $\bar w \not= \tilde w$.
%kf{put figure}

We set
$C_1(w)=\{w,w_0, w_1 \} \subset B_m$, and 
$T_1(w)=[y_n,h_n(w)(y_n)] \cup  \Pi(w)$.
If both $w_0, w_1$ are non-forbidden, this is the 
end of the construction, and put
$C(w)=C_1(w), T(w)=T_1(w)$. Otherwise we go to the second step. 

$T_1(w)$ is a tree like graph in the sense that it is 
the union of two trees, $[y_n,h_n(w)(y_n)]$ and $\Pi(w)$
attached at the point $h_n(w)(y_n)$, and the intersection is 
contained in the $(2\Delta_n+100\delta)$-neighborhood of this point
since $u$ is admissible for $w$. 

{\it Step 2}.
By assumption, at least one of 
$w_0$ and $w_1$ is forbidden.
For example, assume that $w_1$ is forbidden.
Then $w_1$ has  a separator $u_1$ admissible for $w_1$, and an element $w_1' \in B_m$
that satisfy (\ref{forbidden}).
Let $\pi_1$ denote the nearest point projection to the tail
$[h_n(w_1)(y_n), h_n(w_1u_1)(y_n)]$.

Then as in the first step, there is a path (a sequence of points obtained from a shortest expression of $w_1'$ in $\eta(S)$) between $y_n$ and $h_n(w_1')(y_n)$
in $X$,  from which we obtain two elements,
$w_{10}, w_{11} \in B_m$, using  the projection $\pi_1$ from the path to the tail
$[h_n(w_1)(y_n), h_n(w_1u_1)(y_n)]$.
Also, we obtain a tree embedded in $X$:
\begin{align*}
\Pi(w_1)=[h_n(w_1)(y_n), h_n(w_1u_1)(y_n)]
&\cup [h_n(w_{10})(y_n), \pi(h_n(w_{10})(y_n))]
\\
&\cup [h_n(w_{11})(y_n), \pi(h_n(w_{11})(y_n))].
\end{align*}

We point out that $([y_n,h_n(w)(y_n)]\cup \Pi(w) )  \cap N_{10\delta}(\Pi(w_1))$ is 
contained in the $(2\Delta_n+100\delta)$-neighborhood of $h_n(w_1)(y_n)$.
%where $\Pi(w)$ and $\Pi(w_1)$ are connected. 
The reason is because $u$ is admissible to $w$,
and $u_1$ is admissible to $w_1$, while
the tails of $(w,u)$ and $(w_1,u_1)$ are both at
least $100 \Delta_n$ long.

Also, the trisecting points $P_1,Q_1$ and the branch points on the tail $ [h_n(w_1)(y_n),h_n(w_1u_1)(y_n)]$
are out of the $10\Delta_n$-neighborhood of the tail $[h_n(w)(y_n),h_n(wu)(y_n)]$.

It follows 
%that  $[y_n, h_n(w)(y_n)] \cup \Pi(w) \cup \Pi(w_1)$ looks like a tree, and  in particular 
that $w,w_1,w_{10},  w_{11}$ are distinct elements in $B_m$. 
If $w_0$ is not forbidden, then we put
$C_2(w)=C_1(w) \cup \{w_{10}, w_{11} \}$  and 
$$T_2(w)=T_1(w) \cup \Pi(w_1).$$

\begin{figure}[htbp]
\hspace*{-3.3cm}     
\begin{center}
                                                      
   \includegraphics[scale=0.6]{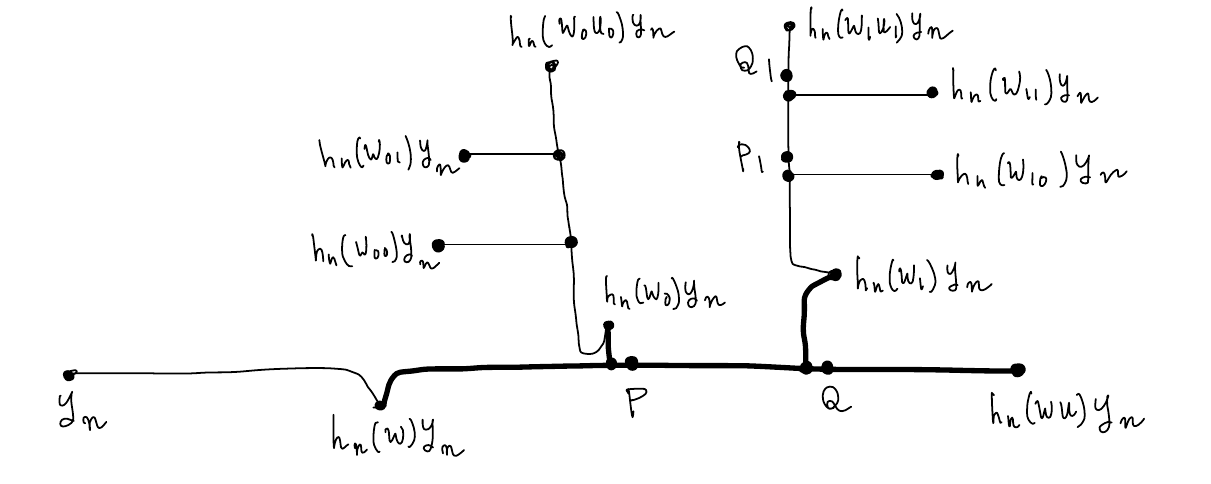}%
   \end{center}
\caption{This is $T_2(w)$ in the case that $w,w_0,w_1$ are forbidden.
$\Pi(w)$ is 
in bold lines. We have two more trees $\Pi(w_0), \Pi(w_1)$.
Those three trees stay away from each other except
for the balls of radius $2\Delta_n$ at the points where they
are connected. 
}
 \label{figure.tree}
\end{figure}

If $w_0$ is forbidden, then we do the same construction as we 
did to $w_1$ and obtain two candidates:
$w_{00}, w_{01} \in B_m$, and a tree $\Pi(w_0)$
which intersects the previously obtained treelike graph
$[y_n, h_n(w)(y_n)] \cup \Pi(w) \cup \Pi(w_1)$
only in  the $(2\Delta_n+100\delta)$-neighborhood of $h_n(w_0)(y_n)$.
In particular,  
$w,w_0,w_{00}, w_{01}, w_1, w_{10}, w_{11}$ (if they exist) are
all distinct. 
See Figure \ref{figure.tree}

We set 
$C_2(w)=C_1(w) \cup \{w_{00}, w_{01}, w_{10}, w_{11} \}$  and 
$$T_2(w)=T_1(w) \cup \Pi(w_0) \cup \Pi(w_1).$$

If all of  the elements $w_{00}, w_{01}$ and
also $w_{10}, w_{11} \in B_m$ (if they exist) 
are non-forbidden, we stop here 
%We obtain 4 (or 3) non-forbidden
%elements from the 3 (or 2) forbidden elements $w, \bar w, \tilde w$.
and set $C(w)=C_2(w), T(w)=T_2(w)$, otherwise we go to the third
step. 
%Here, $\bar w_2, \tilde w_2, \Pi(\tilde w)$ do not exist if $\tilde w=w_2$ is non-forbidden. 

In the third step, we do the same construction to 
the forbidden elements among 
$w_{00},  w_{01}, w_{10}, w_{11} \in B_m$
For  example, if $w_{00}$ is forbidden, then we obtain two
elements $w_{000}, w_{001}$ and a tree $\Pi(w_{00})$.
This tree 
intersects the tree like graph $T_2(w)$ in the $(2\Delta_n+100\delta)$-neighborhood 
of $h_n(w_{00})(y_n)$.
We define $C_3(w)$ from $C_2(w)$ by adding a pair of elements
that are obtained from each forbidden element $w_{**}$ in this step.
Also, we define $T_3(w)$ be adding the trees $\Pi(w_{**})$ 
for forbidden elements $w_{**}$. 

We formally describe the $(N+1)$-th step.
In the $N$-step, we have $C_N(w), T_N(w)$.
If all of the elements in $C_N(w) \backslash C_{N-1}$
are non-forbidden, then we stop and put $C(w)=C_N(w),
T(w)=T_N(w)$. Otherwise, we go to the $(N+1)$-th step. 

{\it $(N+1)$-th step}.
By assumption, at least one element in $C_N(w) \backslash C_{N-1}$,
$w_{i_1 \cdots i_N}$,  is forbidden. Then we do the
same construction as in the step 1, and obtain 
the tree $\Pi(w_{i_1 \cdots i_N})$; two elements, $w_{i_1 \cdots i_N 0}$ and 
$w_{i_1 \cdots i_N 1}$. We add those two elements to $C_N(w)$, 
and also add  $\Pi(w_{i_1 \cdots i_N})$ to $T_N(w)$.
 We do this all forbidden elements in 
 $C_N(w) \backslash C_{N-1}(w)$, and get 
 $C_{N+1}(w)$ and $T_{N+1}(w)$. 
 Note that the trees $\Pi(w_I)$ we added in this step
 are disjoint from each other, moreover, the distance 
 between two $\Pi(w_I), \Pi(w_{I'}$ with $I \not=I'$
 is at least $30\Delta_n$, where $I,I'$ are multi-subscripts of length $N$. 

Since $C_n(w) \subset B_m$ for all $n$, and the set $C_n(w)$
gets bigger at least by two as $n$ increases by one, 
so the process  must end in finite steps since $B_m$ is finite. 
At  the end of the process, we have a set $C(w)$ in $B_m$, 
and a tree like graph $T(w)$. 
We remark that $C(w), T(w)$ depend on the choice of 
separators in the construction, but it is not important. 

No that by construction, in $C(w)$, we have
$$|\{\text{forbidden elements}\} | +1 = |\{\text{non-forbidden elements}\}|.$$
%We add all of $\Pi(w_I)$ for forbidden 
%$w_I \in C(w)$ to $[y_n,h_n(w)(y_n)]$, which 
%is a tree like graph, denoted  by
%$T(w)$.
Also, the orbit of $y_n$ by the elements in 
$C(w)$ appears as vertices of $T(w)$.
%The construction of $C(w),T(w)$ is done. 

{\em Part II}. 
Let $w_1, w_2 \in B_m$ be two forbidden elements,
and we analyze how the sets $C(w_1)$ and $C(w_2)$ are related. 
In the construction of $C(w)$, for each forbidden element, $v$,
we chose a separator $u$ and an element $v' \in B_m$.
The pair $(v,u)$ defines the tail $[h_n(v)(y_n), h_n(vu)(y_n)]$,
which is an arc of $T(w)$.

{\it Claim}. 
Let $w_1, w_2\in B_m$ be  forbidden elements and 
$u_2$ the separator we chose for $w_2$ to construct $T(w_2), C(w_2)$.
Assume  that the tail
 $\tau=(w_2,u_2)$ is not parallel to any of $(v,u)$
that appears to construct $T(w_1),C(w_1)$.
Then, 
$C(w_1) \cap C(w_2)$ is empty or $\{w_2\}$.

Let $P,Q \in \tau$ be the trisecting points of $\Pi(w_2)$.
Then neither of them is in the $10\Delta_n$-neighborhood
of any trisecting point that appears in $T(w_1)$.
Indeed, suppose it was, and let $R$ be a trisecting point
of the tail, $\sigma$, in $T(w_1)$ such that $|P-R|$ or $|Q-R|$
is $\le 10\Delta_n$.
Then, since $T(w_1), T(w_2)$ have $y_n$ as the 
common ``root'' vertex that they start from, 
the intersection of $\sigma$ and 
the $20\delta$-neighborhood of $\tau$ 
is at least $\frac{1}{5}$ of $\sigma$, which 
contradicts the small cancellation property.

This implies that not only the branch points on $\tau$, but also 
all branch points of $T(w_2)$ are 
outside of the $9\Delta_n$-neighborhood of $T(w_1)$.
This is because both $T(w_1), T(w_2)$ are tree like graphs.
It then follows that $C(w_1) \cap C(w_2)$ is empty or 
just $w_2$ because all other points in $C(w_2)$
appear on $T(w_2)$ after the first two branch points on $\tau$,
that are close to $P,Q$.
We showed the claim.

{\em Part III}.
It follows from the claim that there is a finite collection of forbidden elements $w$
in $B_m$ such that the $C(w)$'s are mutually disjoint, 
and that any forbidden element $v \in B_m$ is either contained
in the union of those $C(w)$'s, or $v$ has an admissible separator $u$
such that the pair (tail) $(v,u)$ is  parallel to one of the tails
that appears in one of the $C(w)$'s. 

To see that, order the forbidden elements in $B_m$ as $w_1, w_2, \cdots$.
Choose admissible separators $u_1, u_2, \cdots$.
First, construct $C(w_1)$ using the separators we chose. 
Next, if $w_2$ is already  contained in $C(w_1)$, then disregard it.
Also, if $(w_2,u_2)$ is parallel to one of the tails in $T(w_1)$, then
disregard $w_2$ also. Otherwise, 
construct $C(w_2)$ and keep it. 
 If $w_3$ is contained in $C(w_1)$ or $C(w_2)$; or 
 $(w_3,u_3)$ is parallel to one of the tails of $T(w_1)$ or $T(w_2)$, 
 then disregard $w_3$. Otherwise, construct $C(w_3)$ and keep it,
 and so on. 
 By the claim we have shown, those $C(w)$'s are mutually disjoint
 such that 
for any forbidden element $v \in B_m$,
either $v$ is contained in one of the $C(w)$'s we have, or $(v,u)$
is parallel to one of the tails of the $T(w)$'s. 

To finish the proof,
by Lemma \ref{lemma.tail}, the union of those $C(w)$'s contains at least $\frac{1}{\D}$ of the forbidden 
elements in $B_m$. 
But since in each $C(w)$, there are more non-forbidden 
elements than forbidden elements (by one), 
 we conclude that in $B_m$, 
$$|\{\text{forbidden elements}\} | \le \D |\{\text{non-forbidden elements}\}|.$$
Lemma \ref{4.6} is proved. 
\qed

\subsection{Feasible elements}
By the lemma \ref{4.6}, at least $\frac{1}{\D+1}$ of $B_m(L,\eta(S))$ consists of non-forbidden elements. 

%Again,  for two distinct non-forbidden elements $w,w'$, 
%it is possible that $h_n(w), h_n(w')$ are 
%in the same coset w.r.t. $F(h_n(u_i))$, where $u_i$ is a separator for $h_n$.

Now, we choose a maximal subset in the set of non-forbidden
elements in $B_m(L,\eta(S))$ such that for any two distinct elements 
$w,w'$ in the set, $h_n(w),h_n(w')$ are not 
in the same coset w.r.t. $F(h_n(u_i))$ for any separator $u_i$.
We call an element in this set an {\it adequate} element.
Then, in $B_m(L,\eta(S))$, 

$$\frac{|\{\text{non-forbidden elements}\} |}{\D^4} \le |\{\text{adequate elements}\}|.$$

This is because there are only four separators, $u$,
and $|F(h_n(u))| \le \D$ for each $u$.
Combining this with Lemma \ref{4.6} we get:
\begin{equation}\label{adequate}
\frac{|B_m(L,\eta(S))|}{\D^4(\D+1)} \le |\{\text{adequate elements in } B_m(L,\eta(S))\}|.
\end{equation}

\begin{definition}[feasible elements]
Fix $m$, then $n$. 
An element of the form $w_1u_1 \cdots w_q u_q$
with $w_i\in B_m(L,\eta(S))$
are called {\it feasible} of {\it type} $q$ if all $w_i$ are adequate
and each $u_i$ is admissible for $w_i,w_{i+1}$.
We define $q=0$ for the empty element. 
\end{definition}

For this feasible elements, we consider the following broken geodesic:
%\kf{need $h_n$}
\begin{align*}
\alpha&=[y_n,h_n(w_1)(y_n)] \cup [h_n(w_1)(y_n),h_n(w_1u_1)(y_n)]  \cup \cdots \\
&\cup
[h_n(w_1u_1 \cdots w_q )(y_n), h_n(w_1u_1 \cdots w_q u_q)(y_n)].
\end{align*}
The Hausdorff-distance between $\alpha$ and the geodesic
$[y_n, h_n(w_1u_1 \cdots w_q u_q)(y_n)]$ is at most $2\Delta_n+100\delta$.

\begin{lemma}[cf. \cite{FS},Lemma 2.6]\label{feasible}
For every $q$, the map
$\Phi_n$ is injective on the set of $q$-tuples of adequate elements
in $B_m(L,\eta(S))$.
\end{lemma}
By definition, given a $q$-tuples of adequate elements,
we choose admissible separators (they are not unique) and form a feasible element of type $q$, 
then map it by $\Phi_n$.

\proof
%We first argue the case such that the two elements
%have the same type $q$.
We argue by induction on the type $q$.
If $q=0$ then nothing to prove since the element is empty.

Assume the conclusion holds for $q \ge 0$.
Suppose $(q+1)$-tuples of adequate elements
$w_1, w_2, \cdots, w_{q+1}$ and $w_1', w_2', \cdots, w_{q+1}'$ are
given.
Let $W=w_1u_1w_2u_2 \cdots w_{q+1}u_{q+1}$ and 
$W'=w_1'u_1'w_2'u_2' \cdots w_{q+1}'u_{q+1}'$ be two feasible elements
of type $q+1$.
We assume that $h_n(W)=h_n(W')$ and want to show
$w_i=w_i'$ for all $i$.
%Indeed in the following we will show 
%$W=W'$ in the sense that $w_i=w_i'$ and $u_i=u_i'$ for all $i$.

The elements $W,W'$ define two broken geodesics $\alpha, \alpha'$
between $y_n$ and $h_n(W)(y_n)=h_n(W')(y_n)$.
Look at the initial parts of $\alpha, \alpha'$:
$[y_n,h_n(w_1)(y_n)] \cup [h_n(w_1)(y_n),h_n(w_1u_1)(y_n)]$
and $[y_n,h_n(w_1')(y_n)] \cup [h_n(w_1')(y_n),h_n(w_1'u_1')(y_n)]$.
We first show that $u_1=u_1'$. Suppose not. 
Then it would imply that either $w_1$ or $w_1'$ is forbidden.
We explain the reason. 
First, by  the small cancellation property of the separators
(Lemma \ref{4.2} (iii)), either $[y_n,h_n(w_1)(y_n)] \cup [h_n(w_1)(y_n),h_n(w_1u_1)(y_n)]$ is contained in the $10 \Delta_n$-neighborhood
of $[y_n,h_n(w_1')(y_n)]$, or $[y_n,h_n(w_1')(y_n)] \cup [h_n(w_1')(y_n),h_n(w_1'u_1')(y_n)]$ is contained in the $10 \Delta_n$-neighborhood
of $[y_n,h_n(w_1)(y_n)]$,
since $u_1$ is admissible for $w_1$ and $u_1'$ is admissible for $w_1'$.
Suppose we are in the first case. Then $w_1$ would be forbidden since 
there is a subword of $w_1'$, denoted by $\overline{w_1'}$, 
such that 
 $$|h_n(\overline{w_1'})(y_n)-h_n(w_1 u_1)(y_n)| \le \frac{1}{5}|h_n(w_1)(y_n)-h_n(w_1 u_1)(y_n)|.$$
 Such a subword exists because $L(S_n) \le \frac{\Delta_n}{4}$
 and $|h_n(w_1)(y_n)-h_n(w_1 u_1)(y_n)| \ge 100 \Delta_n$.
 We got a contradiction since $w_1$ is not forbidden.
 Similarly, if we are in the second case, then $w_1'$
 would be forbidden, which is a contradiction. 
 We showed that $u_1=u_1'$.

%$$|h_n(w_1)(y_n)-h_n(w_1' u_1')(y_n)| \le \frac{1}{5}
%|h_n(w_1')(y_n)-h_n(w_1' u_1')(y_n)|;$$
%or $$|h_n(w_1')(y_n)-h_n(w_1 u_1)(y_n)| \le \frac{1}{5}
%|h_n(w_1)(y_n)-h_n(w_1 u_1)(y_n)|.$$
%But then, in the first case, $w_1'$ is forbidden, a contradiction.
%In the second case, $w_1$ is forbidden, a contradiction.
%We showed $u_1=u_1'$

Next we show $w_1=w_1'$. We first show that $h_n(w_1 w_1'^{-1})\in F(h_n(u_1))$.
Suppose not. Then, as in the previous paragraph, the small cancellation property 
of the separators implies that either $w_1$ or $w_1'$
is forbidden, a contradiction. 
We are left with the case that $h_n(w_1 w_1'^{-1}) \in F(h_n(u_1))$.
But, this means  that $h_n(w_1)$ and $h_n(w_1')$ are in the same (right) coset
w.r.t. $F(h_n(u_1))$. Since both $w_1$ and $w_1'$ are adequate, it implies 
$w_1=w_1'$.

Since  $u_1=u_1', w_1=w_1'$, 
it follows that $W_1=w_2u_2 \cdots w_{q+1}u_{q+1}$
and $W_1'=w_2'u_2' \cdots w_{q+1}'u_{q+1}'$
are feasible elements of type $q$
with $h_n(W_1)=h_n(W_1')$.
By the induction hypothesis, we have 
$w_i=w_i'$ for all $i \ge 2$, and we are done. 
%we apply the previous argument to the part $w_1u_1$ and $w_1'u_1'$, 
%we conclude that $w_1=w_1', u_1=u_1'$. Then
%$h_n$ maps $w_2u_2 \cdots w_{q+1}u_{q+1}$ and $w_2'u_2' \cdots w_{q+1}'u_{q+1}'$
%to the same element and they 
%are feasible elements of type $q$. So, they must be equal by the induction hypothesis.
%We conclude that the given two elements are identical.
\qed

%%%%%%%%%%%

%%%%%%%%%%%%%%%%%%%%

\subsection{Proof of Proposition \ref{1.2}}\label{section.4.5}
We prove Proposition \ref{1.2}.
\proof
%We are assuming that the action is $D$-uniformly WPD.
Recall that $D(\epsilon)$ is the WPD function and 
we set $\D=D(100\delta)$. 
Recall that for every $m$, we have
$$\frac{1}{m} \log (|B_m(L,\eta(S))|  \ge \log e(L,\eta(S)).$$
Given $\epsilon>0$, choose and fix a large enough $m$ such that 

 $$\frac{1}{m+b}\left(\log |B_m(L,\eta(S))|-\log (D^4(D+1)\right) \ge \frac{1}{m}\log |B_m(L,\eta(S))| -\epsilon.$$
 The constant $b$ is from Lemma \ref{4.2}.
 Such $m$ exists  since 
 $\lim_{m\to \infty}\frac{\log|B_m(L,\eta(S))|}{m}=\log e(L,\eta(S))$, and $b$ does not depend on $n,m$.
 (We will choose $n$ right after this.)

%\kf{changed 2m+b to 2m}
Choose $n$ large enough such that $h_n$ is injective 
on $B_{2m}(L,\eta(S))$, which 
defines the forbidden elements in $B_m(L,\eta(S))$. Then for all $q$:
$$ |B_{q(m+b)}(G,f_n(S))|\ge \left(\frac{|B_m(L,\eta(S))|}{\D^4(\D+1)}\right)^q $$
because the number of adequate elements in $B_m(L,\eta(S))$
is at least $\frac{|B_m(L,\eta(S))|}{\D^4(\D+1)}$ by the estimate (\ref{adequate}), and $\Phi_n$ is injective
on the set of feasible elements of type $q$ by Lemma \ref{feasible}.

Then 
by the above three inequalities, 
\begin{align*}
 \log(e(G,f_n(S)) &= \lim_{q\to \infty} \frac{1}{q(m+b)}\log|B_{q(m+b)}(G,f_n(S))|
 \\
 &\ge \frac{1}{m+b} \log \left(\frac{|B_m(L,\eta(S))|}{\D^4(\D+1)}\right) \ge  \frac{1}{m} \log (|B_m(L,\eta(S))| -\epsilon
 \\
&\ge \log e(L,\eta(S))- \epsilon.
\end{align*}

Since we have this for all large enough $n$, 
and $\epsilon>0$ is arbitrary,
we have 
$$\lim_{n \to \infty} \log(e(G,f_n(S))\ge \log e(L,\eta(S)).$$
As we said the other direction is trivial, hence
the equality holds. 
\qed

\subsection{A family version}\label{section.family}
We state a family version of Proposition \ref{1.2}.
We do not use this proposition in this paper, but 
it would be useful in the future. 

We explain the setting. 
Let $\delta >0$ be a constant, $M>0$ an integer
and $D(\epsilon)$ a function for WPD.
Suppose we have a family of finitely generated groups
$\{G_\alpha\}$ such that each $G_\alpha$ is not virtually cyclic and acts on a $\delta$-hyperbolic space $X_\alpha$.
Assume that for any 
finite generating set $S$ of any $G_\alpha$, the set $S^M$ contains a hyperbolic 
element on $X_\alpha$ that is $D$-WPD.

Now suppose that $S_n$ is a finite generating set of $G_n$
that is in the family for $n>0$.  Assume that the (infinite) sequence
$e(G_n,S_n)$ is bounded from above. 
Then by Proposition \ref{bound.generators}, there
is a constant $A>0$ that depends only on $M, \delta$ and the function $D$ such that $e(G_n,S_n) \ge A|S_n|^A$, therefore, 
the sequence $|S_n|$ is bounded from above.

As before, passing to a subsequence, we may assume that 
there exists $\ell>0$ with $|S_n|=\ell$ for all $n$. 
This gives a sequence of surjections $f_n:(F,S) \to (G_n, f_n(S))$
with $f_n(S)=S_n$, where $(F,S)$ is a pair of free group $F$ and a free generating set $S$ with $|S|=\ell$. 

Passing to a further subsequence if necessary, 
we may assume that the sequence $f_n$ 
converges to a limit epimorphism $\eta:(F,S) \to (L,\eta)$.

Then we have:
\begin{prop}\label{semi.cont.family}
Assume that there exists an epimorphism $h_n: L \to G_n$ for all $n$ such that  all the homomorphisms 
$\{f_n\}$ factor through the limit epimorphism: $\eta: F \to L$, ie, 
 $f_n=h_n \circ \eta$. 
 Then,
$$\lim_{n \to \infty} e(G_n,f_n(S)) = e(L,\eta(S)).$$

\end{prop}

Before we explain the proof, we point out
that in Proposition \ref{1.2}, the assumption on the existence
of $h_n$ is a consequence of that $G$ is
equationally Noetherian. In the current
setting, we state it as an assumption and 
avoid an assumption related to equational
Noetherianity.

That said, the proof is identical to the argument for Proposition \ref{1.2}. We recall that by the existence of 
$h_n$, we immediately have $\lim_n e(G_n,f_n(S)) \le e(L,\eta(S))$.
The main issue was to show the other inequality
 $\lim_n e(G_n,f_n(S)) \ge e(L,\eta(S))$ 
by constructing separators for  each action of $L$ on $X$, via
$h_n:L \to G_n$.
That argument
applies without change to the current setting
since the constants $\delta, M$ and the function $D(\epsilon)$
are common for all $(G_n,S_n)$ and $X_n$.

\section{Examples}

An obvious example for Theorem \ref{main} is a non-elementary hyperbolic group, $G$. Let $X$ be a Cayley graph of $G$, then 
it is $\delta$-hyperbolic, and the action 
by $G$ is (uniformly) proper, so that 
acylindrical, and non-elementary. The existence of the constant $M$
is known (\cite{Koubi}).
As we said $G$ is equationally Noetherian (\cite{Sela1}, \cite{RW}).
Therefore, $\xi(S)$ is well-ordered. This is proved in \cite{FS}
and we adapted their argument to prove Theorem \ref{main.proof} in this paper. 
In this section we discuss some other examples.

\subsection{Relatively hyperbolic groups}\label{section.5.2}
We treat  a relatively hyperbolic group. For example see \cite{Bo1}
for the definition and basic properties. 
Suppose $G$ is hyperbolic relative to a collection of finitely generated
subgroups $\{P_1, \cdots, P_n\}$. Suppose $G$ is not
virtually cyclic and not equal to any $P_i$. Then it acts properly discontinuously
on a proper $\delta$-hyperbolic space $X$ such that
\cite[Proposition 6.13]{Bo1}:
\begin{enumerate}
\item
There is a $G$-invariant 
collection of points, $\Pi \subset \partial X$, with 
$\Pi/G$ finite. 
For each $i$, there is a point $p_i \in \Pi$ such that the stabilizer
of $p_i$ is $P_i$. The union of the $G$-orbits of the $p_i$'s is $\Pi$. 
\item
For every $r>0$, there is a $G$-invariant collection of horoballs $B(p)$
at each $p\in \Pi$ such that they are $r$-separated, ie,
$d(B(p),B(q)) \ge r$ for every distinct $p,q \in \Pi$.
\item
The action of $G$ on $X \backslash \cup_{p\in \Pi} {\rm int}B(p)$
is co-compact.

\end{enumerate}
A subgroup $H<G$ is called {\it parabolic} if it is infinite, 
fixes a point in $\partial X$, and contains no hyperbolic elements.
The fixed point is unique and called a {\it parabolic point}.
In fact, $\Pi$ is the set of parabolic points
(\cite[Proposition 6.1 and 6.13]{Bo1}).

For this action, we have:

\begin{lemma}[Proposition 5.1\cite{X}]\label{Xie}
Let $G$ and $X$ be as above. Then there exists $M$ such that 
for any finite generating set $S$ of $G$, the set $S^M$
contains a hyperbolic element on $X$.
\end{lemma}

Also we have the following:
\begin{lemma}\label{rel.hyp.wpd}
Let $G$ and $X$ be as above. Then the action is
uniformly WPD.
\end{lemma}

\proof
Given $\epsilon>0$, take a $G$-invariant collection of horoballs
that are $(10\epsilon+100\delta)$-separated in $X$.
Let $\mathcal B$ denote the union of the interior of the horoballs in the collection.
Then since the action of $G$ on $X\backslash \mathcal B$
is co-compact, there exists a constant $D(\epsilon)$ such that 
for any $y \in X\backslash \mathcal B$, the cardinality of the following set is at most $D(\epsilon)$:
$$\{h\in G| |y-h(y)| \le \epsilon + 10 \delta \}.$$
We argue that the action is uniformly WPD w.r.t. $D=D(\epsilon)$.
Let $g\in G$ be hyperbolic with an $10\delta$-axis $\gamma$.
Let $x, y \in \gamma$ with $|x-y| \ge \lambda(g)$.
It suffices to show that the following set
contains at most $D(\epsilon)$ elements:
$$\{h\in G||x-h(x)|\le \epsilon, |y-h(y)|\le \epsilon\}.$$

We divide the case into four:
\\
(1)
$x \not\in \mathcal B$.
Then there are at most $D(\epsilon)$ elements $h \in G$
s.t. $|x-h(x)|\le \epsilon$, and we are done.
\\
(2)
$y \not\in \mathcal B$.
This is same as (1).
\\
(3) $x,y \in \mathcal B$ such that 
$x \in B(p)$ and $y \in B(q)$ with $p\not= q$.
Then, there must be $z \in [x,y]$ with $z
\not\in N_{5\epsilon + 30 \delta}(\mathcal B)$
since the horoballs in $\mathcal B$
are $(10\epsilon + 100 \delta)$-separated. 
Then $|z-h(z)| \le \epsilon +10\delta$.
But there are at most $D(\epsilon)$ such elements.
\\
(4) $x,y \in B(p)$ for some $p$. 
Then $g(x) \in B(q)$ for some $q \not= p$,
since $g$ is hyperbolic and does not preserve
any horoball. Since horoballs are $(10\epsilon + 100 \delta)$-separated, we have $|x-g(x)|\ge 10 \epsilon + 100 \delta$.
So, $\lambda(g) \ge 10 \epsilon + 50 \delta$.
Now, there are two possibilities:
one is that $g(x) \in \gamma$ is between $x$ and $y$. Then, as in (3), there must be $z \in [x,y]$ with $z
\not\in N_{5\epsilon + 30 \delta}(\mathcal B)$,
and we are done. 
The other possibility is that $y$ is between 
$x$ and $g(x)$ on $\gamma$.
But in this case, since $|x-y| \ge \lambda(g)$,
we have $|y-g(x)| \le 50 \delta$. 
Then the distance between $B(p)$ and $B(q)$
is at most $50\delta$ since $y \in B(p)$
and $g(x) \in B(q)$. But it contradicts
the separation of horoballs, so this case does
not happen.
\qed

We quote a theorem, \cite[Theorem D]{GrH}.
\begin{thm}[Equationally Noetherian, \cite{GrH}]\label{rel.hyp.EN}
If $G$ is hyperbolic relative to equationally Noetherian subgroups, 
then $G$ is equationally Noetherian.
\end{thm}

We are ready to state a theorem. 

\begin{thm}[Well-orderedness for relatively hyperbolic groups]\label{rel.hyp.growth}
Let $G$ be a group that is hyperbolic relative to a collection of 
subgroups $\{P_1, \cdots, 
P_n\}$.
Suppose $G$ is not virtually cyclic, and not equal to $P_i$
for any $i$.
Suppose each $P_i$ is finitely generated and equationally Noetherian.
Then $\xi(G)$ is well-ordered. 

\end{thm}

\proof
$G$ is equationally Noetherian by Theorem \ref{rel.hyp.EN}.
Take a hyperbolic space $X$ with a $G$ action 
as above. The action is non-elementary
since $G$ contains a hyperbolic isometry and 
 is not virtually $\Bbb Z$.
Then Theorem \ref{main.proof} applies
%(as we explained $X$ does not have to be a graph)
by Lemma \ref{Xie} and Lemma \ref{rel.hyp.wpd}.
\qed

\subsection{Rank-1 lattices}
There are many examples of relatively hyperbolic groups, but 
we mention one standard family that Theorem \ref{rel.hyp.growth} applies to.

Let $G$ be a lattice in a simple Lie group of rank-1. It is 
always finitely generated and has exponential growth, and in 
fact uniform exponential growth, \cite{EMO}. If it is 
a uniform lattice, then it is hyperbolic, so that
$\xi(G)$ is well-ordered.
We prove the following as an immediate application of
Theorem \ref{rel.hyp.growth}.

\begin{thm}[Rank-1 lattices]\label{lattice}
Let $G$ be one of the following groups:
\begin{enumerate}
\item
A lattice in a simple Lie group of rank-1. 
\item
The fundamental group of a complete Riemannian manifold $M$
of finite volume such that there exist $a,b>0$ with 
$-b^2 \le K \le -a^2 <0$, where $K$ denotes the
sectional curvature. 
\end{enumerate}

Then $\xi(G)$ is well-ordered.

\end{thm}

\proof 
We only need to argue for non-uniform lattices since 
otherwise, $G$ is a non-elementary hyperbolic group
and the conclusion holds. 
Suppose that $G$ is a non-uniform lattice. Then, 
it is known that $G$ is relatively hyperbolic w.r.t.
the the parabolic subgroups, called {\it peripheral subgroups} $\{H_i\}$ that are associated to the cusps,
\cite{Fa}, \cite{Bo1}.
Moreover $G$ is not virtually cyclic, and not equal to any $H_i$. 
Also, those $H_i$ are finitely generated virtually nilpotent groups (see for example,
\cite{Fa}).
It is known that finitely generated virtually nilpotent groups are equationally
Noetherian, \cite{Br}, so that $H_i$ are equationally Noetherian. 
With those facts, Theorem \ref{rel.hyp.growth}
applies to $G$ and the conclusion holds.
\qed

We remark that it is known that if $G$ is linear over a field, then it is equationally Noetherian,
\cite{BMR}. For lattices in a simple Lie group, one can 
apply this result as well to see that $G$ is equationally Noetherian (consider the adjoint representation on 
its Lie algebra. It's faithful since the Lie group is simple).

%
%
%We prove Theorem \ref{lattice}.
%\proof
%By Lemma \ref{lattice.EN}, 
%$G$ is equationally Noetherian.
%Since Proposition \ref{prop.rank1} applies to $G$, 
%Theorem \ref{main} (or \ref{main.proof}) applies to $G$.
%\qed

\subsection{Mapping class groups}\label{section.mcg}
We discuss mapping class groups as another possible application.
%subsection{Mapping class groups}
Let $MCG(\Sigma)$ be the mapping 
class group of a compact oriented surface $\Sigma=\Sigma_{g,p}$
with genus $g$, punctures $p$, and {\it complexity} $c(\Sigma)=3g+p$.
It is known that it is either virtually abelian 
or has exponential growth, and then uniform exponential growth
(\cite{M2}).
Let  $\CC(\Sigma)$ be the {\it curve graph} of $\Sigma$. The graph
$\CC(\Sigma)$ has vertex set representing the non-trivial  homotopy
classes of simple closed curves on $\Sigma$, and edges
joining vertices representing the homotopy classes of disjoint 
curves. 
The group $MCG(\Sigma)$ naturally acts on it by isometries. 

We assume $c(\Sigma) > 4$. Then $MCG(\Sigma)$ has exponential growth.
We recall some facts:
\begin{enumerate}
\item
The graph $\CC(\Sigma)$ is $\delta$-hyperbolic, and 
any pseudo-Anosov element in $MCG(\Sigma)$ acts hyperbolically on $\CC(\Sigma)$,
\cite{MM1}.
\item
The action of $MCG(\Sigma)$ on $\CC(\Sigma)$ is acylindrical, \cite{Bowditch}. It is non-elementary.
\item
There exists $T(\Sigma)>0$ such that for any 
pseudo-Anosov element $g$, 
$\lambda(g) \ge T(\Sigma)>0$ for the action 
on $\CC(\Sigma)$, \cite{MM1}. (This follows from the acylindricity as we pointed out.)
%\kf{put explanation}

\item
There exists $M(\Sigma)$ such that for any finite $S \subset MCG$
such that $\<S\>$ contains a pseudo-Anosov element, then 
$S^M$ contains a pseudo-Anosov element, \cite{M}.

\end{enumerate}

In summary, the action of $MCG(\Sigma)$ on $\CC(\Sigma)$
is known to satisfy all the assumptions of Theorem \ref{main}
except we do not know if $MCG(\Sigma)$ is equationally 
Noetherian or not.
We expect it to hold (for example, see an announcement 
by Daniel Groves in \cite{GrH}), but it does not exist in the literature
yet. 
Once that is verified, it will imply that 
$\xi(MCG(\Sigma))$ is well-ordered
if $c(\Sigma) > 4$.

We remark that for $\Sigma_{1,1}, \Sigma_{1,0}, \Sigma_{0,4}$, the conclusion
holds by \cite{FS} since $MCG(\Sigma)$ is hyperbolic
(a well-known fact, for example, \cite{MM1}).

%\subsection{Questions}
\begin{question}
Let $G=MCG(\Sigma)$. 
Is $\xi(G)$ well-ordered? If so, the infimum is attained. 
 It would be interesting to know
its value and generating sets that attain the minimum
for each $\Sigma$.

\end{question}

\subsection{Three manifold groups}
We discuss three manifold groups. 
Let $M$ be a closed, orientable $3$-manifold. 
$M$ is called {\it irreducible} if $\pi_1(M)$ does not admit
a non-trivial splitting over the trivial group. 

Let $M$ be a closed, orientable, irreducible 3-manifold which is 
not a torus bundle over a circle. 
Then there is a finite collection of embedded disjoint essential tori $T_i$ in $M$ such that 
each connected component of $M \backslash \cup_i T_i$
is geometric, ie, either Seifert fibred, or admitting hyperbolic or Sol-geometry. 
Such a collection of smallest number of tori is called the {\it JSJ decomposition} of $M$. A JSJ-decomposition in this sense exists by the solution of the geometrization conjecture
of 3-manifolds. 
The collection of tori 
is maybe empty. Otherwise we say $M$ has non-trivial 
JSJ-decomposition. 
A non-trivial JSJ-decomposition gives a graph of groups decomposition 
of $\pi_1(M)$ along subgroups isomorphic to $\Z^2$, and $\pi_1(M)$ has exponential growth. Its Bass-Serre tree is called the {\it JSJ-tree}, $T_M$. 
An action of a group $G$ on a tree $T$ is called 
$k$-{\it acylindrical} if for every non-trivial element $g$,
the subtree of fixed points by $g$ is either empty or of diameter at most $k$. 
It is proved that the action of $\pi_1(M)$ on $T_M$ is 
$4$-acylindrical, \cite[Proposition 4.2]{WZ}.
It implies that the action is uniformly $4$-WPD.
It is known that if an action of $G$ on a tree $T$
is $k$-acylindrical for some $k$, then 
it is acylindrical \cite[Lemma 5.2]{MO},
therefore it is uniformly WPD (Lemma \ref{wpd.acyl}).
Moreover, the acylindricity constants $R(\epsilon), N(\epsilon)$ and the uniformly WPD constant $D(\epsilon)$ depend only on $k$ and $\epsilon$. The moreover part easily follows 
from the proof of Lemma 5.2in \cite{MO}.

\begin{thm}\label{3manifold}
Let $M$ be a closed orientable $3$-manifold, and $G=\pi_1(M)$.
If $M$ is one of the following, then $G$ has exponential growth and $\xi(G)$ 
is well-ordered. 
\begin{enumerate}
\item
$M$ is not irreducible such that  $G$ is not isomorphic to $\Z_2 * \Z_2$.
\item
$M$ is irreducible such that it is  not a torus bundle over a circle,
and that it has a non-trivial JSJ-decomposition.
\item
$M$ admits  hyperbolic geometry. 
\item
$M$ is  Seifert fibered such that the base $2$-orbifold is hyperbolic.
\end{enumerate}
\end{thm}

\proof
First, every three manifold group is equationally Noetherian, 
\cite{GHL}.

(1) In this case, $G$ is a non-trivial free product $A*B$. Since it is not $\Z_2 *\Z_2$,
it has exponential growth. $G$ acts on the Bass-Seree tree $T$ of this free 
product. 
Then for any finite generating set $S$, the set $S^2$ contains a
hyperbolic element on $T$, \cite{Se}. The action of $G$ on $T$
is $0$-acylindrical, so that the action is uniformly $D$-WPD for some $D$.
Since $T$ is hyperbolic, Theorem \ref{main.proof} applies with $M=2$.

(2)
In this case, let $T_M$ be the JSJ-tree
of $M$. Then the action of $G$ on $T_M$ is
$4$-acylindrical, so that it is uniformly $D$-WPD
for some $D$.
Theorem \ref{main.proof} applies with $M=2$.

(3)
If $M$ is hyperbolic, then $G$ is a non-elementary, hyperbolic group. Then 
$\xi(G)$ is well-ordered by \cite{FS}. 

(4)
In this case, we have the following exact sequence:
$$ 1 \to \Z \to G \to H \to 1,$$
where $H$ is the orbifold-fundamental group of the base 2-orbifold,
and $\Z$ is the fundamental group of the fiber circle. We denote
this subgroup $C$.
By assumption, $H$ is a non-elementary, hyperbolic group.

We claim that 
$\xi(G)=\xi(H)$.
To see it, let $S$ be a finite generating set of $G$. Let $\bar S$
be the image of $S$ by the projection $G \to H$.
We have $|\bar S^n| \le |S^n|$. 
But it is known that the subgroup $C$ is not distorted in $G$
in the sense that there is a constant $K>0$ such that 
for any $n>0$, we have $|S^n \cap C| \le Kn$, \cite[Proposition 1.2 (2)]{NS}.
It implies that for all $n$, we have 
$|S^n| \le Kn|\bar S^n|$. It follows that $e(G,S)=e(H,\bar S)$, 
therefore $\xi(G) \subset \xi(H)$.
On the other hand, if $S$ is a finite generating set of $H$, 
then there is a finite generating set $\tilde S$ of $G$
which projects to $S$, by lifting each element of $S$ to $G$, then 
adding a generator of $C$, which gives $\tilde S$.
Then as we saw, $e(G,\tilde S)=e(H,S)$, which 
implies $\xi(H) \subset \xi(G)$. We proved the claim.

But by \cite{FS}, $\xi(H)$ is well-ordered, therefore so is 
$\xi(G)$. 
\qed

A torus bundle over a circle either admits Sol-geometry, 
or it is Seifert fibered (see for example, \cite{WZ}).
Therefore the theorem covers all closed, orientable $3$-manifolds
with fundamental groups of exponential growth except 
that $M$ has Sol-geometry.
We leave it as a question.
\begin{question}
Let $M$ a closed, orientable $3$-manifold
that has the geometry of three dimensional solvable group, Sol.
Then is $\xi(\pi_1(M))$ well-ordered ?
Also, in view of Proposition \ref{bound.generators},
is it true that if $|S|\to \infty$ then $e(\pi_1(M),S) \to \infty$?
\end{question}

\section{The set of growth of subgroups}

We discuss the set of growth of subgroups in a finitely generated group $G$.
Define 
$$\Theta(G)=\{e(H,S)| S \subset G, |S| < \infty, H=\<S\>, e(H,S)>1\}.$$
The set $\Theta(G)$ is countable and contains  $\xi(G)$ as a subset.
If $G$ is a hyperbolic group, it is 
known by \cite[Section 5]{FS} that 
$\Theta(G)$ is well-ordered.

\subsection{Subgroups with hyperbolic elements}

As usual, suppose $G$ acts on a $\delta$-hyperbolic space $X$.
We introduce a subset in $\Theta(G)$ as follows:
$$\Theta_X(G)=\{e(H,S)| S \subset G, |S| < \infty, H=\<S\>, e(H,S)>1\},$$
where in addition we only consider $S$ such that  {\it $\< S\>$  contains a hyperbolic element on $X$}.
This set depends on the action on $X$.

\begin{thm}[cf. Theorem \ref{main.proof}]\label{6.2}
Suppose $G$ acts on a $\delta$-hyperbolic space $X$, and $G$ is not virtually cyclic.
Let $D(\epsilon)$ be a function for WPD.
Assume that there exists a constant $M$ such that if $\<S\>$ contains
a hyperbolic element on $X$ for a
finite subset $S \subset G$, then $S^M$ contains a hyperbolic 
element that is $D$-WPD. 
Assume that $G$ is equationally 
Noetherian.
Then, $\Theta_X(G)$ is a well-ordered set.
\end{thm}

\proof
The proof is nearly identical to Theorem \ref{main.proof}.

Let $\{S_n\}$ be a sequence of finite generating sets of  subgroups
of exponential growth in $G$, $\{H_n\}$, such that 
each $H_n=\<S_n\>$ contains a hyperbolic element on $X$ and 
that $\{e(H_n,S_n)\}$ is a strictly decreasing sequence with $\lim_{n \to \infty} e(H_n,S_n)=d$, for some $d\ge 1$.

By our assumption, Proposition \ref{bound.generators} applies to $H_n$.
Therefore, $|S_n|$ is uniformly bounded from above. By passing
to a subsequence we may assume that $|S_n|$ is constant, $|S_n|=\ell$, for the entire sequence. 
 
Let $S_n=\{x_1^{(n)}, \cdots, x_\ell^{(n)}\}$.
Let $F$ be the free group of rank $\ell$ with a free generating set:
$S=\{s_1, \ldots, s_\ell\}$.
For each index $n$, we define a map: $f_n:F \to G$, by setting: $f_n(s_i)=x_i^{(n)}$.
By construction: $e(H_n,S_n)=e(H_n,f_n(S))$.

%We denote $F_\ell$ as $F$. 
Then, as before, the sequence $\{f_n:F \to G\}$ subconverges to a surjective 
homomorphism  $\eta: F \to L$, where $L$ is a limit group over $G$.

By assumption, $G$ is equationally Noetherian.
By the general principle (Lemma \ref{basic}),
there exists an epimorphism $h_n: L \to G$ such that 
by passing to a subsequence we may assume that all the homomorphisms 
$\{f_n\}$ factor through the limit epimorphism: $\eta: F \to L$, ie, 
 $f_n=h_n \circ \eta$. 
  
  Then, we have

\begin{prop}[cf. Proposition \ref{1.2}]\label{6.3}
$$\lim_{n \to \infty} e(H_n,f_n(S))=e(L,\eta(S)).$$

\end{prop}

\begin{remark}
This proposition is a special case
of Proposition \ref{semi.cont.family}, when
all $\delta$-hyperbolic spaces $X_n$
are $X$.
\end{remark}

\proof 
The proof is identical to Proposition \ref{1.2}.
As in the beginning of the argument (Section \ref{section.separator}), for each $n$, we pick one hyperbolic isometry
$g\in S_n^M$ on $X$ that is $D$-WPD to start with. This is possible since $H_n=\<S_n\>$
contains a such hyperbolic element by the definition 
of $\Theta_X(G)$.
The rest is same, and we omit it. 
\qed

We continue the  proof of  the theorem. 
As in the proof of Theorem \ref{main.proof}, 
Proposition \ref{6.3} proves that there is no strictly decreasing sequence of rates of growth, $\{e(H_n,S_n)\}$, since a strictly decreasing sequence
can not approach its upper bound, a contradiction.
 Hence, $\Theta_X(G)$ is well-ordered.
 The theorem is proved. 
\qed

\subsection{Relatively hyperbolic groups}
We apply Theorem \ref{6.2} relatively hyperbolic groups. 
Let $(G,\{P_i \})$ be a relatively hyperbolic group with $P_i$
finitely generated.
Define
$$\Theta_{{\rm non-elem.}}(G)=\{e(H,S)|S\subset G, |S|< \infty,
H=\<S\>, e(H,S)>1\},$$
where in addition we only consider $H$ that is not conjugate into any $P_i$.
This is a subset in $\Theta(G)$.

We first characterize the subgroups $H$ that appear in the 
definition in terms of the action on a $\delta$-hyperbolic space $X$
that we described in Section \ref{section.5.2}.
Fix such $X$ and an action by $G$.

We recall a lemma. This is straightforward from 
the classification of subgroups that act on a hyperbolic space,
\cite[Section 3.1]{G}.
\begin{lemma}\label{subgroup.rel.hyp}
Let $H<G$ be a subgroup. Then the following two 
are equivalent:
\begin{enumerate}
\item
$H$ has an element $g$ that is hyperbolic on $X$, and
$H$ is not virtually $\Bbb Z$.
\item
$H$ is infinite, not virtually $\Bbb Z$ and not conjugate
into any $P_i$. 
\end{enumerate}
\end{lemma}
%This must be well-known but we give a proof.
%\proof
%(1) implies (2) is trivial since a hyperbolic element $g$
%has infinite order, and any conjugate of every $P_i$
%does not contain a hyperbolic element.
%
%To prove the other direction, assume (2).
%Since $H$ is infinite, and $X$ is proper, $H$ has
%non-empty limit set, $\Lambda \subset \partial X$.
%
%If $\Lambda$ consists of one point, $y \in \partial X$, then 
%$H$ is parabolic and $y$ is a parabolic point, so that 
%$y \in \Pi$. 
%But this implies that $H$ is conjugate into one of the $P_i$, impossible. 
%
%If $\Lambda$ consists of two points, then $H$ must 
%contain a hyperbolic element, $g$, and $\< g\> <H$ 
%has finite index, so that $H$ is virtually $\Bbb Z$, impossible. 
%
%We conclude that $\Lambda$ has at least three points, which 
%implies that $H$ is not virtually $\Bbb Z$, and that 
%$H$ contains a hyperbolic isometry. We are done. 
%\qed

The lemma implies:
$$\Theta_{{\rm non-elem.}}(G)=\Theta_X(G).$$

We prove:
\begin{thm}\label{thm.subgroup.rel.hyp}
Let $G$ be a group that is hyperbolic relative to a collection of 
subgroups $\{P_1, \cdots, 
P_n\}$.
Suppose $G$ is not virtually cyclic, and not equal to $P_i$
for any $i$.
Suppose each $P_i$ is finitely generated and equationally Noetherian.
Then $\Theta_{{\rm non-elem.}}(G)$ is well-ordered. 
\end{thm}
\proof
It suffices to argue that $\Theta_X(G)$ is well-ordered. 
For that, we apply Theorem \ref{6.2} to the action of $G$ on $X$.
We already checked the assumptions in the proof of Theorem \ref{rel.hyp.growth}, except that Lemma \ref{Xie} holds 
for subgroups. Namely, there is a constant $M$ such that for any finite 
set $S\subset G$ such that $\<S\>$ contains a hyperbolic 
isometry on $X$, then $S^M$ contains a hyperbolic isometry.
But this is true and the argument is identical to the proof
of \cite[Proposition 5.1]{X}, and we omit it. 
\qed

As an example of Theorem \ref{thm.subgroup.rel.hyp}
we prove:

\begin{thm}\label{6.1}
Let $G$ be a group in Theorem \ref{lattice}.
Then $\Theta(G)$ is a well-ordered set. 

\end{thm}

%We prove Theorem \ref{6.1}.

\proof
As we said in the proof of Theorem \ref{lattice}, $G$ is relatively 
hyperbolic w.r.t. the parabolic subgroups $\{H_i\}$, which are 
associated to the cusps, and $H_i$ are virtually nilpotent. 
By applying Theorem \ref{thm.subgroup.rel.hyp}
to $(G,\{H_i\})$, we have that 
$\Theta_{{\rm non-elem}}(G)$ is well-ordered.
But if a subgroup $H=\<S\>$ is conjugate into one of $H_i$, then 
it is virtually nilpotent, so that $H$ has polynomial growth.
It follows that $\Theta_{{\rm non-elem}}=\Theta(G)$ holds for $G$,
so that $\Theta(G)$ is well-ordered.
\qed

%
%\proof
%We apply  Proposition \ref{prop.rank1} to $G$.
%Let $X$ be as in the proof of the proposition. 
%We claim
%$$\Theta(G)=\Theta_X(G).$$
%It suffices to show that if $S\subset G$ is a finite set and 
%$H=\<S\>$ has exponential growth, then $H$ contains
%a hyperbolic element on $X$.
%But as we said in that proof, by \cite[Theorem 1.13]{BrF}, 
%either $S^M$, which is a subset of $H$, contains
% a hyperbolic element, or $H=\< S\>$ is
% virtually nilpotent, which does not happen since $H$
% has exponential growth. 
% We showed $\Theta(G)=\Theta_X(G)$, which 
% implies that $\Theta(G)$ is well-ordered by Theorem \ref{6.2}.
% \qed

\subsection{Subgroups in mapping class groups}
We discuss mapping class groups.  A subgroup $H<MCG(\Sigma)$
is called {\it large} if it contains two independent 
pseudo-Anosov elements. Such $H$ has exponential growth.

We define:
$$\Theta_{{\rm large}}(MCG(\Sigma))=\{e(H,S)| S \subset MCG(\Sigma), |S| < \infty, \<S\>=H, e(H,S)>1\},$$
where in addition $H < MCG(\Sigma)$ is large. 
Note that $\xi(MCG(\Sigma)) \subset \Theta_{{\rm large}}(MCG(\Sigma))$.
\begin{thm}\label{thm.6.5}
Assume that $G=MCG(\Sigma)$ is equationally Noetherian. Then $\Theta_{{\rm large}}(MCG(\Sigma))$
is well-ordered.

\end{thm}

\proof
We suppress $\Sigma$ and denote $MCG$.
Let $X$ the curve graph of $\Sigma$.
Then, as we said,  Theorem \ref{6.2} applies to the action of 
$MCG$ on $X$.
We conclude that $\Theta_X(MCG)$ is well-ordered.
Now, we claim $\Theta_X(MCG)=\Theta_{{\rm large}}(MCG)$.
Indeed, given $S\subset MCG$, if $H=\<S\>$
contains a hyperbolic element on $X$, then it is a pseudo-Anosov
element, and moreover, from $e(H,S)>1$, $H$ must be large.
We showed $\Theta_X(MCG) \subset \Theta_{{\rm large}}(MCG)$.
On the other hand, for $S\subset MCG$ if $H=\<S\>$
is large in $MCG$, then $H$ contains hyperbolic isometries on $X$, so that 
$\Theta_{{\rm large}}(MCG) \subset \Theta_X(MCG)$.
\qed

It is natural to ask the following question. 
To deal with a non-large subgroup, considering the action 
on the curve graph does not seem to be enough. 
\begin{question}
Is $\Theta(MCG(\Sigma))$ well-ordered ?
\end{question}

\section{Finiteness}
\subsection{Finiteness of equal growth generating sets}
If $G$ is a hyperbolic group, it is known by \cite[Section 3]{FS}
that for $\rho \in \xi(G)$, there are only finitely many 
generating sets $S$ of $G$, up to $\Aut(G)$, such that 
$\rho=e(G, S)$.
We discuss this issue.

\begin{thm}[Finiteness. cf. Theorem 3.1 in \cite{FS}]\label{7.1}
Suppose a finitely generated group $G$ acts on a $\delta$-hyperbolic space $X$ and 
$G$ is not virtually cyclic.
Let $D(\epsilon)$ be a WPD-function. 
Assume that there exists a constant $M$ such that if $S$ is a finite generating set of $G$, then $S^M$ contains a hyperbolic 
element that is $D$-WPD.
Assume that $G$ is equationally 
Noetherian.

Then for any $\rho \in \xi(G)$, up to the action on $\Aut (G)$,
there are at most finitely many finite generating set $S$
such that $e(G,S)=\rho$.
\end{thm}

\proof  We argue by contradiction. Suppose that there are infinitely many finite sets of generators $\{S_n\}$ that satisfy: $e(G,S_n)=\rho$, and no pair of generating sets $S_n$ is equivalent 
under the action of the automorphism group $\Aut(G)$. As in the proof of Theorem \ref{main.proof}, by Proposition \ref{bound.generators},
 the cardinality of 
 the 
generating sets $\{S_n\}$ is bounded, so we may pass to a subsequence that have a fixed cardinality $\ell$. Hence, each generating set $S_n$ corresponds
to an epimorphism, $f_n:F \to G$, where $S$ is a fixed free generating set of $F$, and $f_n(S)=S_n$.

By passing to a further subsequence, we may assume that the sequence of epimorphisms $\{f_n\}$ converges to a limit group $L$ with $\eta:F \to L$ the associated  quotient map. 
As in the proof of Theorem \ref{main.proof}, since $G$ is equationally Noetherian, by Lemma \ref{basic},
for large $n$, $f_n=h_n \circ \eta$, where $h_n:L \to G$ is an epimorphism.
In particular, $S_n=h_n (\eta(S))$. We pass to a further subsequence such that for every $n$, $f_n=h_n \circ \eta$.

Since for every index $n$, $h_n$ is an epimorphism from $L$ onto $G$ that maps $\eta(S)$ to $f_n(S)$, $e(G,f_n(S)) \leq e(L,\eta(S))$. 
We prove:

\begin{prop} [cf. Proposition 3.2 \cite{FS}]\label{2.2}
%Assume in addition that $G$ does not contain a non-trivial, finite 
%normal subgroup.
If  $\ker(h_{n_0})$ is infinite for some $n_0$, then 
%an epimorphism $h_n$ is not an isomorphism, then:
%Then for this $n$,  the generating sets $\{g_n(S)\}$ of $G$ (from the remaining subsequence that factor through the limit
%group $L$)  satisfy: 
$e(G,f_{n_0}(S))<e(L,\eta(S))$. 
\end{prop}

We postpone the proof of the proposition until the next
section, and proceed. 
We prove a lemma.
\begin{lemma}\label{radical.limit.group}
The group $L$ contains a finite normal subgroup $N=N_L$ that 
contains all finite normal subgroups in $L$, such that $|N|\le 2D(100\delta)$.

\end{lemma}
We recall one fact we use in the proof. 
If a finitely generated group $G$ acts on a $\delta$-hyperbolic space $X$  such that $G$
is not virtually cyclic and $G$ contains
a hyperbolic element on $X$ that is $D$-WPD,
then
$G$ contains a maximal finite normal subgroup $N<G$. Moreover, $|N| \le 2 D(100\delta)$. 
We sometimes denote $N$ by $N_G$. 

The existence of such $N$ is known for an acylindrically hyperbolic group \cite[Theorem 6.14]{DGO}, 
and the same proof applies to our setting, which we briefly recall. Indeed, $N$ is the intersection of $E(g)$ for all hyperbolic elements $g\in G$
on $X$. It is obvious that $N$ is normal.
By assumption, there must be a hyperbolic
and WPD element, $g$. Also, there is another 
element $h$ such that $g$ and $h$ are 
independent. Then, by Proposition \ref{prop.elementary}, $E(g) \cap E(h)$ is finite, 
so that $N$ is finite. 
On the other hand, if $N'$ is a finite normal
subgroup in $G$, then for every hyperbolic
element $g \in G$,  there is $n>0$
such that $N'$ is contained in the centralizer
of $g^n$, so that $N' < E(g)$.
It implies that $N' <N$. We showed that $N$  is maximal.

Lastly, to see $|N|\le 2D(100\delta)$, 
consider the exact sequence
$1 \to F(g) \to  E(g) \to C \to 1$
for the hyperbolic and $D$-WPD element $g$.
Recall that $|F(g)| \le D(100\delta)$
if  $C$ is cyclic. From this we have
$|N|\le 2D(100\delta)$.

We prove the lemma. 
\proof
Let $N<L$ be a finite normal subgroup. 
Since $h_n$ is surjective, $h_n(N) <G$
is a finite normal subgroup, therefore $h_n(N) < N_G$
 for any $h_n$.
 Also, for sufficiently large $n$, the surjection $h_n:L \to G$
is injective on $N$. 
But since $|N_G| \le 2D(100\delta)$, we have $|N|\le 2D(100\delta)$. 

If $N_1, N_2 < L$ are two finite normal subgroups, 
then $N_1N_2 $ is a finite normal subgroup. Combined with 
the fact in the previous paragraph, there must be the maximal 
finite normal subgroup $N_L$ in $L$ with $|N_L| \le 2D(100\delta)$. 
\qed

By Proposition \ref{1.2}, $\lim_{n \to \infty} e(G,S_n)=e(L,\eta(S))$. By our assumption, for every index $n$, $e(G,S_n)=\rho$. Hence, $e(L,\eta(S))=\rho$, so that for every $n$, $e(G,S_n)=e(L,\eta(S))$.

It follows from Proposition \ref{2.2} that 
for every $n$, $\ker(h_n)$ is finite.
Since $\ker(h_n)$ is a normal subgroup in $L$, 
by Lemma \ref{radical.limit.group}, $\ker(h_n) < N_L$. 
Since $N_L$ is a finite group, there are only finitely many 
possibilities for $\ker(h_n)$.
It follows that there must be $N_0 < N_L$ 
such that $\ker(h_n)=N_0$ for infinitely many $n$.

The map $h_n$ induces an isomorphism from $L/\ker(h_n)$
to $G$. Notice that this gives an isomorphism
from $(L/\ker(h_n), \eta(S))$ to $(G,S_n)$ since 
$h_n$ gives a bijection between $\eta(S)$ and $S_n$. 
(Here, we may assume that each $S_n$ consists of distinct 
elements, so that no two elements in $\eta(S)$ are identifies
by $h_n$.)
But this implies that  $(L/N_0, \eta(S))$
is isomorphic to $(G,S_n)$ for infinitely many
$n$ by $h_n$, ie, those $(G,S_n)$ are isomorphic to 
each other. This is a contradiction since all of them
must be non-isomorphic. 
Theorem \ref{7.1} is proved.
\qed

\subsection{Idea of the proof of Proposition \ref{2.2}}
We prove Proposition \ref{2.2}.
The argument is long and complicated, but the main idea is same as the proof of \cite[Proposition 3.2]{FS},
and we adapt it to our setting. 
Also, the proof is similar to the proof of Proposition \ref{1.2}, which 
also follows the counterpart in the paper \cite{FS}.
The difference between this paper and \cite{FS} is that 
while they use the action of the limit group $L$ on a limit object,
called a limit tree, while in our paper we use the actions of $L$
on $X$ induced from the maps $h_n:L \to G$. But this approach is already  taken 
in the proof of Proposition \ref{1.2}.

So, rather than giving a full formal proof, 
we first explain the strategy of the proof, then 
give all definitions and intermediate claims, which 
appear in the proof of \cite[Proposition 3.2]{FS},
then explain the part where we need to make
technical modifications, most of which already appeared in 
Section \ref{section.4}.
One advantage of not using the action of $L$ on a limit object is 
that one does not need to deal with the degeneration of the action
on the limit object. A trade-off  is that we need to keep attention 
to the various constants related to the actions induced by $h_n$
through the argument.

{\it Strategy of the proof}.
We start with an informal description of the idea. The constant $n_0$
is given as the assumption, which gives the homomorphism 
$h_{n_0}:(L,\eta(S)) \to (G,f_{n_0}(S))$ with an infinite kernel.
To show that $e(G,f_{n_0}(S)) < e(L,\eta(S))$, 
we will produce not only infinitely many (which is obvious
by the assumption), but ``exponentially many'' elements
in the preimage of $g$ by $h_{n_0}$ for each $g \in G$.
Those elements are given by the map
$\phi_n$. 
They are exponentially many in terms of the word length 
of $g$ w.r.t. $f_{n_0}(S)$. See the estimate (\ref{ad3}).

In the proof two constants $m$ and $n$ appear.
They will be chosen and fixed  around the end of the proof. The constant $m>1$
is first chosen. It will be used to measure the gap
between $e(G,f_{n_0}(S))$ and $e(L,\eta(S))$. The constant $n$ depends on $m$, so that  $h_n$ is injective on 
$B_{2m}(L,\eta(S))$.
Also, a positive integer  $q$ is used to make the word length of 
$g$, which is $mq$ longer and take a limit at the end of the proof. 

We now explain more concretely. 
Similar to $B_m(G,S_{n})$, the ball of radius $m$ in $Cayley(G,S_{n})$
centered at the identity, 
let $\Sph_m(G,S_n)$ denote the sphere of radius $m$.
It is an elementary fact that if $e(G,S_n) >1$, then 
$$e(G,S_n) = \limsup_{m\to \infty} |\Sph_m(G,S_n)|^{\frac{1}{m}}.$$

Given $m>0$, for a large enough $n>0$ depending on $m$, 
we will define a ``map'' $\phi_n$:
$$\phi_n: \Sph_{mq}(G,S_{n_0}) \to B_{q(m+2b)}(L,\eta(S))$$
for all $q>0$, 
where $b$ is a constant that does not depend on $n,m,q$. 
The map $\phi_n$ is similar to  the map $\Phi_n$ in the section \ref{section.forbidden}, but 
strictly speaking $\phi_n$ is not a map, 
but $\phi_n(g)$ is a finite set of elements in $B_{q(m+2b)}(L,\eta(S))$
for each $g$. But we abuse the notation and call them maps
in the following account. 

It satisfies for every $g \in \Sph_{mq}(G,S_{n_0})$ to have the following two properties.
Set $\D=D(200\delta)$.
\begin{enumerate}
\item[(i).]
$$\left( \frac{m-1}{\D^4(\D-1)}\right)^q \le |\phi_n(g)|, $$
see the estimate (\ref{ad2}). 
\item[(ii).]
$h_{n_0} \circ \phi_n (g)=g,$
which implies that for $g\not=h$, we have $\phi_n(g) \cap \phi_n(h) =\emptyset$.
\end{enumerate}

Once we have such a map $\phi_n$, we argue as follows:
fix a (large) $m$.
% and let us pretend for all $q$ (in fact we only have it 
%asymptotically with $q \to \infty$):
%$$\frac{1}{mq} \log(|\Sph_{mq}(G,S_n)|)=\log e(G,S_n).$$
%On the other hand,
 Since $\phi_n(g) \subset B_{q(m+2b)}(L,\eta(S))$, we have from (i) and (ii) that 
$$|\Sph_{mq}(G,S_{n_0})| \left( \frac{m-1}{\D^4(\D-1)}\right)^q  \le |B_{q(m+2b)}(L,\eta(S))|.$$
Taking $\log$, dividing by $mq$, and letting $q\to \infty$, 
we have as the limsup
$$\log e(G,S_{n_0}) + \frac{\log (m-1)-\log(\D^4(\D-1))}{m} \le \frac{m+2b}{m}
\log e(L,\eta(S)).$$
Since $b$ does not depend on $m$, choosing $m$ large
enough, this inequality implies
$$\log e(G,S_{n_0}) < \log e(L,\eta(S)).$$

Roughly speaking, the construction of $\phi_n$
is as follows. 
As in the construction of $\Phi_n$, we first construct separators. 
To define separators in $L$, we use a non-trivial element $r_n \in \ker h_{n_0}$.
Separators will be products of conjugates of $r_n$, 
so that they are also in $\ker h_{n_0}$, which will 
imply the property (ii) in the above. 
The separators depend on $n$. 
%For each $g$, the set $\phi_n(g)$ are actually $(m-1)^q$ words in $(L,\eta(S))$, 
%so that some of the $(m-1)^q$ elements in $L$ we get may be
%same as elements in $L$. Nevertheless, we prove that 
%$\phi(g)$ gives at least $\left(\frac{m-1}{D^4(D-1)}\right)^q$ distinct elements in $L$,
%which are called {\it feasible elements}, 
%which is enough for our purpose if we modify the above estimates.

For each $n$, the map $h_n:(L,\eta(S)) \to (G,S_n)$ gives a canonical bijection between $\eta(S)$ and $S_n=f_n(S_n)$. This gives a bijection 
between the words on $\eta(S)$ (not elements in $L$)  and the words on $S_n$.

Let $m,q>0$ be integers. We will fix $m$ and let $q \to \infty$ later. 
Given an element $g \in \Sph_{mq}(G,S_{n_0})$,
we choose a word $w(g)$ of length $mq$ on $S_{n_0}$
that represents  $g$. We divide
$w(g)$ into $q$ subwords of length $m$.
As we said, each subword of length $m$ canonically gives a word of 
length $m$ on $\eta(S)$ via the map $h_{n_0}$. 
We further subdivide each of the subwords of length $m$ on $\eta(S)$ into two
words of length $k$ and $m-k$. We choose $k$ to satisfy 
$1 \le k \le m-1$.
In this way, for each choice of  a $q$-tuple of such $k$'s, 
we divided the word on $\eta(S)$ corresponding to $w(g)$ into $2q$ subwords. There are $(m-1)^q$ ways
to subdivide it. 

To each of such subdivision, we insert separators to the $(2q-1)$ break points
and obtain an element in $B_{q(m+2b)}(L,\eta(S))$
since the word length of each separator is at most $b$. We obtain 
$(m-1)^q$ such elements.
Since separators are in $\ker(h_{n_0})$, 
those elements (words) are mapped to $g$ by $h_{n_0}$. 

But we do not know if they are all distinct as elements in $L$, but
we will show that there are at least $(\frac{m-1}{\D^4(\D-1)})^q$
elements that are distinct. They are called {\it feasible elements}. 
This collection of feasible elements
is denoted by $\phi_n(g)$.
They are significantly many, so that $e(L,\eta(S))$
is strictly larger than $e(G,S_{n_0})$ as we computed in the above.

Lastly, to show that those feasible elements in $L$ defined for each $g$ are 
distinct, as in Section 4, we let them act on the space $X$ via 
the map $h_n$ for a large enough $n$ (We choose $n$ 
such that  $h_n$ is injective on $B_{2m}(L,\eta(S))$.) We then argue that 
the images of the base point $y_n \in X$ by those elements are distinct.
In the paper \cite{FS}, 
they use a limit tree $Y$ on which $L$ acts to argue that 
feasible elements are distinct. Here, we use the action on $X$.
This difference 
already appeared in Section 4.

\subsection{Proof}
We prove Proposition \ref{2.2}.
%Note that the action of $G$ on $X$ is faithful. Indeed, let 
%$N<G$ be the kernel of the action. Then uniform WPD implies that 
%$N$ is finite. Then it must be trivial since $N$ is normal. 
\proof

As in Section 4.1, set $\D=D(200\delta)$, and for each $n$  let $y_n \in X$ be a point where $L(S_n^{2\D M})$
is achieved and put
$$\Delta_n=100\delta + 4\D L(S_n^{2\D M}).$$
We consider germs w.r.t. $\Delta_n$.

As we did in Lemma \ref{4.2}, in the next lemma we construct separators 
$u_i \in G$, which give 
$\hat u_1, \hat u_2, \hat u_3, \hat u_4 \in L$
by pulling them back by $h_n$. 
In addition to the properties in Lemma \ref{4.2} they satisfy 
$h_{n_0}(\hat u_i)=1$ in $G$.
We state it as a lemma. The constant $b$ is different from the constant $b$ in Lemma \ref{4.2}, and it depends on $\ker(h_{n_0})$
but not on $n$.

\begin{lemma}[cf Lemma 3.3 and Lemma 3.9 in \cite{FS}]\label{lemma7.4}
Suppose $\ker(h_{n_0})$ is infinite. 
Then there exists a constant $b$ with the following property:
if $n$ is sufficiently large, then there are elements $u_1, u_2, u_3, u_4 
\in S_n^b$ that satisfy the conditions (i), (ii) and (iii)  in 
Lemma \ref{4.2}, and in addition to that, 
all $u_i$ satisfy
\begin{itemize}
\item[(iv)] $u_i \in h_n(\ker h_{n_0})$.
\end{itemize}
Moreover, those elements are such that there are elements
$\hat u_1, \hat u_2, \hat u_3, \hat u_4 \in \eta(S)^b$ with
$h_n(\hat u_i)=u_i$ and $u_i \in \ker (h_{n_0})$ for each $i$.

\end{lemma}

We remark that the moreover part is rather than an additional 
property, but immediate from the construction, ie, we first 
construct $\hat u_i$ then map them by $h_n$. 
We will call both $u_i$ and $\hat u_i$ separators. 

\proof
In the proof there will be several constants $b_i$, for which we do not try to give explicit  values. The important property of those constants is that they do not depend on $n$.

First, since $\ker(h_{n_0})$ is infinite, choose
distinct elements $r_1, \cdots, r_{\D+1} \in L$ that are in the kernel.
Let $b_1$ be the maximum of the word lengths of the $r_i$ in terms of $\eta(S)$.
If $n$ is large enough, then the image of those $\D+1$ elements by $h_n$
are all distinct. From now on, we only consider such $n$.
In the following we fix each such $n$ and argue. 
We have $\D+1$ distinct elements $\{h_n(r_i)\}$ and the word lengths of those w.r.t. $S_n$ are bounded by $b_1$.

Secondly, choose an element $g_n \in S_n^M$
that is hyperbolic on $X$ such that its $10\delta$-axis is at distance at most $10\delta$ from 
the point $y_n$. 
By assumption, such an element exists. 
Also there is $s \in S_n$ such that $g_n$ and $sg_ns^{-1}$ are independent. 
Choose such $s$. 

Thirdly, we choose one element, $r$, from the $r_i$'s as follows: if there is $r_i$ such that \\
(I) The element $h_n(r_i)$ is hyperbolic on $X$, 
\\
then choose one of such $r_i$ and set $r=r_i$. Otherwise choose $r_i$ with \\
(II) $h_n(r_i) \not\in F(g_n)$
\\
and set $r=r_i$.
This is clearly possible since $|F(g_n)| \le \D$. 

Note that the element $r$ depends on $n$.
From now on we suppress $n$ and write $g_n$ as $g$.

Now we divide the case into two depending on (I) or (II) in the above.
Suppose we are in the case (I). 
We consider the power $g^k$ with $k =60\D$, then we have (see Lemma \ref{2.4})
$$\<g^k,sg^ks^{-1}\>=\<g^k\> * \<sg^k s^{-1}\>.$$
Note that we have 
$$\lambda(g^k)  \le \Delta_n-100\delta.$$
This is because since $k=60\D$ and $g \in S_n^M$,
we have 
$$\lambda(g^k) =30 \lambda(g^{2\D}) \le 30 L(S_n^{2\D M}) \le \Delta_n-100\delta.$$ 
The last inequality is by 
$\D \ge 10$ and $\Delta_n=100\delta + 4\D L(S_n^{2\D M}).$

Recall that the axes of $g^k, sg^ks^{-1}$
are at at most $40\delta + L(S_n^{2M\D})$ from $y_n$.

In the proof of Lemma \ref{4.2}, we set $w=g^k, z=sg^ks^{-1}$
and produce $u_i$ as words on $w, z$,  but 
this time we take into account the germs of the element $h_n(r)$
and choose $z,w \in \<g^k\> * \<sg^k s^{-1}\>$
as follows.

Notice that 6 elements $g^k,g^{-k},sg^ks^{-1}, sg^{-k}s^{-1},
g^ksg^k s^{-1}, g^ksg^{-k} s^{-1},$ define six germs at $y_n$
that are mutually opposite since $\lambda(g^k) \le \Delta_n-100\delta$
as we noted. 
%(We use this inequality to deal with the last two germs.)
From the six, choose four distinct germs that are opposite to the germs
for $h_n(r)$ and $h_n(r)^{-1}$. If those two germs are empty, 
then ignore this condition. 
Denote those four germs as $\gamma_1, \gamma_2, \gamma_3, \gamma_4$.
Now choose $w, z \in \<g^k,sg^ks^{-1}\>$ such that 
the germ for $w, w^{-1}, z, z^{-1}$ is equivalent to $\gamma_1, \gamma_3, \gamma_2, \gamma_4$,
respectively,
such that the axes of $w,z$ are at at most $60\delta+L(S_n^{2\D M})$ from $y_n$.

We also arrange that the axes of $w,z$ are not parallel to each other, and no element of $G$ flips the axes
of $w$ or $z$. This is achieved by a similar technique to the one we used
to prove Lemma \ref{2.4}, so we do not repeat. 
As in the proof of Lemma \ref{4.2}, there exists a constant $b_2$ that depends only on $\delta, M$
and  $\D=D(200\delta)$ such that the word lengths of $w,z$ w.r.t. $S_n$
are bounded by $b_2$. The constant $b_2$ does not depend on $n$, nor  the 
choice of $r$.

Now, choose $\hat z, \hat w \in L$ with $h_n(\hat z)=z, h_n(\hat w)=w$,
whose word length w.r.t. $\eta(S)$ are also bounded by $b_2$. 
Using $\hat z, \hat w,r$ we define $\hat u_i \in L$
as follows: 

\begin{align*}
\hat u_1 & = \hat w r\hat w  ^{-1} \cdot \hat zr\hat z^{-1} \cdot \hat w^2 r\hat w^{-2}  \cdot  \hat zr\hat z^{-1} \cdots \hat w^{19} r\hat w^{-19
}  \cdot  \hat  zr\hat  z^{-1}  \cdot  \hat w^{-20} r\hat w^{20}, \\
\hat u_2 & = \hat  w^{21}r\hat  w^{-21}  \cdot \hat  z^{-1}r\hat  z  \cdot \hat  w^{22} r\hat  w^{-22}  \cdot \hat  z^{-1}r\hat  z \cdots    \hat  w^{40} r\hat  w^{-40} \cdot z^{-1}r\hat  z, \\
\hat u_3 & = \hat z r \hat  z^{-1}  \cdot \hat  w^{-41}r\hat  w^{41}  \cdot  \hat  zr\hat  z^{-1}  \cdot  \hat  w^{-42} r\hat  w^{42}  \cdot  \hat  zr\hat  z^{-1} \cdots   \hat  zr\hat  z^{-1}  \cdot  \hat  w^{-60} r\hat  w^{60}, \\
\hat u_4& = \hat z r \hat  z^{-1}  \cdot \hat  w^{61}r\hat  w^{-61}  \cdot  \hat  zr\hat  z^{-1}  \cdot  \hat  w^{62} r\hat  w^{-62}  \cdots  \hat  z r\hat  z^{-1}  \cdot  \hat  w^{80} r\hat  w^{-80} \cdot  \hat  z^{-1}r\hat  z. 
\end{align*}

The word length of $\hat u_i$ in terms of $\eta(S)$ is at most $b_3$,
which does not depend on $n$.
Indeed, we may set $b_3=21 b_1 + 1420 b_2$. 

Clearly, $\hat u_i$ are in $\ker h_{n_0}$ since 
they are products of conjugates of $r \in \ker h_{n_0}$. 
Finally, define $u_i=h_n(\hat u_i) \in S_n^{b_3}$.
Then, they satisfy the property (iv). 

We need to check those elements
satisfy the other properties, (i), (ii) and (iii) in Lemma \ref{4.2}.
Regarding (i), the germ $\gamma_1$ is the germ of $w$, 
the germ $\gamma_3$ is the germ of $w^{-1}$, 
the germ $\gamma_2$ is the germ of $z$, and 
the germ $\gamma_4$ is the germ of $z^{-1}$. 
Then $u_1,u_2,u_3,u_4$ satisfy (i). 
The property (iii) is a consequence of that the axes of $w$ and $z$
are not parallel to each other. 
We skip details of the arguments since it is similar to Lemma \ref{4.2}.
We  point out that some of the argument slightly differs 
depending if the germs for $h_n(r), h_n(r)^{-1}$
are defined or empty. 
In the empty case, we use the property that 
the axes of $w, z$ are not flipped by any element
of $G$.

In conclusion,  those are desired elements,
so in this case we take $b=b_3$ and we are done.

Suppose we are in the case (II) when we chose $r$. 
In this case, we replace $r$ with another element $r'$
such that $h_n(r')$ is hyperbolic on $X$. 
We explain how we produce such an $r'$. 
 First choose an element $\hat g \in L$ with 
$h_n(\hat g)=g$, where the word length of $\hat g$ is at most $M$ in terms
of $\eta(S)$. 
Then we consider an element of the form:
$$r' =r \hat g^Q r \hat g^{-Q}.$$
We will show that if $Q \ge 40\D $, 
$h_n(r')$ is hyperbolic. See Lemma \ref{lemma.Q}.
Also, since $r \in \ker h_{n_0}$, we have 
$r' \in \ker h_{n_0}$.
%and $L(h_n(r')) \ge 100 \delta$.
%Moreover we may assume that $Q$
%is bounded by a constant, $b_4$, uniformly over all $n$. 
It is a well-known method to produce a hyperbolic element
as a product of two non-hyperbolic elements, and
we postpone an explanation on this, and 
proceed. 

But, then the word length of $r'$ in terms of $\eta(S)$ is at most $b_5$, where $b_5=2b_1+60\D M$, which 
does not depend on $n$.

With the element $r'$ we repeat the same argument as we did for $r$
in the case (I) and obtain desired $u_i$ with a bound on the word length, 
uniformly over all $n$, which will 
finish the proof of Lemma \ref{lemma7.4}.

We now explain some details on how to produce $r'$ from $r$.
For an element $k \in Isom(X)$, define a set
$$Min(k)=\{x \in X| |x-k(x)| \le 100 \delta\}.$$
This is a $k$-invariant set. 
We state a few standard facts on this set.
\begin{enumerate}
\item[(M1)]
If $k$ is not hyperbolic then $Min(k)$  is not empty, 
since if it is empty, then $L(k) > 100 \delta$, which implies $k$ is hyperbolic, 
a contradiction. 
\item[(M2)]
If $k$ is not hyperbolic, then 
for any point $y \in X$, we have
$$|y-k(y)| \ge 2(d(y,Min(k))- 400 \delta).$$
%\item
%The set $Min(k)$ is a $500 \delta$-quasi-convex set in the sense that a geodesic joining any two points in the set is in the $500 \delta$-neighborhood of the set. 
\end{enumerate}

We prove (M2). 
If $d(y,Min(k)) \le 400\delta$, then nothing to show
so suppose $d(y,Min(k)) > 400\delta$.
Let $x\in Min(k)$ be a point with $d(x,y)=d(Min(k),y)$.

Let $z \in [x,y]$ be the point with $|x-z|=350\delta$.
We claim that $z \not\in N_{10\delta}([k(x),k(y)])$
and $k(z) \not\in N_{10\delta}([x,y])$.
To prove the first claim, suppose not, ie,
$z \in N_{10\delta}([k(x),k(y)])$.
Let $v,w \in [x,y]$ be with $|x-v|=150\delta$
and $|x-w|=200\delta$.
Then we have $k([v,w]) \subset N_{10\delta}
([x,y])$ since $|x-k(x)| \le 100 \delta$
and $z \in N_{10\delta}([k(x),k(y)])$.
Then we have $|v-k(v)| \le 50\delta$
since $k$ is not hyperbolic.
(Otherwise, $v$ is ``pushed'' along the geodesic $[x,y]$ by at least $40\delta$ by $k$, which 
implies that $k$ is hyperbolic, impossible.)
But it implies that $v \in Min(k)$, which is a 
contradiction.
We showed $z \not\in N_{10\delta}([k(x),k(y)])$.

By the same argument, we can show 
$k(z) \not\in N_{10\delta}([x,y])$.

Having those two claims, we have
$$|y-k(y)| \ge |y-z| + |z-k(z)| + |k(z)-k(y)|
- 100\delta \ge 2(|y-x|-350\delta) - 100\delta
=2(|y-x|-400\delta).$$
We showed (M2).

%Let $z\in [x,y]$ be the last point in $N_{20\delta}([x,k(x)]\cup [k(x),k(z)])$
%when moving from $x$ to $y$. Then $|x-z| \le 300 \delta$,
%since if not, then 
%$z \in N_{20\delta}([k(x),k(z)])$.
%Take $w\in [x,z]$ wtih $|w-z|=150\delta$.
%Then $[w,z] \subset N_{30\delta}([k(x),k(y])$.
%But since $|z-k(z)| > 100 \delta$ (since $z$ is not in the set), 
%we find that $k$ is hyperbolic ($[w,z]$ is ``pushed'' along
%$[k(x),k(y)]$), a contradiction.
%We showed $|x-z| \le 300 \delta$.
%On the other hand, by the way  we chose $z$,
%we have $$|y-k(y)| \ge |y-z| +|z-k(z)| +|k(z)-k(y)| - 100 \delta
%\ge 2(|y-x|-300\delta) - 100\delta.$$ We showed (3). 

\begin{figure}[htbp]
\hspace*{-3.3cm}     
\begin{center}
                                                      
   \includegraphics[scale=0.6]{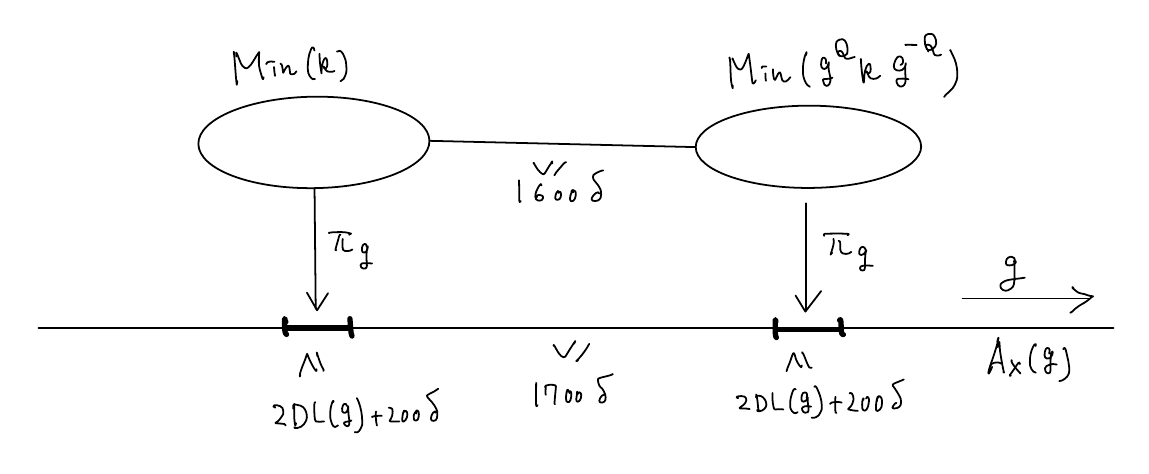}%
   \end{center}
\caption{The Min sets and the projection to an axis}

\label{fig.min}
\end{figure}

We go back to the explanation. 
For a hyperbolic isometry $g$ and its $10\delta$-axis $Ax(g)$, we consider the nearest points projection in $X$ to $Ax(g)$.
We denote the projection by $\pi_{g}$.
For every point $x\in X$, although $\pi_g(x)$
is not a point, 
the diameter of $\pi_g(x)$ is bounded by $100 \delta$.
The following lemma is well-known. See Figure \ref{fig.min}.
\begin{lemma}[Bounded projection]\label{lemma.projection}
If $g\in G$ is hyperbolic and $D$-WPD; and $k \in G \backslash E(g)$ and $k$ is not 
hyperbolic, then the image 
of $Min(k)$ in $Ax(g)$ by the projection $\pi_g$ is bounded by $2 D(100\delta) L(g) + 200 \delta$
in diameter. 
\end{lemma}

We  prove the lemma for readers' convenience. 
\proof
Let $x,y \in Min(k)$ and suppose 
$p \in \pi_g(x), q \in \pi_g(y)$.
Assume that $|p-q|$ is larger than $200 \delta$,
since otherwise there is nothing to show. 
Then $[x,p] \cup [p,q] \cup [q,y]$ is a uniform 
quasi-geodesic and every point on it is moved by $k$
by at most $200 \delta$.
Moreover, each point $x$ on $[p,q]$ with $|x-p|, |x-q| \ge 50 \delta$ is 
moved by $k$ by at most $20\delta$ since 
$k$ is not hyperbolic. 

Now, consider the points $p',q' \in [p,q]$
with $|p-p'|=|q-q'|=50\delta$.
Then since $g$ is $D$-WPD, we have $|p'-q'| \le 2D(100\delta) L(g) +100\delta$
by Lemma \ref{2.1} (2). This is because, 
otherwise, $k \in E(g)$, impossible. 
It follows that  
$|p-q| \le 2D(100\delta) L(g) +200\delta$.
\qed

The following lemma is also standard, and this is what we need for our purpose. See Figure \ref{fig.min}.
\begin{lemma}[Producing hyperbolic element]\label{lemma.Q}
For $g,k$ as in the  lemma \ref{lemma.projection}, the element 
$k g^Q k  g^{-Q}$ is hyperbolic
if $Q\ge 40D(100\delta)$.
\end{lemma}
We also give a brief proof. 
\proof
Set $\D=D(100\delta)$. (Or one can set $\D=D(200\delta)$ as usual. 
It does not matter since $D(100\delta) \le D(200\delta)$.)
In general, $Ax(g)$ is not exactly $g$-invariant, 
but if $L(g) \ge 10\delta$ by definition.
We also know that $L(g) \ge 50\delta/\D$
by Lemma \ref{lower.bound} (1) for our $g$.
So, if necessary, by replacing $g$ by $g^\D$,
we may assume that $Ax(g)$ is $g$-invariant.
For simplicity, in the following argument, we assume
that $Ax(g)$ is $g$-invariant. 

Consider the set $Min(g^Q kg^{-Q})$, which 
is equal to the set $g^Q (Min(k))$.
We consider the projection of those two sets
by $\pi_g$. Then since $Ax(g)$ is $g$-invariant, 
the projection $\pi_g$ is $g$-equivariant.
It implies that
$\pi_g(Min(g^Q kg^{-Q}))=  g^Q (\pi_g (Min(k)))$.

Then by Lemma \ref{lemma.projection}, 
the distance between $\pi_g(Min(g^Q kg^{-Q}))$
and $\pi_g (Min(k)))$ is 
at least
$$Q L(g) - (2 \D L(g) + 200 \delta)
\ge 38 \D L(g) - 200\delta \ge 1700 \delta$$
since $L(g) \ge 50\delta/\D$.
It follows that  the distance between $Min(g^Q kg^{-Q})$
and $Min(k)$ is at least $1600 \delta$.
It follows that the product of $k$ and $g^Q k g^{-Q}$ 
is hyperbolic (this is a well-known fact in $\delta$-hyperbolic geometry, ie, if the distance between 
$Min(a)$ and $Min(b)$ is at least $1000\delta$, then 
$ab$ is hyperbolic since both $Min(a), Min(b)$ satisfy the property (M2)).
\qed

This finishes the explanation 
for the part to produce $r'$ from $r$,
and the case (ii) is done. 
We proved Lemma \ref{lemma7.4}.
\qed

We go back to the proof of Proposition 
\ref{2.2}.
With Lemma \ref{lemma7.4}, the rest is very similar to \cite{FS}.
%$n_0$ is given. 
% We choose $n$ large enough, which we will do later. 
Fix $n$ that is large enough to apply Lemma \ref{lemma7.4}.
We explain how we define the map $\phi_n$. 
%Toward defining the map $\phi_n$ in the strategy, for each $m$, 
%we define a map
%$$\Sph_m(G,S_{n_0}) \to B_{m+b}(L,\eta(S)),$$
%where for each $g$ we assign a set of $m-1$ elements. 
%This map depends on $n$.

Let $g \in \Sph_m(G,S_{n_0})$.
Choose a shortest representative word $w(g)$ of length $m$ on $S_{n_0}$ for $g$.
By the bijection $h_{n_0}$ between $S_{n_0}$ and $\eta(S)$, 
$w(g)$ canonically gives  a word $w$ of length $m$ on $\eta(S)$.

From this word $w$, we construct a collection of elements in $L$.
Given a positive integer $k$ with $1 \leq k \leq m-1$, we  
divide the word $w$ into a prefix of length $k$, and a suffix of length $m-k$. The prefix corresponds to an element in $L$ that 
we denote $w^k_p$, and
the suffix corresponds to an element in $L$ that we denote $w^k_s$.

Now, from the four separators we constructed in Lemma \ref{lemma7.4}
we choose a separator $\hat u$ for $h_n$ such that $\hat u$ is 
admissible for $w^k_p$ and $u^{-1}$ is admissible 
for $(w^k_s)^{-1}$, after we map them to $G$ by $h_n$. 
We are writing $\hat u$ instead of $u$ to indicate that 
the separator is in $L$. 
%in $X$, terminates in a germ that is in the orbit of a germ of $y_n$ in $Y$.
%The interval $[y_n,w^k_s(y_n)]$ starts in a germ that is in the orbit of a germ of $y_n$ in $Y$. With the pair $w^k_p,w^k_s$ we associate an element
%$v_{i,j}$, that was constructed in lemma \ref{2.3}, that does not start with the germ that 
%$[y_n,w^k_p(y_n)]$ terminates with, and does not end with the germ that 
%$[y_n,w^k_s(y_n)]$ starts with. 

To the pair $w^k_p,w^k_s$, we associate the following element in $L$:
$w_p^k \hat u  w_s^k \in B_{m+b}(L,\eta(S))$.
$$w \leadsto w_p^k w_s^k \leadsto w_p^k \hat u w_s^k$$
Note that $h_{n_0}(w_p^k \hat u  w_s^k)=h_{n_0}(w)$ for all $k$
since $\hat u \in \ker(h_{n_0})$.

In this way, we obtain $m-1$ ``words'' in $L$ from $w$,
but possibly, some of them represent the same elements in $L$.
 To address this issue, we define a subcollection of words,
 called  {\it forbidden} words 
(in somewhat similar way to what we did in Section \ref{section.4}).

Given $m$, take $n$ large enough such that $h_n$ 
is injective on $B_{2m}(L,\eta(S))$ from now on. 
%\kf{chose n large w.r.t m. check}
\begin{definition}[Forbidden words. cf. Def 3.4 in \cite{FS}]\label{7.4}
Let $w$ be a word of length $m$ on $\eta(S)$ in the above explanation. 
We say that a word $w^k_p \hat u w^k_s$, from the collection that is built from $w$, is {\it forbidden} if there exists $f$, $1 \leq f \leq m$ such that:
$$d_X(h_n(w^k_p \hat u) (y_n),h_n(w^f_p)(y_n)) \leq \frac {1} {5} d_X(y_n, h_n(\hat u)(y_n)).$$
\end{definition}

%With a word $w \in L$ that is associated with a subword of length $m$ of a word in the regular language that the automata that is associated with
%$(\Gamma,S_{\frac{30\delta}{T}})$ produces, we have constructed $m-1$ words of the form $w^k_pv_{i,j}w^k_s$. It is further possible
%to bound the number of the forbidden words of that form.
We give a bound on the number of forbidden words.

\begin{lemma}[cf. Lemma 3.5, \cite{FS}]\label{7.5}
For $m$ and  each word $w$ as in Definition \ref{7.4}, 
%Let $w \in L$ be associated with a subword of length $m$ of a word in the regular language produced by the finite automata that is associated with
%$(\Gamma,S_{\frac{30\delta}{T}})$. Then 
there are at most $\frac {1} {\D+1} m$ forbidden words of the form: $w^k_p \hat u w^k_s$ for $k=1,\ldots,m-1$.
\end{lemma}

In the proof of this lemma we use the assumption that $h_n$
is injective on $B_{2m}(L,\eta(S))$ as we did 
in the proof of Lemma \ref{4.6}.

\proof
The strategy of the proof is same as the proof of Lemma \ref{4.6}.
In there, we ran the argument in $B_m(L,\eta(S))$, but here,
we do it in the set $Z(w)=\{ w^k_p \hat u w^k_s| k=1, \cdots, m-1\}$.
Namely, if $w^k_p \hat u w^k_s\in Z(w)$ is forbidden for some $k$, then
there are two other elements in $Z(w)$ that are candidates for
non-forbidden elements. Then if at least one of them is forbidden,
then there are two other elements in $Z(w)$ that are candidates
for non-forbidden elements, and so on.
We omit details. 
\qed

Then we have:
\begin{lemma}[cf. Lemma 3.6 \cite{FS}]\label{7.6}
%Let $w \in L$ be an element that is associated with a  subword of length $m$ of a word in the regular language that is associated with 
%$(\Gamma,S_{\frac{30\delta}{T}})$. 
For $m$ and the word $w$ as above, 
%Then 
the non-forbidden words: $w^k_p \hat v_{i,j} w^k_s$, for  all
$k$, $1 \leq k \leq m-1$, are  distinct elements in $L$. 
\end{lemma}

\proof
The proof is nearly identical to the proof of Lemma 3.6 in \cite{FS}, 
and we omit it.
\qed

As in Section \ref{section.4}, we define adequate elements in $Z(w)$.
We choose a maximal subset in the set of non-forbidden
elements in $Z(w)$ such that for any two distinct
elements $z_1, z_2$ in the subset, $h_n(z_1), h_n(z_2)$
are not in the same coset w.r.t. $F(h_n(u_i))$
for any separator $u_i$.
We call those elements {\it adequate elements}.
This notion depends on $n$.
Then, as before,
$$\frac{|\{\text{non-forbidden elements}\} |}{\D^4} \le |\{\text{adequate elements}\}|.$$

This is because, as we explained,  there are only four separators, $\hat u$,
and $|F(h_n(\hat u))| \le \D$.
In conclusion, since $|Z(w)|=m-1$, 
\begin{equation}\label{ad2}
\frac{m-1}{\D^4(\D+1)} \le |\{\text{adequate elements in } Z(w)\}|.
\end{equation}

Using  adequate elements, we construct a collection of {\it feasible} words in $L$. 

\begin{definition}[Feasible words in $L$. cf. Definition 3.7 \cite{FS}]\label{7.7}
Let $m, q$ be positive integers. Let $w$ be
a word of length $mq$ on $\eta(S)$ that is associated to 
an element $g\in \Sph_{mq}(G,S_{n_0})$ as in the above discussion. 
We present $w$ as a concatenation of $q$ subwords of length $m$: $w=w(1) \ldots w(q)$.

Then for any choice of integers: $k_1,\ldots,k_q$ with $1 \leq k_t \leq m-1$,
and  $t=1,\ldots,q$, for which all the elements, $w(t)^{k_t}_p v^t w(t)^{k_t}_s$, 
are adequate (here we drop the ``hat'' from $v^t$ deliberately
although they are in $L$ to avoid confusion since we want to use it right in the below) , we associate a {\it feasible} word (of type $q$)  on $\eta(S)$
(in $L$):
$$w(1)^{k_1}_p v^1 w(1)^{k_1}_s \, \hat v^1 \,   
w(2)^{k_2}_pv^2 w(2)^{k_2}_s \, \hat v^2  \, \ldots 
w(q)^{k_q}_pv^q w(q)^{k_q}_s,$$ 
where for each $t$, $1 \leq t \leq q-1$, $\hat v^t $ is one of the separators 
from Lemma \ref{4.2} such that 
$\hat v^t$ is admissible for $w(t)^{k_t}_s$
and $(\hat v^t)^{-1}$ is admissible for $\left(w(t+1)^{k_{t+1}}_p\right)^{-1}$.
\end{definition}

Finally we define the map $\phi_n$. 
Suppose positive integers $m,q$ are given.
Then for $g \in \Sph_{mq}(G, S_{n_0})$, choose one shortest
representative $w(g)$ of length $mq$ on $S_{n_0}$, which defines a word $\tilde w(g)$
of length $mq$ on $\eta(S)$ as in the definition \ref{7.7}. From $\tilde w(g)$ we produce 
feasible words on $\eta(S)$, which define {\it feasible} elements in $L$.
Note that those elements are in $B_{q(m+2b)}(L,\eta(S))$
and mapped to $g$ by $h_{n_0}$.
We denote this collection as $\phi_n(g)$.
%This is the map $\phi$ that appeared in the strategy.

%Note that for distinct $g_1,g_2 \in \Sph_{mq}(G, S_{n_0})$,
%the sets $\phi(g_1), \phi(g_2)$ are disjoint since 
%each of them mapped to $g_1, g_2$, respectively, by $h_{n_0}$.

%We collect those feasible elements obtained from all $g \in \Sph_{mq}(G, S_{n_0})$,
%and denote the set by $\phi(g) \subset B_{q(m+2b)}(L,\eta(S))$.
 %\kf{explain the lower bound of $\phi_n(g)$}

We have :

\begin{lemma}[cf. Lemma 3.8 \cite{FS}]\label{7.8}
For any positive integers $m,q$, the feasible elements
in the collection $\phi_n(g)$
we obtain for each $g \in \Sph_{mq}(G, S_{n_0})$
are all distinct in $L$.

Moreover, all the feasible elements 
obtained from all the elements $g$ in $\Sph_{mq}(G, S_{n_0})$
are all distinct in $L$.
\end{lemma}

In this lemma, $n$ must be large enough
in the sense that Lemma  \ref{lemma7.4} applies and also, 
for the given $m$, the map $h_n$ is injective
on $B_{2m}(L,\eta(S))$ (cf. Lemma \ref{lemma.tail}). 
We will choose such $n$ in the proof. 

\proof
The proof of the first sentence  is similar to the proof of Lemma \ref{feasible}.
%\kf{check the lemma number}
 So we omit it.
(cf. the proof of 
Lemma 3.8 \cite{FS}.)
Then the moreover part immediately follows. 
\qed

We note that for each $g$, we have
\begin{equation}\label{ad3}
\left(\frac{m-1}{\D^4(\D-1)}\right)^q \le \phi_n(g).
\end{equation}
This is from the lower bound (\ref{ad2}) on the number of 
adequate elements and the lemma.

We finish the proof of the proposition \ref {2.2}.
We review the setting.
The constant $n_0$ is given
to start with. We want to show $e(G,f_{n_0}(S)) < e(L,\eta(S))$.
%We choose and fix $1\not=r_{n_0} \in \ker(h_{n_0})<L$.
Then the constant $b$ is given by Lemma \ref{lemma7.4}, which does not depend on $n$.
Choose $n$ large enough so that we can apply Lemma \ref{lemma7.4}.

Then choose $m$ such that 
$$\log (m-1)> 2b \log (e(L,\eta(S)) + \log (\D^4(\D-1)).$$
This implies:
$$\log (m-1)> 2b \log (e(G,S_{n_0})) + \log (\D^4(\D-1)).$$

Then choose $n>0$ larger if necessary   such that 
$h_n$ is injective on $B_{2m}(L,\eta(S))$.
We need this to apply Lemma \ref{7.6}.

Now, for all $q>0$, 
by combining the lower bound (\ref{ad3}) and (the moreover part of)
Lemma \ref{7.8}, 
we have 
$$|\Sph_{mq}(G,S_{n_0})| \left(\frac{m-1}{\D^4(\D-1)}\right)^q \le |B_{q(m+2b)}(L,\eta(S))|.$$
From this,
\begin{align*}
\log e(L,\eta(S))  & \ge \lim_{q \to \infty} \frac {\log ( |\Sph_{mq}(G,S_{n_0})| (\frac {m-1} {\D^4(\D-1)})^q)} {q(m+2b)} 
\\
& = \lim_{q \to \infty} \frac {\log (|\Sph_{mq}(G,S_{n_0})|)} {q(m+2b)}  + \frac { q(\log (m-1) - \log(\D^4(\D-1)))} {q(m+2b)}
\\
&=
\log(e(G,S_{n_0}))  \frac {m} {m+2b} + \frac {\log(m-1) - \log (\D^4(\D-1))} {m+2b}
\\
&> \log (e(G,S_{n_0})).
\end{align*}
The last inequality is by the way we chose $m$. 
Hence, $e(L,\eta(S)) > e(G,S_{n_0})$.
We proved Proposition \ref{2.2}.
\qed

\subsection{Family version}
%\kf{check}
We state a family version of Proposition \ref{2.2}.
Let $\delta, M$ be constants and $D(\epsilon)$
a function for WPD. 
Suppose $X_n$ is $\delta$-hyperbolic, 
a group $G_n$ acts on $X_n$ and $S_n$ is a finite
generating set of $G_n$ such that $S_n^M$ contains a hyperbolic 
element on $X_n$ that is $D$-WPD.

Suppose that $|S_n|=\ell$ for all $n$, and let $F$
be a free group of rank $\ell$ with a free
generating set $S$. Let $f_n:(F,S)\to (G_n,S_n)$ be a surjection  with a bijection $f_n(S)=S_n$. 
Assume that the sequence $\{f_n\}$ converges to $\eta:(F,S)
\to (L,\eta(S))$.

Also,  assume that for all $n$, there is 
a surjection  $h_n:(L,\eta(S)) \to (G_n,S_n)$
with $f_n=h_n \circ \eta$.
As in Proposition \ref{semi.cont.family}, we put
this as an assumption. 
Then we have the following generalizing
Proposition \ref{2.2}. The proof is identical
and we omit it. 
\begin{prop}\label{family.2.2}
Assume that for all $n$, there is 
a surjection  $h_n:(L,\eta(S)) \to (G_n,S_n)$
with $f_n=h_n \circ \eta$.
If  $\ker(h_{n_0})$ is infinite for some $n_0$, then 
$$e(G_{n_0},f_{n_0}(S))<e(L,\eta(S)).$$ 
\end{prop}

\subsection{Finiteness for $\Theta_X(G)$}
A finiteness result similar to  Theorem \ref{7.1}
holds for subgroups. It is known for hyperbolic
groups, \cite[Theorem 5.3]{FS}.

Let $S_1, S_2 \subset G$ be two finite subsets.
Let $H_i=\<S_i\> < G$ be the subgroup generated by $S_i$.
We say $(H_1,S_1)$ and $(H_2,S_2)$ are {\it isomorphic}
if there is a bijection between $S_1, S_2$ that induces
an isomorphism between $H_1, H_2$.

\begin{thm}[Finiteness in the subgroups case]\label{7.10}
Assume the same condition on $G$ as in Theorem \ref{7.1}.
Moreover, we assume that if $S$ is a finite set of $G$ such that 
$\< S \>$ contains a hyperbolic isometry on $X$, then
$S^M$ contains a hyperbolic isometry that is $D$-WPD.
Let $\rho \in \Theta_X(G)$, 
Then there are at  most finitely many $(H,S)$, up to isomorphism,
such  that $S\subset G$ is finite, $H=\<S\>$, $H$ contains
a hyperbolic element on $X$, and $\rho=e(H,S)$.
\end{thm}

The proof is nearly identical to the proof of Theorem \ref{7.1}
and we only need to modify the setting from the entire group $G$
to subgroups. 
\proof
To argue by contradiction, let $\rho \in \Theta_X(G)$ and suppose that  there are infinitely many distinct, up to 
isomorphism, $(H_n, S_n)$ with $e(H_n,S_n)=\rho$ such that $H_n$
contains a hyperbolic isometry on $X$.
Note that by assumption, $S_n^M$ contains a hyperbolic
isometry that is $D$-WPD.

As before, by Proposition \ref{bound.generators}, passing to a subsequence, one may assume 
that there is $\ell$ with $|S_n|=\ell$ for all $n$. 
Then we obtain $f_n:F \to G$ with $f_n(S)=S_n$,
where $F$ is the free group on $S$ with $|S|=\ell$.
Then, passing to a subsequence again, $f_n$ converges
to a limit group $(L, \eta(S))$ with $\eta:F \to L$.
%A finitely generated subgroup of an equationally Noetherian group
%is equationally Noetherian (immediate from the definition),
%therefore all $H_n$ are equationally Noetherian. 
Since $G$ is equationally Noetherian, by Lemma \ref{basic}, passing to a further subsequence, we may assume that there are $h_n:L \to H_n <G$ with $h_n \circ \eta=f_n$ for all $n$.

First, by Proposition \ref{6.3}, we have that
$e(H_n,S_n)=e(L,\eta(L))$ for all $n$.

On the other hand, 
we prove a version of  Proposition \ref{2.2} for subgroups:
if $\ker(h_{n_0})$ is infinite for some $n_0$,
then $e(H_{n_0},S_{n_0}) < e(L,\eta(S))$.
The proof is same and we only outline it.
As before set $\D=D(100\delta)$.
First, since $H_{n}$ contains a hyperbolic 
isometry on $X$ that is $D$-WPD, 
each $H_n$ contains the maximal finite normal subgroup, which 
we denote by $N_{H_n}$ with $|N_{H_n}| \le 2\D$,
(by Lemma \ref{radical.limit.group}).

The key step is to prove a lemma similar to Lemma \ref{lemma7.4}.
%\kf{lemma number}
The argument is same. We use that $S_n^M$ has
a hyperbolic and $D$-WPD element for all $n$. 
Then, 
as before, for all sufficiently large $n$, there exists an element $r(n) \in \ker(h_{n_0})$ such that 
$h_n(r(n))$ is hyperbolic on $X$, and that the word
length of $r(n)$ in terms of $\eta(S)$ is bounded uniformly on $n$.
Then then rest is same as proving the lemma. 
Once we have the lemma, the rest is same to 
show the proposition. 
(We point out that this is a special case of 
Proposition \ref{family.2.2}, where $X_n$ are common. 
But we did not describe the details of the argument for that.)

Combining those two, we conclude that 
$\ker(h_n)$ is finite for all $n$. 

%Then we proceed as before. Namely, we prove that $L$
%contains the maximal finite normal subgroup $N_L$ with $|N_L|\le D$.
%This is a counterpart to Lemma \ref{radical.limit.group}.
%The argument is same since we have $|N_{H_n}| \le D$. 

Finally, there are only finitely many possibilities
for $\ker(h_n)$ since it is contained in $N_L$
and $|N_L|\le 2\D$. 
It implies that the desired finiteness for $(H_n, S_n)$ holds as before. 
\qed

\subsection{Examples}
We give some examples of Theorem \ref{7.1} and Theorem \ref{7.10}.
%\kf{subgorups too}
We start with relatively hyperbolic groups.

\begin{thm}[Finiteness for relatively hyperbolic groups]\label{7.11}
Let $G$ be a group that is hyperbolic relative to a collection of 
subgroups $\{P_1, \cdots, 
P_n\}$.
Suppose $G$ is not virtually cyclic, and not equal to $P_i$
for any $i$.
Suppose each $P_i$ is finitely generated and equationally Noetherian.
Then for each $\rho \in \xi(G)$ there
are at most finitely many finite generating
sets $S_n$ of $G$, up to $\Aut(G)$, s.t. $e(G,S_n)=\rho$.

Moreover, for each $\rho \in \Theta_{{\rm non-elem.}}(G)$, there are
at most finitely many $(H_n,S_n)$, up to isomorphism,
s.t. $e(H_n,S_n)=\rho$, where $S_n \subset G$ is finite and 
$H_n=\<S_n\>$ is not conjugate into any $P_i$.
\end{thm}

\proof
Let $X$ be a hyperbolic space on which $G$ acts as we explained in Section \ref{section.5.2}. We also verified that all the assumption of Theorem 
\ref{7.1} for $G$ and the action of $G$ on $X$. 
Recall that the action of $G$ on $X$ is uniformly WPD by Lemma \ref{rel.hyp.wpd}. 
It implies the first part of the theorem.

For the moreover part, we apply Theorem \ref{7.10}.
As Lemma \ref{subgroup.rel.hyp} shows, $\Theta_{{\rm non-elem.}}(G)
=\Theta_X(G)$, which implies the conclusion.
\qed

Theorem \ref{7.11} immediately implies the following as
Theorem \ref{thm.subgroup.rel.hyp} implies Theorem \ref{6.1}:
\begin{thm}[Finiteness for lattices]\label{lattice.finite}
Let $G$ be a group in Theorem \ref{lattice}.
Then for each $\rho \in \xi(G)$ there
are at most finitely many finite generating
sets $S_n$, up to $\Aut(G)$, s.t. $e(G,S_n)=\rho$.

Moreover, for each $\rho \in \Theta(G)$, there are
at most finitely many $(H_n,S_n)$, up to isomorphism,
s.t. $e(H_n,S_n)=\rho$, where $S_n \subset G$ is finite and 
$H_n=\<S_n\>$.
\end{thm}

%\proof
%As in the proof of Proposition \ref{prop.rank1}, 
%let $X=Cayley(G,S_0 \bigcup \cup_iH_i)$.
%Then since $X/G$ has only one vertex, the assumption 
%of Theorem \ref{7.1} is satisfied.
%
%For the moreover part, the conclusion holds
%for $\Theta_X(G)$ by Theorem \ref{7.10}.
%But, as we said in the proof of Theorem
%\ref{6.1}, $\Theta_X(G)=\Theta(G)$, which implies the 
%conclusion. 
%\qed

Lastly we record the following (potential) example:

\begin{thm}[Finiteness for MCG]
Let $MCG=MCG(\Sigma)$ be the mapping class group
of a compact orientable surface $\Sigma$. Assume that 
it is equationally Noetherian. 

Then for each $\rho \in \xi(MCG)$ there
are at most finitely many finite generating
sets $S_n$, up to $\Aut(MCG)$, such that $e(MCG,S_n)=\rho$.

Moreover, for each $\rho \in \Theta_{{\rm large}}(MCG)$, 
there are
at most finitely many $(H_n,S_n)$, up to isomorphism,
such that $e(H_n,S_n)=\rho$, where $S_n \subset MCG$ is finite and 
$H_n=\<S_n\>$ is a large subgroup.
\end{thm}

\proof
As we explained in Section \ref{section.mcg}, the action of $MCG(\Sigma)$ on the curve graph $X=\CC(\Sigma)$
satisfies the assumption of Theorem \ref{7.1}.
It is uniformly WPD. 

For the moreover part, the conclusion holds
for $\Theta_X(MCG)$ by Theorem \ref{7.10}.
But, as we said in the proof of Theorem \ref{thm.6.5}, we have $\Theta_X(MCG)=\Theta_{{\rm large}}(MCG)$, 
so that the conclusion holds. 
\qed


\begin{thebibliography}{BMR2}

%\bibitem{Arzhantseva-Lysenok1}
%G. N. Arzhantseva and I. G. Lysenok, Growth tightness for word hyperbolic groups, Math. Zeitschrift 241 (2002), 597-611.

\bibitem[AL]{AL}
G. N. Arzhantseva and I. G. Lysenok, A lower bound on the growth of word hyperbolic groups, Jour. of the LMS 73 (2006), 109-125.



\bibitem[BMR]{BMR}
Gilbert Baumslag,  Alexei Myasnikov,   Vladimir Remeslennikov, 
Algebraic geometry over groups. I. Algebraic sets and ideal theory.
J. Algebra 219 (1999), no. 1, 16–-79. 

\bibitem[BMR2]{BMR2}
Gilbert Baumslaga, Alexei Myasnikov, Vitaly Roman'kov.
Two Theorems about Equationally Noetherian Groups.
Journal of Algebra
Volume 194, Issue 2, Pages 654-664





\bibitem[BeF]{BeF}

Mladen Bestvina, Koji  Fujiwara, Bounded cohomology of subgroups of mapping class groups. Geom. Topol. 6 (2002), 69-–89. 

\bibitem[BeF2]{BeF2}
Mladen Bestvina, Koji Fujiwara.
Handlebody subgroups in a mapping class group.
{\it In the Tradition of Ahlfors–Bers, VII}. 
Contemporary Mathematics Volume 696, 2017, Pages 30-50.
AMS.

\bibitem[Bo1]{Bo1}
B. H. Bowditch, 
Relatively hyperbolic groups. 
Internat. J. Algebra Comput. 22 (2012), no. 3, 66 pp.

\bibitem[Bo]{Bowditch}
Brian H. Bowditch, 
Tight geodesics in the curve complex. 
Invent. Math. 171 (2008), no. 2, 281-–300.


\bibitem[BrF]{BrF}
E. Breuillard, K. Fujiwara.
On the joint spectral radius for isometries of non-positively curved spaces and uniform growth.
to appear in Annales de l'institut Fourier.



\bibitem[Br]{Br}
Roger M. Bryant, 
The verbal topology of a group.
J. Algebra 48 (1977), no. 2, 340–-346. 




\bibitem[DGO]{DGO}
F. Dahmani, V. Guirardel, D. Osin, 
Hyperbolically embedded subgroups and rotating families in groups acting on hyperbolic spaces.
Mem. Amer. Math. Soc. 245 (2017), no. 1156, 152 pp.

\bibitem[De]{De}

 Thomas Delzant, Kaehler groups, $\Bbb R$-trees, and holomorphic families of Riemann surfaces. Geom. Funct. Anal. 26 (2016), no. 1, 160-–187. 



\bibitem[DrS]{DrS}
Cornelia Dru\c{ţ}u, Mark Sapir, Tree-graded spaces and asymptotic cones of groups. With an appendix by Denis Osin and Mark Sapir. Topology 44 (2005), no. 5, 959-–1058.

\bibitem[EMO]{EMO}
Alex Eskin, Shahar Mozes, Hee Oh, 
On uniform exponential growth for linear groups.
Invent. Math. 160 (2005), no. 1, 1-–30. 

\bibitem[Fa]{Fa}
B.  Farb, Relatively hyperbolic groups. Geom. Funct. Anal. 8 (1998), no. 5, 810-–840. 

\bibitem[F]{F}
K. Fujiwara, Subgroups generated by two pseudo-Anosov elements in a mapping class group. I. Uniform exponential growth. {\it Groups of diffeomorphisms}, 283–-296, Adv. Stud. Pure Math., 52, Math. Soc. Japan, 2008.



\bibitem[FS]{FS}
K. Fujiwara, Z. Sela.
The rates of growth in a hyperbolic group.
to appear in Invent. Math.

\bibitem[G]{G}
M. Gromov. 
Hyperbolic groups. Essays in group theory, 75–263, Math. Sci. Res. Inst. Publ., 8, Springer, New York, 1987. 

\bibitem[GHL]{GHL}
Daniel Groves, Michael Hull, Hao Liang.
Homomorphisms to 3-manifold groups.
arXiv:2103.06095






\bibitem[GrH]{GrH}
 D. Groves, M. Hull,  Homomorphisms to acylindrically hyperbolic groups I: Equationally noetherian groups and families. Trans. Amer. Math. Soc. 372 (2019), no. 10, 7141-–7190.

\bibitem[Gu]{Guba}
Victor S. Guba. Equivalence of infinite systems of equations in free groups and semigroups to finite subsystems. Mat. Zametki, 48(3):321-–324, 1986. 



\bibitem[K]{Koubi}
M. Koubi,  Croissance uniforme dans les groupes hyperboliques, Ann. Inst. Fourier 48 (1998), 1441--1453.




\bibitem[Ma]{M}
Johanna Mangahas,  A recipe for short-word pseudo-Anosovs. Amer. J. Math. 135 (2013), no. 4, 1087–-1116.

\bibitem[Ma2]{M2}
Johanna Mangahas, 
Uniform uniform exponential growth of subgroups of the mapping class group. 
Geom. Funct. Anal. 19 (2010), no. 5, 1468-–1480.

\bibitem[MM1]{MM1}
 Howard A. Masur, Yair N. Minsky,  Geometry of the complex of curves. I. Hyperbolicity. Invent. Math. 138 (1999), no. 1, 103-–149.




\bibitem[MO]{MO} 
Ashot Minasyan, Denis Osin.
Acylindrical hyperbolicity of groups acting on trees.
Mathematische Annalen 362,  (2015), 1055–1105

\bibitem[NS]{NS}
Hoang Thanh Nguyen, Hongbin Sun, 
Subgroup distortion of 3-manifold groups. 
Trans. Amer. Math. Soc. 373 (2020), no. 9, 6683–-6711. 

\bibitem[O]{O}
D. Osin.
 Acylindrically hyperbolic groups. Trans. Amer. Math. Soc. 368 (2016), no. 2, 851–-888. 


\bibitem[RW]{RW}
C. Reinfeldt, R. Weidmann, 
Makanin-Razborov diagrams for hyperbolic groups.
Ann. Math. Blaise Pascal 26 (2019), no. 2, 119-–208. 




\bibitem[Sa]{Sa}
A. Sambusetti. Growth tightness of free and amalgamated products. Ann. Sci. Ecole Norm. Sup. serie 35 (2002), no. 4, 477–-488.


\bibitem[S]{Sela1}
Z. Sela,  Diophantine geometry over groups I: Makanin-Razborov diagrams, Publication de la IHES 93 (2001), 31-105.

%\bibitem{Sela2}
%Z. Sela,  Word equations I: Pairs and their Makanin-Razborov diagrams, preprint, arXiv:1607.05431. 
%
%\bibitem{Sela3}
%Z. Sela,  Non-commutative algebraic geometry I: Monomial equations with a single variable, preprint, arXiv:1906.02049. 

%\bibitem{Sela4}
%Z. Sela,  Endomorphisms of hyperbolic groups  I: The Hopf property, Topology 38 (1999), 301-321.

\bibitem[Se]{Se}
Jean-Pierre Serre, Trees. Translated from the French by John Stillwell. Springer-Verlag, Berlin-New York, 1980.

\bibitem[Si]{Si} A. Sisto.
Projections and relative hyperbolicity.
Enseign. Math. (2) 59 (2013), no. 1-2, 165-–181.


\bibitem[W]{Wilson}
J. S. Wilson,  On exponential growth and uniformly exponential growth for groups,
Inventiones math. 155 (2004), 287-303.

\bibitem[WZ]{WZ}
Henry Wilton, Pavel Zalesskii, Profinite detection of 3-manifold decompositions. Compos. Math. 155 (2019), no. 2, 246–259. 

\bibitem[X]{X}
Xiangdong Xie, 
Growth of relatively hyperbolic groups. 
Proc. Amer. Math. Soc. 135 (2007), no. 3, 695–-704.

\end{thebibliography}
\end{document}